\documentclass{article}
\usepackage{amsmath}
\usepackage{amssymb}
\usepackage{amsthm}
\usepackage{cite}
\usepackage[arrow,curve,matrix]{xy}
\begin{document}
%%%%%%%%%%%%%%%%%%%%%%%%%%%%%%%%%%%%%%%%%%%%%%%%%%%%%%%%%%%%%%%%%%%%%%%
%%%%%%%%%%%%%%%%%%%%%%%%%     Macros      %%%%%%%%%%%%%%%%%%%%%%%%%%%%%
%%%%%%%%%%%%%%%%%%%%%%%%%%%%%%%%%%%%%%%%%%%%%%%%%%%%%%%%%%%%%%%%%%%%%%%
\def\e#1\e{\begin{equation}#1\end{equation}}
\def\ea#1\ea{\begin{align}#1\end{align}}
\def\eq#1{{\rm(\ref{#1})}}
\theoremstyle{plain}% default
\newtheorem{thm}{Theorem}[section]
\newtheorem{lem}[thm]{Lemma}
\newtheorem{prop}[thm]{Proposition}
\newtheorem{cor}[thm]{Corollary}
\theoremstyle{definition}
\newtheorem{dfn}[thm]{Definition}
\newtheorem{ass}[thm]{Assumption}
\newtheorem{quest}[thm]{Question}
\newtheorem{rem}[thm]{Remark}
\newtheorem{ex}[thm]{Example}
\def\na{{\rm na}}
\def\stk{{\rm stk}}
\def\ind{{\rm \kern.05em ind}}
\def\fin{{\rm fin}}
\def\Ho{{\rm Ho}}
\def\Po{{\rm Po}}
\def\uni{{\rm uni}}
\def\rk{{\rm rk}}
\def\vi{{\rm vi}}
\def\opp{{\rm op}}
\def\dim{\mathop{\rm dim}\nolimits}
\def\cha{\mathop{\rm char}}
\def\Im{\mathop{\rm Im}}
\def\GL{\mathop{\rm GL}\nolimits}
\def\Spec{\mathop{\rm Spec}\nolimits}
\def\Stab{\mathop{\rm Stab}\nolimits}
\def\Sch{\mathop{\rm Sch}\nolimits}
\def\coh{\mathop{\rm coh}}
\def\Hom{\mathop{\rm Hom}\nolimits}
\def\Iso{\mathop{\rm Iso}\nolimits}
\def\Aut{\mathop{\rm Aut}}
\def\End{\mathop{\rm End}}
\def\ad{\mathop{\rm ad}}
\def\Mor{\mathop{\rm Mor}\nolimits}
\def\CF{\mathop{\rm CF}\nolimits}
\def\CFi{\mathop{\rm CF}\nolimits^{\rm ind}}
\def\LCF{\mathop{\rm LCF}\nolimits}
\def\dLCF{{\dot{\rm LCF}}\kern-.1em\mathop{}\nolimits}
\def\SF{\mathop{\rm SF}\nolimits}
\def\SFa{\mathop{\rm SF}\nolimits_{\rm al}}
\def\SFai{\mathop{\rm SF}\nolimits_{\rm al}^{\rm ind}}
\def\uSF{\mathop{\smash{\underline{\rm SF\!}\,}}\nolimits}
\def\uSFi{\mathop{\smash{\underline{\rm SF\!}\,}}\nolimits^{\rm ind}}
\def\oSF{\mathop{\bar{\rm SF}}\nolimits}
\def\oSFa{\mathop{\bar{\rm SF}}\nolimits_{\rm al}}
\def\oSFai{{\ts\bar{\rm SF}{}_{\rm al}^{\rm ind}}}
\def\uoSF{\mathop{\bar{\underline{\rm SF\!}\,}}\nolimits}
\def\uoSFa{\mathop{\bar{\underline{\rm SF\!}\,}}\nolimits_{\rm al}}
\def\uoSFi{\mathop{\bar{\underline{\rm SF\!}\,}}\nolimits_{\rm al}^{\rm
ind}}
\def\LSF{\mathop{\rm LSF}\nolimits}
\def\LSFa{\mathop{\rm LSF}\nolimits_{\rm al}}
\def\uLSF{\mathop{\smash{\underline{\rm LSF\!}\,}}\nolimits}
\def\dLSF{{\dot{\rm LSF}}\kern-.1em\mathop{}\nolimits}
\def\doLSF{{\dot{\bar{\rm LSF}}}\kern-.1em\mathop{}\nolimits}
\def\duoLSF{{\dot{\bar{\underline{\rm LSF\!}\,}}}\kern-.1em\mathop{}
\nolimits}
\def\ouLSF{{\bar{\underline{\rm LSF\!}\,}}\kern-.1em\mathop{}
\nolimits}
\def\dLSFi{{\dot{\rm LSF}}\kern-.1em\mathop{}\nolimits^{\rm ind}}
\def\dLSFa{{\dot{\rm LSF}}\kern-.1em\mathop{}\nolimits_{\rm al}}
\def\doLSFa{{\dot{\bar{\rm LSF}}}\kern-.1em\mathop{}\nolimits_{\rm al}}
\def\dLSFai{{\dot{\rm LSF}}\kern-.1em\mathop{}\nolimits^{\rm ind}_{\rm
al}}
\def\duLSF{{\dot{\underline{\rm LSF\!}\,}}\kern-.1em\mathop{}\nolimits}
\def\duLSFi{{\dot{\underline{\rm LSF\!}\,}}\kern-.1em\mathop{}
\nolimits^{\rm ind}}
\def\oLSF{\mathop{\bar{\rm LSF}}\nolimits}
\def\oLSFa{\mathop{\bar{\rm LSF}}\nolimits_{\rm al}}
\def\oLSFai{\mathop{\bar{\rm LSF}}\nolimits_{\rm al}^{\rm ind}}
\def\uoLSF{\mathop{\bar{\underline{\rm LSF\!}\,}}\nolimits}
\def\uoLSFa{\mathop{\bar{\underline{\rm LSF\!}\,}}\nolimits_{\rm al}}
\def\uoLSFi{\mathop{\bar{\underline{\rm LSF\!}\,}}\nolimits_{\rm
al}^{\rm ind}}
\def\Ext{\mathop{\rm Ext}\nolimits}
\def\id{\mathop{\rm id}\nolimits}
\def\ha{{\ts\frac{1}{2}}}
\def\Obj{\mathop{\rm Obj\kern .1em}\nolimits}
\def\fObj{\mathop{\mathfrak{Obj}\kern .05em}\nolimits}
\def\modA{\text{\rm mod-$A$}}
\def\modKQ{\text{\rm mod-$\K Q$}}
\def\modKQI{\text{\rm mod-$\K Q/I$}}
\def\nilKQ{\text{\rm nil-$\K Q$}}
\def\nilKQI{\text{\rm nil-$\K Q/I$}}
\def\bs{\boldsymbol}
\def\ge{\geqslant}
\def\le{\mathord{\leqslant}}
\def\pr{{\mathop{\preceq}\nolimits}}
\def\npr{{\mathop{\npreceq}\nolimits}}
\def\tl{\trianglelefteq\nobreak}
\def\ps{\precsim\nobreak}
\def\ls{{\mathop{\lesssim\kern .05em}\nolimits}}
\def\bF{{\mathbin{\mathbb F}}}
\def\N{{\mathbin{\mathbb N}}}
\def\Z{{\mathbin{\mathbb Z}}}
\def\Q{{\mathbin{\mathbb Q}}}
\def\C{{\mathbin{\mathbb C}}}
\def\K{{\mathbin{\mathbb K\kern .05em}}}
\def\KP{{\mathbin{\mathbb{KP}}}}
\def\cC{{\mathbin{\mathcal C}}}
\def\A{{\mathbin{\mathcal A}}}
\def\B{{\mathbin{\mathcal B}}}
\def\F{{\mathbin{\mathcal F}}}
\def\G{{\mathbin{\mathcal G}}}
\def\H{{\mathbin{\mathcal H}}}
\def\L{{\mathbin{\mathcal L}}}
\def\M{{\mathcal M}}
\def\cP{{\mathbin{\mathcal P}}}
\def\cQ{{\mathbin{\mathcal Q}}}
\def\cR{{\mathbin{\mathcal R}}}
\def\g{{\mathbin{\mathfrak g}}}
\def\h{{\mathbin{\mathfrak h}}}
\def\n{{\mathbin{\mathfrak n}}}
\def\fD{{\mathbin{\mathfrak D}}}
\def\fE{{\mathbin{\mathfrak E}}}
\def\fF{{\mathbin{\mathfrak F}}}
\def\fG{{\mathbin{\mathfrak G}}}
\def\fH{{\mathbin{\mathfrak H}}}
\def\fM{{\mathbin{\mathfrak M}}}
\def\fR{{\mathbin{\mathfrak R}}}
\def\fS{{\mathbin{\mathfrak S}}}
\def\fT{{\mathbin{\mathfrak T}}}
\def\fU{{\mathbin{\mathfrak U\kern .05em}}}
\def\fV{{\mathbin{\mathfrak V}}}
\def\sIp{{\smash{\sst(I,\pr)}}}
\def\sIt{{\smash{\sst(I,\tl)}}}
\def\sJp{{\smash{\sst(J,\pr)}}}
\def\sKt{{\smash{\sst(K,\tl)}}}
\def\sIl{{\smash{\sst(I,\ls)}}}
\def\sJl{{\smash{\sst(J,\ls)}}}
\def\al{\alpha}
\def\be{\beta}
\def\ga{\gamma}
\def\de{\delta}
\def\bde{\bar\delta}
\def\io{\iota}
\def\ep{\epsilon}
\def\la{\lambda}
\def\ka{\kappa}
\def\th{\theta}
\def\ze{\zeta}
\def\si{\sigma}
\def\om{\omega}
\def\De{\Delta}
\def\La{\Lambda}
\def\Om{\Omega}
\def\Ga{\Gamma}
\def\Si{\Sigma}
\def\Th{\Theta}
\def\Up{\Upsilon}
\def\ts{\textstyle}
\def\sst{\scriptscriptstyle}
\def\sm{\setminus}
\def\bu{\bullet}
\def\op{\oplus}
\def\ot{\otimes}
\def\bigop{\bigoplus}
\def\bigot{\bigotimes}
\def\iy{\infty}
\def\ra{\rightarrow}
\def\ab{\allowbreak}
\def\longra{\longrightarrow}
\def\dashra{\dashrightarrow}
\def\hookra{\hookrightarrow}
\def\lt{\ltimes}
\def\el{{\mathbin{\ell\kern .08em}}}
\def\t{\times}
\def\ci{\circ}
\def\ti{\tilde}
\def\md#1{\vert #1 \vert}
%%%%%%%%%%%%%%%%%%%%%%%%%%%%%%%%%%%%%%%%%%%%%%%%%%%%%%%%%%%%%%%%%%%%%%%
%%%%%%%%%%%%%%%%%%%%%     Text of paper    %%%%%%%%%%%%%%%%%%%%%%%%%%%%
%%%%%%%%%%%%%%%%%%%%%%%%%%%%%%%%%%%%%%%%%%%%%%%%%%%%%%%%%%%%%%%%%%%%%%%
\title{Configurations in abelian categories. II. \\
Ringel--Hall algebras}
\author{Dominic Joyce}
\date{}
\maketitle

\begin{abstract}
This is the second in a series on {\it configurations\/} in an
abelian category $\A$. Given a finite poset $(I,\pr)$, an
$(I,\pr)$-{\it configuration\/} $(\si,\io,\pi)$ is a finite
collection of objects $\si(J)$ and morphisms $\io(J,K)$ or
$\pi(J,K):\si(J)\ra\si(K)$ in $\A$ satisfying some axioms, where
$J,K\subseteq I$. Configurations describe how an object $X$
in $\A$ decomposes into subobjects.

The first paper defined configurations and studied moduli spaces
of $(I,\pr)$-configurations in $\A$, using the theory of Artin
stacks. It showed well-behaved moduli stacks $\fObj_\A,
\fM(I,\pr)_\A$ of objects and configurations in $\A$ exist
when $\A$ is the abelian category $\coh(P)$ of coherent
sheaves on a projective scheme $P$, or $\modKQ$ of
representations of a quiver~$Q$.

Write $\CF(\fObj_\A)$ for the vector space of $\Q$-valued
constructible functions on the stack $\fObj_\A$. Motivated by the
idea of {\it Ringel--Hall algebras}, we define an associative
multiplication $*$ on $\CF(\fObj_\A)$ using pushforwards and
pullbacks along 1-morphisms between configuration moduli stacks, so
that $\CF(\fObj_\A)$ is a $\Q$-{\it algebra}. We also study
representations of $\CF(\fObj_\A)$, the Lie subalgebra
$\CFi(\fObj_\A)$ of functions supported on indecomposables, and
other algebraic structures on~$\CF(\fObj_\A)$.

Then we generalize all these ideas to {\it stack functions}
$\uSF(\fObj_\A)$, a universal generalization of constructible
functions, containing more information. When $\Ext^i(X,Y)=0$ for all
$X,Y\in\A$ and $i>1$, or when $\A=\coh(P)$ for $P$ a Calabi--Yau
3-fold, we construct ({\it Lie}) {\it algebra morphisms} from stack
algebras to explicit algebras, which will be important in the
sequels on invariants counting $\tau$-semistable objects in~$\A$.
\end{abstract}

\section{Introduction}
\label{ai1}

This is the second in a series of papers \cite{Joyc3,Joyc4,Joyc5}
developing the concept of {\it configuration} in an abelian category
$\A$. Given a finite partially ordered set (poset) $(I,\pr)$, we
define an $(I,\pr)$-{\it configuration} $(\si,\io,\pi)$ in $\A$ to
be a collection of objects $\si(J)$ and morphisms $\io(J,K)$ or
$\pi(J,K):\si(J)\ra\si(K)$ in $\A$ satisfying certain axioms, where
$J,K$ are subsets of~$I$. Configurations are a tool for describing
{\it how an object\/ $X$ in $\A$ decomposes into subobjects}.

The first paper \cite{Joyc3} defined configurations and developed
their basic properties, and studied moduli spaces of
$(I,\pr)$-configurations in $\A$, using the theory of {\it Artin
stacks}. It proved well-behaved moduli stacks $\fM(I,\pr)_\A$ exist
when $\A$ is the abelian category of coherent sheaves on a
projective scheme $P$, or of representations of a quiver~$Q$.

This paper develops versions of {\it Ringel--Hall algebras}
\cite{Ring2,Ring3} in the framework of configurations and
Artin stacks. The idea of Ringel--Hall algebras is to make
a $\Q$-algebra $\H$ from an abelian category $\A$. In the
simplest version, the isomorphism classes $[X]$ of objects
$X\in\A$ form a basis for $\H$ with multiplication $[X]*[Z]=
\sum_{\sst[Y]}g^{\sst Y}_{\sst XZ}[Y]$, where $g^{\sst Y}_{\sst XZ}$
is the `number' of exact sequences $0\ra X\ra Y\ra Z\ra 0$ in $\A$.
The important point is that $*$ is {\it associative}.

Ringel--Hall type algebras are defined in four main contexts:
\begin{itemize}
\setlength{\itemsep}{0pt}
\setlength{\parsep}{0pt}
\item {\bf Counting subobjects over finite fields}, as in
Ringel~\cite{Ring2,Ring3}.
\item {\bf Perverse sheaves on moduli spaces} are used by
Lusztig~\cite{Lusz}.
\item {\bf Homology of moduli spaces}, as in Nakajima~\cite{Naka}.
\item {\bf Constructible functions on moduli spaces} are
used by Lusztig \cite[\S 10.18--\S 10.19]{Lusz},
Nakajima \cite[\S 10]{Naka}, Frenkel, Malkin and Vybornov
\cite{FMV}, Riedtmann \cite{Ried} and others.
\end{itemize}
In the first half of the paper we follow the latter path. After some
background on stacks and configurations in \S\ref{ai2} and
\S\ref{ai3}, we begin in \S\ref{ai4} with a detailed account of {\it
Ringel--Hall algebras $\CF(\fObj_\A)$ of constructible functions on
Artin stacks}, using the constructible functions theory developed
in~\cite{Joyc1}.

A distinctive feature of our treatment is the use of {\it
configurations}. Working with configuration moduli stacks
$\fM(I,\pr)_\A$ and 1-morphisms between them makes the
proofs more systematic, and also suggests new ideas. In
particular we construct representations of Ringel--Hall
algebras in a way that appears to be new, and define
bialgebras and other algebraic structures. The only other
paper known to the author using stacks in this way is the
brief sketch in Kapranov and Vasserot \cite[\S 3]{KaVa},
but stacks appear to be the most natural setting.

The second half of our paper \S\ref{ai5}--\S\ref{ai6} studies
various Ringel--Hall algebras of {\it stack functions}
$\uSF(\fObj_\A),\uSF(\fObj_\A,\Up,\La),\ldots$. Stack functions are
a universal generalization of constructible functions on stacks
introduced in \cite{Joyc2}, which contain much more information than
constructible functions. When $\Ext^i(X,Y)=0$ for all $X,Y\in\A$ and
$i>1$, \S\ref{ai6} constructs interesting {\it algebra morphisms}
$\Phi^{\sst\La},\Psi^{\sst\La}, \Psi^{\sst\La^\ci},\Psi^{\sst\Om}$
from stack algebras $\uSF(\fObj_\A),\oSFa(\fObj_\A,*,*)$ to certain
explicit algebras $A(\A,\La,\chi),B(\A,\La,\chi),B(\A,\La^\ci,\chi),
C(\A,\Om,\chi)$. When $\A=\coh(P)$ for $P$ a Calabi--Yau 3-fold, the
same techniques give a {\it Lie algebra morphism}~$\Psi^{\sst\Om}:
\oSFai(\fObj_\A,\Th,\Om)\ra C^\ind(\A,\Om,\chi)$.

These ideas will be applied in the sequels \cite{Joyc4,Joyc5}. Given
a {\it stability condition} $(\tau,T,\le)$ on $\A$, we will define
stack functions $\bde_{\rm ss}^\al(\tau)$ in $\SFa(\fObj_\A)$
parametrizing $\tau$-semistable objects in class $\al$. These
satisfy many {\it identities} in the stack algebra $\SFa(\fObj_\A)$.
Applying $\Phi^{\sst\La},\ldots,\Psi^{\sst\Om}$ to $\bde_{\rm
ss}^\al(\tau)$ yields {\it invariants} of $\A,(\tau,T,\le)$ in
$A(\A,\La,\chi),\ldots,C(\A,\Om,\chi)$, with interesting
transformation laws.
\medskip

\noindent{\it Acknowledgements.} I thank Tom Bridgeland, Frances
Kirwan, Ian Grojnowksi, Alastair King, Richard Thomas and Burt
Totaro for useful conversations. I held an EPSRC Advanced Research
Fellowship whilst writing this paper.

\section{Background material}
\label{ai2}

We begin with some background material on Artin stacks,
constructible functions, stack functions, and motivic invariants,
drawn mostly from~\cite{Joyc1,Joyc2}.

\subsection{Introduction to Artin $\K$-stacks}
\label{ai21}

Fix an algebraically closed field $\K$ throughout. There
are four main classes of `spaces' over $\K$ used in
algebraic geometry, in increasing order of generality:
\begin{equation*}
\text{$\K$-varieties}\subset
\text{$\K$-schemes}\subset
\text{algebraic $\K$-spaces}\subset
\text{algebraic $\K$-stacks}.
\end{equation*}

{\it Algebraic stacks} (also known as Artin stacks) were
introduced by Artin, generalizing {\it Deligne--Mumford stacks}.
For a good introduction to algebraic stacks see G\'omez
\cite{Gome}, and for a thorough treatment see Laumon and
Moret-Bailly \cite{LaMo}. We make the convention that all
algebraic $\K$-stacks in this paper are {\it locally of
finite type}, and $\K$-substacks are {\it locally closed}.

Algebraic $\K$-stacks form a 2-{\it category}. That is, we have {\it
objects} which are $\K$-stacks $\fF,\fG$, and also two kinds of
morphisms, 1-{\it morphisms} $\phi,\psi:\fF\ra\fG$ between
$\K$-stacks, and 2-{\it morphisms} $A:\phi\ra\psi$ between
1-morphisms. An analogy to keep in mind is a 2-category of
categories, where objects are categories, 1-morphisms are functors
between the categories, and 2-morphisms are isomorphisms (natural
transformations) between functors.

We define the set of $\K$-{\it points} of a stack.

\begin{dfn} Let $\fF$ be a $\K$-stack. Write $\fF(\K)$ for the set of
2-isomorphism classes $[x]$ of 1-morphisms $x:\Spec\K\ra\fF$.
Elements of $\fF(\K)$ are called $\K$-{\it points}, or {\it
geometric points}, of $\fF$. If $\phi:\fF\ra\fG$ is a 1-morphism
then composition with $\phi$ induces a map of
sets~$\phi_*:\fF(\K)\ra\fG(\K)$.

For a 1-morphism $x:\Spec\K\ra\fF$, the {\it stabilizer group}
$\Iso_\K(x)$ is the group of 2-morphisms $x\ra x$. When $\fF$ is an
algebraic $\K$-stack, $\Iso_\K(x)$ is an {\it algebraic $\K$-group}.
We say that $\fF$ {\it has affine geometric stabilizers} if
$\Iso_\K(x)$ is an affine algebraic $\K$-group for all
1-morphisms~$x:\Spec\K\ra\fF$.

As an algebraic $\K$-group up to isomorphism, $\Iso_\K(x)$
depends only on the isomorphism class $[x]\in\fF(\K)$ of $x$
in $\Hom(\Spec\K,\fF)$. If $\phi:\fF\ra\fG$ is a 1-morphism,
composition induces a morphism of algebraic $\K$-groups
$\phi_*:\Iso_\K([x])\ra\Iso_\K\bigr(\phi_*([x])\bigr)$,
for~$[x]\in\fF(\K)$.
\label{ai2def1}
\end{dfn}

One important difference in working with 2-categories rather than
ordinary categories is that in diagram-chasing one only requires
1-morphisms to be 2-{\it isomorphic} rather than {\it equal}. The
simplest kind of {\it commutative diagram} is:
\begin{equation*}
\xymatrix@R=6pt{
& \fG \ar@{=>}[d]^{\,F} \ar[dr]^\psi \\
\fF \ar[ur]^\phi \ar[rr]_\chi && \fH, }
\end{equation*}
by which we mean that $\fF,\fG,\fH$ are $\K$-stacks,
$\phi,\psi,\chi$ are 1-morphisms, and $F:\psi\ci\phi\ra\chi$ is a
2-isomorphism. Usually we omit $F$, and mean
that~$\psi\ci\phi\cong\chi$.

\begin{dfn} Let $\phi:\fF\ra\fH$, $\psi:\fG\ra\fH$ be 1-morphisms
of $\K$-stacks. Then one can define the {\it fibre product stack\/}
$\fF\t_{\phi,\fH,\psi}\fG$, or $\fF\t_\fH\fG$ for short, with
1-morphisms $\pi_\fF,\pi_\fG$ fitting into a commutative diagram:
\e
\begin{gathered}
\xymatrix@R=-4pt{
& \fF \ar[dr]^\phi \ar@{=>}[dd] \\
\fF\t_\fH\fG
\ar[dr]_{\pi_\fG} \ar[ur]^{\pi_\fF} && \fH.\\
& \fG \ar[ur]_\psi \\
}
\end{gathered}
\label{ai2eq1}
\e
A commutative diagram
\begin{equation*}
\xymatrix@R=-4pt{
& \fF \ar[dr]^\phi \ar@{=>}[dd] \\
\fE
\ar[dr]_\eta \ar[ur]^\th && \fH\\
& \fG \ar[ur]_\psi \\
}
\end{equation*}
is a {\it Cartesian square} if it is isomorphic to \eq{ai2eq1}, so
there is a 1-isomorphism $\fE\cong\fF\t_\fH\fG$. Cartesian squares
may also be characterized by a universal property.
\label{ai2def2}
\end{dfn}

\subsection{Constructible functions on stacks}
\label{ai22}

Next we discuss {\it constructible functions} on $\K$-stacks,
following \cite{Joyc1}. For this section we need $\K$ to have {\it
characteristic zero}.

\begin{dfn} Let $\fF$ be an algebraic $\K$-stack. We call
$C\subseteq\fF(\K)$ {\it constructible} if $C=\bigcup_{i\in I}
\fF_i(\K)$, where $\{\fF_i:i\in I\}$ is a finite collection of
finite type algebraic $\K$-substacks $\fF_i$ of $\fF$. We call
$S\subseteq\fF(\K)$ {\it locally constructible} if $S\cap C$
is constructible for all constructible~$C\subseteq\fF(\K)$.

A function $f:\fF(\K)\ra\Q$ is called {\it constructible} if
$f(\fF(\K))$ is finite and $f^{-1}(c)$ is a constructible set
in $\fF(\K)$ for each $c\in f(\fF(\K))\sm\{0\}$. A function
$f:\fF(\K)\ra\Q$ is called {\it locally constructible} if
$f\cdot\de_C$ is constructible for all constructible
$C\subseteq\fF(\K)$, where $\de_C$ is the characteristic
function of $C$. Write $\CF(\fF)$ and $\LCF(\fF)$ for the
$\Q$-vector spaces of $\Q$-valued constructible and
locally constructible functions on~$\fF$.
\label{ai2def3}
\end{dfn}

Following \cite[Def.s~4.8, 5.1 \& 5.5]{Joyc1} we define {\it
pushforwards} and {\it pullbacks} of constructible functions along
1-morphisms.

\begin{dfn} Let $\fF$ be an algebraic $\K$-stack with affine
geometric stabilizers and $C\subseteq\fF(\K)$ be constructible.
Then \cite[Def.~4.8]{Joyc1} defines the {\it na\"\i ve Euler
characteristic} $\chi^\na(C)$ of $C$. It is called {\it na\"\i ve}
as it takes no account of stabilizer groups. For $f\in\CF(\fF)$,
define $\chi^\na(\fF,f)$ in $\Q$ by
\begin{equation*}
\chi^\na(\fF,f)=\ts\sum_{c\in f(\fF(\K))\sm\{0\}}c\,\chi^\na
\bigl(f^{-1}(c)\bigr).
\end{equation*}

Let $\fF,\fG$ be algebraic $\K$-stacks with affine geometric
stabilizers, and $\phi:\fF\ra\fG$ a representable 1-morphism.
Then for any $x\in\fF(\K)$ we have an injective morphism
$\phi_*:\Iso_\K(x)\ra\Iso_\K\bigl(\phi_*(x)\bigr)$ of affine
algebraic $\K$-groups. The image $\phi_*\bigl(\Iso_\K(x)\bigr)$
is an affine algebraic $\K$-group closed in $\Iso_\K\bigl(
\phi_*(x)\bigr)$, so the quotient $\Iso_\K\bigl(\phi_*(x)\bigr)
/\phi_*\bigl(\Iso_\K(x)\bigr)$ exists as a quasiprojective
$\K$-variety. Define a function $m_\phi:\fF(\K)\ra\Z$ by
$m_\phi(x)=\chi\bigl(\Iso_\K(\phi_*(x))/\phi_*(\Iso_\K(x))
\bigr)$ for~$x\in\fF(\K)$.

For $f\in\CF(\fF)$, define $\CF^\stk(\phi)f:\fG(\K)\ra\Q$ by
\begin{equation*}
\CF^\stk(\phi)f(y)=\chi^\na\bigl(\fF,m_\phi\cdot f\cdot
\de_{\phi_*^{-1}(y)}\bigr) \quad\text{for $y\in\fG(\K)$,}
\end{equation*}
where $\de_{\smash{\phi_*^{-1}(y)}}$ is the characteristic function
of $\phi_*^{-1}(\{y\})\subseteq\fG(\K)$ on $\fG(\K)$. Then
$\CF^\stk(\phi):\CF(\fF)\ra\CF(\fG)$ is a $\Q$-linear map
called the {\it stack pushforward}.

Let $\th:\fF\ra\fG$ be a finite type 1-morphism. If $C\subseteq
\fG(\K)$ is constructible then so is $\th_*^{-1}(C)\subseteq\fF(\K)$.
It follows that if $f\in\CF(\fG)$ then $f\ci\th_*$ lies in $\CF(\fF)$.
Define the {\it pullback\/} $\th^*:\CF(\fG)\ra\CF(\fF)$ by $\th^*(f)=
f\ci\th_*$. It is a linear map.
\label{ai2def4}
\end{dfn}

Here \cite[Th.s~5.4, 5.6 \& Def.~5.5]{Joyc1} are some properties
of these.

\begin{thm} Let\/ $\fE,\fF,\fG,\fH$ be algebraic $\K$-stacks with
affine geometric stabilizers, and\/ $\be:\fF\ra\fG$, $\ga:\fG\ra\fH$
be $1$-morphisms. Then
\ea
\CF^\stk(\ga\ci\be)&=\CF^\stk(\ga)\ci\CF^\stk(\be):\CF(\fF)\ra\CF(\fH),
\label{ai2eq2}\\
(\ga\ci\be)^*&=\be^*\ci\ga^*:\CF(\fH)\ra\CF(\fF),
\label{ai2eq3}
\ea
supposing $\be,\ga$ representable in \eq{ai2eq2}, and of finite type
in \eq{ai2eq3}. If
\e
\begin{gathered}
\xymatrix{
\fE \ar[r]_\eta \ar[d]^\th & \fG \ar[d]_\psi \\
\fF \ar[r]^\phi & \fH
}
\end{gathered}
\quad
\begin{gathered}
\text{is a Cartesian square with}\\
\text{$\eta,\phi$ representable and}\\
\text{$\th,\psi$ of finite type, then}\\
\text{the following commutes:}
\end{gathered}
\quad
\begin{gathered}
\xymatrix@C=35pt{
\CF(\fE) \ar[r]_{\CF^\stk(\eta)} & \CF(\fG) \\
\CF(\fF) \ar[r]^{\CF^\stk(\phi)} \ar[u]_{\th^*}
& \CF(\fH). \ar[u]^{\psi^*}
}
\end{gathered}
\label{ai2eq4}
\e
\label{ai2thm1}
\end{thm}

As discussed in \cite[\S 3.3]{Joyc1} for the $\K$-scheme case,
equation \eq{ai2eq2} is {\it false} for algebraically closed fields
$\K$ of characteristic $p>0$. This is our reason for restricting to
$\K$ of characteristic zero in \S\ref{ai4}. In \cite[\S 5.3]{Joyc1}
we extend Definition \ref{ai2def4} and Theorem \ref{ai2thm1} to {\it
locally constructible functions}.

\subsection{Stack functions}
\label{ai23}

{\it Stack functions} are a universal generalization of
constructible functions introduced in \cite[\S 3]{Joyc2}. Here
\cite[Def.~3.1]{Joyc2} is the basic definition. Throughout $\K$ is
algebraically closed of arbitrary characteristic, except when we
specify~$\cha\K=0$.

\begin{dfn} Let $\fF$ be an algebraic $\K$-stack with affine
geometric stabilizers. Consider pairs $(\fR,\rho)$, where $\fR$ is a
finite type algebraic $\K$-stack with affine geometric stabilizers
and $\rho:\fR\ra\fF$ is a 1-morphism. We call two pairs
$(\fR,\rho)$, $(\fR',\rho')$ {\it equivalent\/} if there exists a
1-isomorphism $\io:\fR\ra\fR'$ such that $\rho'\ci\io$ and $\rho$
are 2-isomorphic 1-morphisms $\fR\ra\fF$. Write $[(\fR,\rho)]$ for
the equivalence class of $(\fR,\rho)$. If $(\fR,\rho)$ is such a
pair and $\fS$ is a closed $\K$-substack of $\fR$ then
$(\fS,\rho\vert_\fS)$, $(\fR\sm\fS,\rho\vert_{\fR\sm\fS})$ are pairs
of the same kind.

Define $\uSF(\fF)$ to be the $\Q$-vector space generated by
equivalence classes $[(\fR,\rho)]$ as above, with for each closed
$\K$-substack $\fS$ of $\fR$ a relation
\e
[(\fR,\rho)]=[(\fS,\rho\vert_\fS)]+[(\fR\sm\fS,\rho\vert_{\fR\sm\fS})].
\label{ai2eq5}
\e
Define $\SF(\fF)$ to be the $\Q$-vector space generated by
$[(\fR,\rho)]$ with $\rho$ {\it representable}, with the same
relations \eq{ai2eq5}. Then~$\SF(\fF)\subseteq\uSF(\fF)$.
\label{ai2def5}
\end{dfn}

Elements of $\uSF(\fF)$ will be called {\it stack functions}. In
\cite[Def.~3.2]{Joyc2} we relate $\CF(\fF)$ and~$\SF(\fF)$.

\begin{dfn} Let $\fF$ be an algebraic $\K$-stack with affine
geometric stabilizers, and $C\subseteq\fF(\K)$ be constructible.
Then $C=\coprod_{i=1}^n\fR_i(\K)$, for $\fR_1,\ldots,\fR_n$ finite
type $\K$-substacks of $\fF$. Let $\rho_i:\fR_i\ra\fF$ be the
inclusion 1-morphism. Then $[(\fR_i,\rho_i)]\in\SF(\fF)$. Define
$\bde_C=\ts\sum_{i=1}^n[(\fR_i,\rho_i)]\in\SF(\fF)$. We think of
this stack function as the analogue of the characteristic function
$\de_C\in\CF(\fF)$ of $C$. Define a $\Q$-linear map
$\io_\fF:\CF(\fF)\ra\SF(\fF)$ by $\io_\fF(f)=\ts\sum_{0\ne c\in
f(\fF(\K))}c\cdot\bde_{f^{-1}(c)}$. For $\K$ of characteristic zero,
define a $\Q$-linear map $\pi_\fF^\stk:\SF(\fF)\ra\CF(\fF)$ by
\begin{equation*}
\pi_\fF^\stk\bigl(\ts\sum_{i=1}^nc_i[(\fR_i,\rho_i)]\bigr)=
\ts\sum_{i=1}^nc_i\CF^\stk(\rho_i)1_{\fR_i},
\end{equation*}
where $1_{\fR_i}$ is the function 1 in $\CF(\fR_i)$. Then
\cite[Prop.~3.3]{Joyc2} shows $\pi_\fF^\stk\ci\io_\fF$ is the
identity on $\CF(\fF)$. Thus, $\io_\fF$ is injective and
$\pi_\fF^\stk$ is surjective. In general $\io_\fF$ is far from
surjective, and $\uSF,\SF(\fF)$ are much larger than~$\CF(\fF)$.
\label{ai2def6}
\end{dfn}

All the operations of constructible functions in \S\ref{ai22} extend
to stack functions.

\begin{dfn} Define {\it multiplication} `$\,\cdot\,$' on $\uSF(\fF)$ by
\e
[(\fR,\rho)]\cdot[(\fS,\si)]=[(\fR\t_{\rho,\fF,\si}\fS,\rho\ci\pi_\fR)].
\label{ai2eq6}
\e
This extends to a $\Q$-bilinear product $\uSF(\fF)\t\uSF(\fF)\ra
\uSF(\fF)$ which is commutative and associative, and $\SF(\fF)$ is
closed under `$\,\cdot\,$'. Let $\phi:\fF\!\ra\!\fG$ be a 1-morphism
of algebraic $\K$-stacks with affine geometric stabilizers. Define
the {\it pushforward\/} $\phi_*:\uSF(\fF)\!\ra\!\uSF(\fG)$~by
\e
\phi_*:\ts\sum_{i=1}^mc_i[(\fR_i,\rho_i)]\longmapsto
\ts\sum_{i=1}^mc_i[(\fR_i,\phi\ci\rho_i)].
\label{ai2eq7}
\e
If $\phi$ is representable then $\phi_*$ maps $\SF(\fF)\!\ra\!
\SF(\fG)$. For $\phi$ of finite type, define {\it pullbacks}
$\phi^*:\uSF(\fG)\!\ra\!\uSF(\fF)$,
$\phi^*:\SF(\fG)\!\ra\!\SF(\fF)$~by
\e
\phi^*:\ts\sum_{i=1}^mc_i[(\fR_i,\rho_i)]\longmapsto
\ts\sum_{i=1}^mc_i[(\fR_i\t_{\rho_i,\fG,\phi}\fF,\pi_\fF)].
\label{ai2eq8}
\e
The {\it tensor product\/} $\ot\!:\!\uSF(\fF)\!\t\!\uSF(\fG)\!\ra
\!\uSF(\fF\!\t\!\fG)$ or $\SF(\fF)\!\t\!\SF(\fG)\!\ra\!
\SF(\fF\!\t\!\fG)$~is
\e
\bigl(\ts\sum_{i=1}^mc_i[(\fR_i,\rho_i)]\bigr)\!\ot\!
\bigl(\ts\sum_{j=1}^nd_j[(\fS_j,\si_j)]\bigr)\!=\!\ts
\sum_{i,j}c_id_j[(\fR_i\!\t\!\fS_j,\rho_i\!\t\!\si_j)].
\label{ai2eq9}
\e
\label{ai2def7}
\end{dfn}

Here \cite[Th.~3.5]{Joyc2} is the analogue of Theorem~\ref{ai2thm1}.

\begin{thm} Let\/ $\fE,\fF,\fG,\fH$ be algebraic $\K$-stacks with
affine geometric stabilizers, and\/ $\be:\fF\ra\fG$, $\ga:\fG\ra\fH$
be $1$-morphisms. Then
\begin{align*}
(\ga\!\ci\!\be)_*\!&=\!\ga_*\!\ci\!\be_*:\uSF(\fF)\!\ra\!\uSF(\fH),&
(\ga\!\ci\!\be)_*\!&=\!\ga_*\!\ci\!\be_*:\SF(\fF)\!\ra\!\SF(\fH),\\
(\ga\!\ci\!\be)^*\!&=\!\be^*\!\ci\!\ga^*:\uSF(\fH)\!\ra\!\uSF(\fF),&
(\ga\!\ci\!\be)^*\!&\!=\!\be^*\!\ci\!\ga^*:\SF(\fH)\!\ra\!\SF(\fF),
\end{align*}
for $\be,\ga$ representable in the second equation, and of finite
type in the third and fourth. If\/ $f,g\in\uSF(\fG)$ and\/ $\be$ is
finite type then $\be^*(f\cdot g)=\be^*(f)\cdot\be^*(g)$. If
\begin{equation*}
\begin{gathered}
\xymatrix@R=15pt{
\fE \ar[r]_\eta \ar[d]^{\,\th} & \fG \ar[d]_{\psi\,} \\
\fF \ar[r]^\phi & \fH
}
\end{gathered}
\quad
\begin{gathered}
\text{is a Cartesian square with}\\
\text{$\th,\psi$ of finite type, then}\\
\text{the following commutes:}
\end{gathered}
\quad
\begin{gathered}
\xymatrix@C=35pt@R=10pt{
\uSF(\fE) \ar[r]_{\eta_*} & \uSF(\fG) \\
\uSF(\fF) \ar[r]^{\phi_*} \ar[u]_{\,\th^*}
& \uSF(\fH). \ar[u]^{\psi^*\,}
}
\end{gathered}
\end{equation*}
The same applies for $\SF(\fE),\ldots,\SF(\fH)$ if\/ $\eta,\phi$ are
representable.
\label{ai2thm2}
\end{thm}

In \cite[Prop.~3.7 \& Th.~3.8]{Joyc2} we relate pushforwards and
pullbacks of stack and constructible functions
using~$\io_\fF,\pi_\fF^\stk$.

\begin{thm} Let\/ $\K$ have characteristic zero, $\fF,\fG$ be
algebraic $\K$-stacks with affine geometric stabilizers, and\/
$\phi:\fF\ra\fG$ be a $1$-morphism. Then
\begin{itemize}
\setlength{\itemsep}{0pt}
\setlength{\parsep}{0pt}
\item[{\rm(a)}] $\phi^*\!\ci\!\io_\fG\!=\!\io_\fF\!\ci\!\phi^*:
\CF(\fG)\!\ra\!\SF(\fF)$ if\/ $\phi$ is of finite type;
\item[{\rm(b)}] $\pi^\stk_\fG\ci\phi_*=\CF^\stk(\phi)\ci\pi_\fF^\stk:
\SF(\fF)\ra\CF(\fG)$ if\/ $\phi$ is representable; and
\item[{\rm(c)}] $\pi^\stk_\fF\ci\phi^*=\phi^*\ci\pi_\fG^\stk:
\SF(\fG)\ra\CF(\fF)$ if\/ $\phi$ is of finite type.
\end{itemize}
\label{ai2thm3}
\end{thm}

In \cite[\S 3]{Joyc2} we extend all the material on $\uSF,\SF(\fF)$
to {\it local stack functions} $\uLSF,\LSF(\fF)$, the analogues of
locally constructible functions. The main differences are in which
1-morphisms must be of finite type.

\subsection{Motivic invariants of Artin stacks}
\label{ai24}

In \cite[\S 4]{Joyc2} we extend {\it motivic} invariants of
quasiprojective $\K$-varieties to Artin stacks. We need the
following data,~\cite[Assumptions 4.1 \& 6.1]{Joyc2}.

\begin{ass} Suppose $\La$ is a commutative $\Q$-algebra with
identity 1, and
\begin{equation*}
\Up:\{\text{isomorphism classes $[X]$ of quasiprojective
$\K$-varieties $X$}\}\longra\La
\end{equation*}
a map for $\K$ an algebraically closed field, satisfying the
following conditions:
\begin{itemize}
\setlength{\itemsep}{0pt}
\setlength{\parsep}{0pt}
\item[(i)] If $Y\subseteq X$ is a closed subvariety then
$\Up([X])=\Up([X\sm Y])+\Up([Y])$;
\item[(ii)] If $X,Y$ are quasiprojective $\K$-varieties then
$\Up([X\!\t\!Y])\!=\!\Up([X])\Up([Y])$;
\item[(iii)] Write $\el=\Up([\K])$ in $\La$, regarding $\K$ as a
$\K$-variety, the affine line (not the point $\Spec\K$). Then $\el$
and $\el^k-1$ for $k=1,2,\ldots$ are invertible in~$\La$.
\end{itemize}
Suppose $\La^\ci$ is a $\Q$-subalgebra of $\La$ containing the image
of $\Up$ and the elements $\el^{-1}$ and
$(\el^k+\el^{k-1}+\cdots+1)^{-1}$ for $k=1,2,\ldots$, but {\it
not\/} containing $(\el-1)^{-1}$. Let $\Om$ be a commutative
$\Q$-algebra, and $\pi:\La^\ci\ra\Om$ a surjective $\Q$-algebra
morphism, such that $\pi(\el)=1$. Define
\begin{equation*}
\Th:\{\text{isomorphism classes $[X]$ of quasiprojective
$\K$-varieties $X$}\}\longra\Om
\end{equation*}
by $\Th=\pi\ci\Up$. Then~$\Th([\K])=1$.
\label{ai2ass}
\end{ass}

We chose the notation `$\el$' as in motivic integration $[\K]$ is
called the {\it Tate motive} and written $\mathbb L$. We have
$\Up\bigl([\GL(m,\K)]\bigr)=\el^{m(m-1)/2}\prod_{k=1}^m(\el^k-1)$,
so (iii) ensures $\Up([\GL(m,\K)])$ is invertible in $\La$ for all
$m\ge 1$. Here \cite[Ex.s 4.3 \& 6.3]{Joyc2} is an example of
suitable $\La,\Up,\ldots$; more are given in~\cite[\S 4.1 \& \S
6.1]{Joyc2}.

\begin{ex} Let $\K$ be an algebraically closed field. Define
$\La=\Q(z)$, the algebra of rational functions in $z$ with
coefficients in $\Q$. For any quasiprojective $\K$-variety $X$, let
$\Up([X])=P(X;z)$ be the {\it virtual Poincar\'e polynomial\/} of
$X$. This has a complicated definition in \cite[Ex.~4.3]{Joyc2}
which we do not repeat, involving Deligne's weight filtration when
$\cha\K=0$ and the action of the Frobenius on $l$-adic cohomology
when $\cha\K>0$. If $X$ is smooth and projective then $P(X;z)$ is
the ordinary Poincar\'e polynomial $\sum_{k=0}^{2\dim X}b^k(X)z^k$,
where $b^k(X)$ is the $k^{\rm th}$ Betti number in $l$-adic
cohomology, for $l$ coprime to $\cha\K$. Also~$\el=P(\K;z)=z^2$.

Let $\La^\ci$ be the subalgebra of $P(z)/Q(z)$ in $\La$ for which
$z\pm 1$ do not divide $Q(z)$. Here are two possibilities for
$\Om,\pi$. Assumption \ref{ai2ass} holds in each case.
\begin{itemize}
\setlength{\itemsep}{0pt}
\setlength{\parsep}{0pt}
\item[(a)] Set $\Om=\Q$ and $\pi:f(z)\mapsto f(-1)$. Then
$\Th([X])=\pi\ci\Up([X])$ is the {\it Euler characteristic} of~$X$.
\item[(b)] Set $\Om=\Q$ and $\pi:f(z)\mapsto f(1)$. Then
$\Th([X])=\pi\ci\Up([X])$ is the {\it sum of the virtual Betti
numbers} of~$X$.
\end{itemize}
\label{ai2ex}
\end{ex}

We need a few facts about {\it algebraic $\K$-groups}. A good
reference is Borel \cite{Bore}. Following Borel, we define a
$\K$-{\it variety} to be a $\K$-scheme which is reduced, separated,
and of finite type, but {\it not\/} necessarily irreducible. An
algebraic $\K$-group is then a $\K$-variety $G$ with identity $1\in
G$, multiplication $\mu:G\t G\ra G$ and inverse $i:G\ra G$ (as
morphisms of $\K$-varieties) satisfying the usual group axioms. We
call $G$ {\it affine} if it is an affine $\K$-variety. {\it
Special\/} $\K$-groups are studied by Serre and Grothendieck
in~\cite[\S\S 1, 5]{Chev}.

\begin{dfn} An algebraic $\K$-group $G$ is called {\it special\/} if
every principal $G$-bundle is locally trivial. Properties of special
$\K$-groups can be found in \cite[\S\S 1.4, 1.5 \& 5.5]{Chev} and
\cite[\S 2.1]{Joyc2}. In \cite[Lem.~4.6]{Joyc2} we show that if
Assumption \ref{ai2ass} holds and $G$ is special then $\Up([G])$ is
invertible in~$\La$.
\label{ai2def8}
\end{dfn}

In \cite[Th.~4.9]{Joyc2} we extend $\Up$ to Artin stacks, using
Definition~\ref{ai2def8}.

\begin{thm} Let Assumption \ref{ai2ass} hold. Then there exists a
unique morphism of\/ $\Q$-algebras $\Up':\uSF(\Spec\K)\ra\La$ such
that if\/ $G$ is a special algebraic $\K$-group acting on a
quasiprojective $\K$-variety $X$ then~$\Up'\bigl(\bigl[[X/G]
\bigr]\bigr)=\Up([X])/\Up([G])$.
\label{ai2thm4}
\end{thm}

Thus, if $\fR$ is a finite type algebraic $\K$-stack with affine
geometric stabilizers the theorem defines $\Up'([\fR])\in\La$.
Taking $\La,\Up$ as in Example \ref{ai2ex} yields the {\it virtual
Poincar\'e function} $P(\fR;z)=\Up'([\fR])$ of $\fR$, a natural
extension of virtual Poincar\'e polynomials to stacks. Clearly,
Theorem \ref{ai2thm4} only makes sense if $\Up([G])^{-1}$ exists for
all special $\K$-groups $G$. This excludes the Euler characteristic
$\Up=\chi$, for instance, since $\chi([\K^\t])=0$ is not invertible.
We overcome this in \cite[\S 6]{Joyc2} by defining a finer extension
of $\Up$ to stacks that keeps track of maximal tori of stabilizer
groups, and allows $\Up=\chi$. This can then be used with $\Th,\Om$
in Assumption~\ref{ai2ass}.

\subsection{Stack functions over motivic invariants}
\label{ai25}

In \cite[\S 4--\S 6]{Joyc2} we integrate the stack functions of
\S\ref{ai23} with the motivic invariant ideas of \S\ref{ai24} to
define more stack function spaces.

\begin{dfn} Let Assumption \ref{ai2ass} hold, and $\fF$ be an
algebraic $\K$-stack with affine geometric stabilizers. Consider
pairs $(\fR,\rho)$, with equivalence, as in Definition
\ref{ai2def5}. Define $\uSF(\fF,\Up,\La)$ to be the $\La$-module
generated by equivalence classes $[(\fR,\rho)]$, with the following
relations:
\begin{itemize}
\setlength{\itemsep}{0pt}
\setlength{\parsep}{0pt}
\item[(i)] Given $[(\fR,\rho)]$ as above and $\fS$ a closed $\K$-substack
of $\fR$ we have $[(\fR,\rho)]=[(\fS,\rho\vert_\fS)]+[(\fR\sm\fS,
\rho\vert_{\fR\sm\fS})]$, as in~\eq{ai2eq5}.
\item[(ii)] Let $\fR$ be a finite type algebraic $\K$-stack with
affine geometric stabilizers, $U$ a quasiprojective $\K$-variety,
$\pi_\fR:\fR\t U\ra\fR$ the natural projection, and $\rho:\fR\ra\fF$
a 1-morphism. Then~$[(\fR\t U,\rho\ci\pi_\fR)]=
\Up([U])[(\fR,\rho)]$.
\item[(iii)] Given $[(\fR,\rho)]$ as above and a 1-isomorphism
$\fR\cong[X/G]$ for $X$ a quasiprojective $\K$-variety and $G$ a
special algebraic $\K$-group acting on $X$, we have
$[(\fR,\rho)]=\Up([G])^{-1}[(X,\rho\ci\pi)]$, where
$\pi:X\ra\fR\cong[X/G]$ is the natural projection 1-morphism.
\end{itemize}
Define a $\Q$-linear projection $\Pi^{\Up,\La}_\fF:\uSF(\fF)\ra
\uSF(\fF,\Up,\La)$ by
\begin{equation*}
\Pi^{\Up,\La}_\fF:\ts\sum_{i\in I}c_i[(\fR_i,\rho_i)]\longmapsto
\ts\sum_{i\in I}c_i[(\fR_i,\rho_i)],
\end{equation*}
using the embedding $\Q\subseteq\La$ to regard $c_i\in\Q$ as an
element of~$\La$.

We also define variants of these: spaces $\uoSF,\oSF(\fF,\Up,\La)$,
$\uoSF,\oSF(\fF,\Up,\La^\ci)$ and $\uoSF,\oSF(\fF,\Th,\Om)$, which
are the $\La,\La^\ci$- and $\Om$-modules respectively generated by
$[(\fR,\rho)]$ as above, with $\rho$ representable for
$\oSF(\fF,*,*)$, and with relations (i),(ii) above but (iii)
replaced by a finer, more complicated relation
\cite[Def.~5.17(iii)]{Joyc2}. There are natural projections
$\Pi^{\Up,\La}_\fF,\bar\Pi^{\Up,\La}_\fF,\bar\Pi^{\Up,\La^\ci}_\fF,
\bar\Pi^{\Th,\Om}_\fF$ between various of the spaces. We can also
define {\it local stack function} spaces~$\uLSF,\uoLSF,\oLSF
(\fF,*,*)$.
\label{ai2def9}
\end{dfn}

In \cite{Joyc2} we give analogues of Definitions \ref{ai2def6} and
\ref{ai2def7} and Theorems \ref{ai2thm2} and \ref{ai2thm3} for these
spaces. For the analogue of $\pi^\stk_\fF$, suppose ${\rm
X}:\La^\ci\ra\Q$ or ${\rm X}:\Om\ra\Q$ is an algebra morphism with
${\rm X}\ci\Up([U])=\chi([U])$ or ${\rm X}\ci\Th([U])=\chi([U])$ for
varieties $U$, where $\chi$ is the Euler characteristic. Define
$\bar\pi^\stk_\fF:\oSF(\fF,\Up,\La^\ci)\ra\CF(\fF)$ or
$\bar\pi^\stk_\fF:\oSF(\fF,\Th,\Om)\ra\CF(\fF)$ by
\begin{equation*}
\bar\pi_\fF^\stk\bigl(\ts\sum_{i=1}^nc_i[(\fR_i,\rho_i)]\bigr)=
\ts\sum_{i=1}^n{\rm X}(c_i)\CF^\stk(\rho_i)1_{\fR_i}.
\end{equation*}
The operations `$\,\cdot\,$'$,\phi_*,\phi^*,\ot$ on
$\uSF(*,\Up,\La),\ldots,\oSF(*,\Th,\Om)$ are given by the same
formulae. The important point is that \eq{ai2eq6}--\eq{ai2eq9} are
compatible with the relations defining $\uSF(*,\Up,\La),\ldots,
\oSF(*,\Th,\Om)$, or they would not be well-defined.

In \cite[Prop.~4.14]{Joyc2} we identify $\uSF(\Spec\K,\Up,\La)$. The
proof involves showing that $\Up'$ in Theorem \ref{ai2thm4} is
compatible with Definition \ref{ai2def9}(i)--(iii) and so descends
to $\Up':\uSF(\Spec\K,\Up,\La)\ra\La$, which is the inverse
of~$i_\La$.

\begin{prop} The map $i_\La:\La\ra\uSF(\Spec\K,\Up,\La)$ taking
$i_\La:c\mapsto c[\Spec\K]$ is an isomorphism of algebras.
\label{ai2prop1}
\end{prop}

Here \cite[Prop.s 5.21 \& 5.22]{Joyc2} is a useful way of
representing these spaces.

\begin{prop} $\uoSF,\oSF(\fF,\Up,\La)$, $\uoSF,\oSF(\fF,\Up,\La^\ci)$
and\/ $\uoSF,\oSF(\fF,\Th,\Om)$ are generated over $\La,\La^\ci$
and\/ $\Om$ respectively by elements $[(U\t[\Spec\K/T],\rho)]$,
for\/ $U$ a quasiprojective $\K$-variety and\/ $T$ an algebraic
$\K$-group isomorphic to $(\K^\t)^k\t K$ for $k\ge 0$ and\/ $K$
finite abelian.

Suppose $\sum_{i\in I}c_i[(U_i\t[\Spec\K/T_i],\rho_i)]=0$ in one of
these spaces, where $I$ is finite set, $c_i\in\La,\La^\ci$ or $\Om$,
$U_i$ is a quasiprojective $\K$-variety and\/ $T_i$ an algebraic
$\K$-group isomorphic to $(\K^\t)^{k_i}\t K_i$ for $k_i\ge 0$ and\/
$K_i$ finite abelian, with\/ $T_i\not\cong T_j$ for $i\ne j$. Then
$c_j[(U_j\t[\Spec\K/T_j],\rho_j)]=0$ for all\/~$j\in I$.
\label{ai2prop2}
\end{prop}

In \cite[\S 5.2]{Joyc2} we define operators
$\Pi^\mu,\Pi^\vi_n,\hat\Pi^\nu_\fF$ on $\uSF(\fF),\uoSF(\fF,*,*)$
(but not on $\uSF(\fF,\Up,\La)$). Very roughly speaking, $\Pi^\vi_n$
projects $[(\fR,\rho)]\in\uSF(\fF)$ to $[(\fR_n,\rho)]$, where
$\fR_n$ is the $\K$-substack of points $r\in\fR(\K)$ whose
stabilizer groups $\Iso_\K(r)$ have {\it rank\/} $n$, that is,
maximal torus~$(\K^\t)^n$.

Unfortunately, it is more complicated than this. The right notion is
not the actual rank of stabilizer groups, but the {\it virtual
rank}. This is a difficult idea which treats $r\in\fR(\K)$ with
nonabelian stabilizer group $G=\Iso_\K(r)$ as a linear combination
of points with `virtual ranks' in the range $\rk\,C(G)\le
n\le\rk\,G$. Effectively this {\it abelianizes stabilizer groups},
that is, using virtual rank we can treat $\fR$ as though its
stabilizer groups were all abelian, essentially tori $(\K^\t)^n$.
These ideas will be key tools in~\S\ref{ai5}--\S\ref{ai6}.

Here is a way to interpret the spaces of Definition \ref{ai2def9},
explained in \cite{Joyc2}. In \S\ref{ai22}, pushforwards
$\CF^\stk(\phi):\CF(\fF)\ra\CF(\fG)$ are defined by `integration'
over the fibres of $\phi$, using the Euler characteristic $\chi$ as
measure. In the same way, given $\La,\Up$ as in Assumption
\ref{ai2ass} we could consider $\La$-valued constructible functions
$\CF(\fF)_\La$, and define a pushforward $\phi_*:\CF(\fF)_\La\ra
\CF(\fG)_\La$ by `integration' using $\Up$ as measure, instead of
$\chi$. But then $(\psi\ci\phi)_*=\psi_*\ci\phi_*$ may no longer
hold, as this depends on properties of $\chi$ on
non-Zariski-locally-trivial fibrations which are false for other
$\Up$ such as virtual Poincar\'e polynomials.

The space $\uSF(\fF,\Up,\La)$ is very like $\CF(\fF)_\La$ with
pushforwards $\phi_*$ defined using $\Up$, but satisfies
$(\psi\ci\phi)_*=\psi_*\ci\phi_*$ and other useful functoriality
properties. So we can use it as a substitute for $\CF(\fF)$. The
spaces $\uoSF,\oSF(\fF,*,*)$ are similar but also keep track of
information on the maximal tori of stabilizer groups.

\section{Background on configurations from \cite{Joyc3}}
\label{ai3}

We recall in \S\ref{ai31} the main definitions and results on
$(I,\pr)$-configurations that we will need later, and in
\S\ref{ai32} some important facts on moduli stacks of
configurations. For motivation and other results see \cite{Joyc3},
and for background on abelian categories, see Gelfand and
Manin~\cite{GeMa}.

\subsection{Basic definitions}
\label{ai31}

Here is some notation for finite posets, taken from \cite[Def.s~3.2
\& 4.1]{Joyc3}.

\begin{dfn} A {\it finite partially ordered set\/} or {\it
finite poset\/} $(I,\pr)$ is a finite set $I$ with a partial
order $I$. Define $J\subseteq I$ to be an {\it f-set\/} if
$i\in I$ and $h,j\in J$ and $h\pr i\pr j$ implies $i\in J$.
Define $\F_\sIp$ to be the set of f-sets of $I$. Define
$\G_\sIp$ to be the subset of $(J,K)\in\F_\sIp\t\F_\sIp$
such that $J\subseteq K$, and if $j\in J$ and $k\in K$
with $k\pr j$, then $k\in J$. Define $\H_\sIp$ to be the
subset of $(J,K)\in\F_\sIp\t\F_\sIp$ such that
$K\subseteq J$, and if $j\in J$ and $k\in K$ with
$k\pr j$, then~$j\in K$.
\label{ai3def1}
\end{dfn}

We define $(I,\pr)$-{\it configurations},~\cite[Def.~4.1]{Joyc3}.

\begin{dfn} Let $(I,\pr)$ be a finite poset, and use the
notation of Definition \ref{ai3def1}. Define an $(I,\pr)$-{\it
configuration} $(\si,\io,\pi)$ in an abelian category $\A$ to be
maps $\si:\F_\sIp\!\ra\!\Obj(\A)$, $\io:\G_\sIp\!\ra\!\Mor(\A)$, and
$\pi:\H_\sIp\!\ra\!\Mor(\A)$, where
\begin{itemize}
\setlength{\itemsep}{0pt}
\setlength{\parsep}{0pt}
\item[(i)] $\si(J)$ is an object in $\A$ for $J\in\F_\sIp$,
with~$\si(\emptyset)=0$.
\item[(ii)] $\io(J,K):\si(J)\!\ra\!\si(K)$ is injective
for $(J,K)\!\in\!\G_\sIp$, and~$\io(J,J)\!=\!\id_{\si(J)}$.
\item[(iii)] $\pi(J,K)\!:\!\si(J)\!\ra\!\si(K)$ is surjective
for $(J,K)\!\in\!\H_\sIp$, and~$\pi(J,J)\!=\!\id_{\si(J)}$.
\end{itemize}
These should satisfy the conditions:
\begin{itemize}
\setlength{\itemsep}{0pt}
\setlength{\parsep}{0pt}
\item[(A)] Let $(J,K)\in\G_\sIp$ and set $L=K\sm J$. Then the
following is exact in~$\A$:
\begin{equation*}
\xymatrix@C=40pt{ 0 \ar[r] &\si(J) \ar[r]^{\io(J,K)} &\si(K)
\ar[r]^{\pi(K,L)} &\si(L) \ar[r] & 0. }
\end{equation*}
\item[(B)] If $(J,K)\in\G_\sIp$ and $(K,L)\in\G_\sIp$
then~$\io(J,L)=\io(K,L)\ci\io(J,K)$.
\item[(C)] If $(J,K)\in\H_\sIp$ and $(K,L)\in\H_\sIp$
then~$\pi(J,L)=\pi(K,L)\ci\pi(J,K)$.
\item[(D)] If $(J,K)\in\G_\sIp$ and $(K,L)\in\H_\sIp$ then
\begin{equation*}
\pi(K,L)\ci\io(J,K)=\io(J\cap L,L)\ci\pi(J,J\cap L).
\end{equation*}
\end{itemize}

A {\it morphism} $\al:(\si,\io,\pi)\ra(\si',\io',\pi')$ of
$(I,\pr)$-configurations in $\A$ is a collection of morphisms
$\al(J):\si(J)\ra\si'(J)$ for each $J\in\F_\sIp$ satisfying
\begin{align*}
\al(K)\ci\io(J,K)&=\io'(J,K)\ci\al(J)&&
\text{for all $(J,K)\in\G_\sIp$, and}\\
\al(K)\ci\pi(J,K)&=\pi'(J,K)\ci\al(J)&&
\text{for all $(J,K)\in\H_\sIp$.}
\end{align*}
It is an {\it isomorphism} if $\al(J)$ is an isomorphism for
all~$J\in\F_\sIp$.
\label{ai3def2}
\end{dfn}

Here \cite[Def.s~5.1, 5.2]{Joyc3} are two ways to construct new
configurations.

\begin{dfn} Let $(I,\pr)$ be a finite poset and $J\in\F_\sIp$. Then
$(J,\pr)$ is also a finite poset, and $\F_\sJp,\G_\sJp,\H_\sJp\!
\subseteq\!\F_\sIp,\G_\sIp,\H_\sIp$. Let $(\si,\io,\pi)$ be an
$(I,\pr)$-configuration in an abelian category $\A$. Define the
$(J,\pr)$-{\it subconfiguration} $(\si',\io',\pi')$ of
$(\si,\io,\pi)$ by $\si'\!=\!\si\vert_{\F_\sJp}$,
$\io'\!=\!\io\vert_{\G_\sJp}$ and~$\pi'\!=\!\pi\vert_{\H_\sJp}$.

Let $(I,\pr),(K,\tl)$ be finite posets, and $\phi:I\!\ra\!K$ be
surjective with $i\pr j$ implies $\phi(i)\!\tl \!\phi(j)$. Using
$\phi^{-1}$ to pull subsets of $K$ back to $I$ maps
$\F_\sKt,\G_\sKt,\ab\H_\sKt\!\ra\!\F_\sIp,\G_\sIp,\H_\sIp$. Let
$(\si,\io,\pi)$ be an $(I,\pr)$-configuration in an abelian category
$\A$. Define the {\it quotient\/ $(K,\tl)$-configuration}
$(\ti\si,\ti\io,\ti\pi)$ by
$\ti\si(A)\!=\!\si\bigl(\phi^{-1}(A)\bigr)$ for $A\!\in\!\F_\sKt$,
$\ti\io(A,B)\!=\!\io\bigl(\phi^{-1}(A),\ab \phi^{-1}(B)\bigr)$ for
$(A,B)\!\in\!\G_\sKt$, and
$\ti\pi(A,B)\!=\!\pi\bigl(\phi^{-1}(A),\phi^{-1}(B)\bigr)$
for~$(A,B)\!\in\!\H_\sKt$.
\label{ai3def3}
\end{dfn}

\subsection{Moduli stacks of configurations}
\label{ai32}

Here are our initial assumptions.

\begin{ass} Fix an algebraically closed field $\K$. (Throughout
\S\ref{ai4} we will require $\K$ to have {\it characteristic zero}.)
Let $\A$ be an abelian category with $\Hom(X,Y)=\Ext^0(X,Y)$ and
$\Ext^1(X,Y)$ finite-dimensional $\K$-vector spaces for all
$X,Y\in\A$, and all composition maps $\Ext^i(Y,Z)\t\Ext^j(X,Y)
\ra\Ext^{i+j}(X,Z)$ bilinear for $i,j,i+j=0$ or 1. Let $K(\A)$
be the quotient of the Grothendieck group $K_0(\A)$ by some
fixed subgroup. Suppose that if $X\in\Obj(\A)$ with $[X]=0$ in
$K(\A)$ then~$X\cong 0$.

To define moduli stacks of objects or configurations in $\A$, we
need some {\it extra data}, to tell us about algebraic families of
objects and morphisms in $\A$, parametrized by a base scheme $U$. We
encode this extra data as a {\it stack in exact categories} $\fF_\A$
on the 2-{\it category of\/ $\K$-schemes} $\Sch_\K$, made into a
{\it site} with the {\it \'etale topology}. The $\K,\A,K(\A),\fF_\A$
must satisfy some complex additional conditions \cite[Assumptions
7.1 \& 8.1]{Joyc3}, which we do not give. \label{ai3ass}
\end{ass}

We define some notation,~\cite[Def.~7.3]{Joyc3}.

\begin{dfn} We work in the situation of Assumption \ref{ai3ass}.
Define
\begin{equation*}
\bar C(\A)=\bigl\{[X]\in K(\A):X\in\A\bigr\}\subset K(\A).
\end{equation*}
That is, $\bar C(\A)$ is the collection of classes in $K(\A)$ of
objects $X\in\A$. Note that $\bar C(\A)$ is closed under addition,
as $[X\op Y]=[X]+[Y]$. In \cite{Joyc4,Joyc5} we shall make much use
of $C(\A)=\bar C(\A)\sm\{0\}$. We think of $C(\A)$ as the `positive
cone' and $\bar C(\A)$ as the `closed positive cone' in $K(\A)$,
which explains the notation. For $(I,\pr)$ a finite poset and
$\ka:I\ra\bar C(\A)$, define an $(I,\pr,\ka)$-{\it configuration} to
be an $(I,\pr)$-configuration $(\si,\io,\pi)$ with
$[\si(\{i\})]=\ka(i)$ in $K(\A)$ for all~$i\in I$.
\label{ai3def4}
\end{dfn}

In the situation above, we define the following $\K$-stacks
\cite[Def.s 7.2 \& 7.4]{Joyc3}:
\begin{itemize}
\setlength{\itemsep}{0pt}
\setlength{\parsep}{0pt}
\item The {\it moduli stacks} $\fObj_\A$ of {\it objects in} $\A$,
and $\fObj^\al_\A$ of {\it objects in $\A$ with class $\al$ in}
$K(\A)$, for each $\al\in\bar C(\A)$. They are algebraic $\K$-stacks,
locally of finite type, with $\fObj_\A^\al$ an open and closed
$\K$-substack of $\fObj_\A$. The underlying geometric spaces
$\fObj_\A(\K),\fObj_\A^\al(\K)$ are the sets of isomorphism
classes of objects $X$ in $\A$, with $[X]=\al$ for~$\fObj_\A^\al(\K)$.
\item The {\it moduli stacks\/} $\fM(I,\pr)_\A$ of $(I,\pr)$-{\it
configurations} and $\fM(I,\pr,\ka)_\A$ of $(I,\pr,\ka)$-{\it
configurations in} $\A$, for all finite posets $(I,\pr)$ and
$\ka:I\ra\bar C(\A)$. They are algebraic $\K$-stacks, locally of
finite type, with $\fM(I,\pr,\ka)_\A$ an open and closed
$\K$-substack of $\fM(I,\pr)_\A$. Write
$\M(I,\pr)_\A,\M(I,\pr,\ka)_\A$ for the underlying geometric spaces
$\fM(I,\pr)_\A(\K),\fM(I,\pr,\ka)_\A(\K)$. Then $\M(I,\pr)_\A$ and
$\M(I,\pr,\ka)_\A$ are the {\it sets of isomorphism classes of\/
$(I,\pr)$- and\/ $(I,\pr,\ka)$-configurations in} $\A$,
by~\cite[Prop.~7.6]{Joyc3}.
\end{itemize}
Each stabilizer group $\Iso_\K([X])$ or $\Iso_\K\bigl([(\si,
\io,\pi)]\bigr)$ in $\fObj_\A$ or $\fM(I,\pr)_\A$ is the group of
invertible elements in the finite-dimensional $\K$-algebra $\End(X)$
or $\End\bigl((\si,\io,\pi)\bigr)$. Thus $\fObj_\A,
\fObj_\A^\al,\fM(I,\pr)_\A,\fM(I,\pr,\ka)_\A$ have {\it affine
geometric stabilizers}, which is required to use the results
of~\S\ref{ai22}.

In \cite[Def.~7.7 \& Prop.~7.8]{Joyc3} we define 1-morphisms of
$\K$-stacks, as follows:
\begin{itemize}
\setlength{\itemsep}{0pt}
\setlength{\parsep}{0pt}
\item For $(I,\pr)$ a finite poset, $\ka:I\ra\bar C(\A)$ and
$J\in\F_\sIp$, we define $\bs\si(J):\fM(I,\pr)_\A\ra\fObj_\A$
or $\bs\si(J):\fM(I,\pr,\ka)_\A\ra\fObj_\A^{
\sum_{j\in J}\ka(j)}$. The induced maps $\bs\si(J)_*:\M(I,\pr)_\A
\ra\fObj_\A(\K)$ or $\M(I,\pr,\ka)_\A\ra\fObj_\A^{\sum_{j\in J}\ka(j)}
(\K)$ act by~$\bs\si(J)_*:[(\si,\io,\pi)]\mapsto[\si(J)]$.
\item For $(I,\pr)$ a finite poset, $\ka:I\ra\bar C(\A)$ and
$J\in\F_\sIp$, we define the $(J,\pr)$-{\it subconfiguration
$1$-morphism} $S(I,\pr,J):\fM(I,\pr)_\A\ra\fM(J,\pr)_\A$ or
$S(I,\pr,J):\fM(I,\pr,\ka)_\A\ra\fM(J,\pr,\ka\vert_J)_\A$.
The induced maps $S(I,\pr,J)_*$ act by $S(I,\pr,J)_*:
[(\si,\io,\pi)]\mapsto[(\si',\io',\pi')]$, where
$(\si,\io,\pi)$ is an $(I,\pr)$-configuration in $\A$,
and $(\si',\io',\pi')$ its $(J,\pr)$-subconfiguration.
\item Let $(I,\pr)$, $(K,\tl)$ be finite posets, $\ka:I\ra\bar C(\A)$,
and $\phi:I\ra K$ be surjective with $i\pr j$ implies $\phi(i)
\tl\phi(j)$ for $i,j\in I$. Define $\mu:K\ra\bar C(\A)$ by $\mu(k)=
\sum_{i\in\phi^{-1}(k)}\ka(i)$. We define the {\it quotient\/
$(K,\tl)$-configuration\/ $1$-morphisms}
\ea
&Q(I,\pr,K,\tl,\phi):\fM(I,\pr)_\A\ra\fM(K,\tl)_\A,
\label{ai3eq1}\\
&Q(I,\pr,K,\tl,\phi):\fM(I,\pr,\ka)_\A\ra\fM(K,\tl,\mu)_\A.
\label{ai3eq2}
\ea
The induced maps $Q(I,\pr,K,\tl,\phi)_*$ act by
$Q(I,\pr,K,\tl,\phi)_*:[(\si,\io,\pi)]\mapsto[(\ti\si,\ti\io,\ti\pi)]$,
where $(\si,\io,\pi)$ is an $(I,\pr)$-configuration in $\A$, and
$(\ti\si,\ti\io,\ti\pi)$ its quotient $(K,\tl)$-configuration
from~$\phi$.
\end{itemize}

Here \cite[Th.~8.4]{Joyc3} are some properties of these 1-morphisms:

\begin{thm} In the situation above:
\begin{itemize}
\setlength{\itemsep}{0pt}
\setlength{\parsep}{0pt}
\item[{\rm(a)}] $Q(I,\pr,K,\tl,\phi)$ in \eq{ai3eq1} and\/ \eq{ai3eq2}
are representable, and\/ \eq{ai3eq2} is finite type.
\item[{\rm(b)}] $\bs\si(I):\fM(I,\pr,\ka)_\A\longra\fObj_\A^{\ka(I)}$ is
representable and of finite type, and\/ $\bs\si(I):\fM(I,\pr)_\A\longra
\fObj_\A$ is representable.
\item[{\rm(c)}] $\prod_{i\in I}\bs\si(\{i\}):\fM(I,\pr)_\A\longra
\prod_{i\in I}\fObj_\A$ is of finite type.
\end{itemize}
\label{ai3thm}
\end{thm}

In \cite[\S 9--\S 10]{Joyc3} we define the data $\A,K(\A),\fF_\A$ in
some large classes of examples, and prove that Assumption
\ref{ai3ass} holds in each case.

\section{Algebras of constructible functions on $\fObj_\A$}
\label{ai4}

We now generalize the idea of Ringel--Hall algebras to
{\it constructible functions on stacks}. Let Assumption
\ref{ai3ass} hold. We show how to make the $\Q$-vector
space $\CF(\fObj_\A)$ of constructible functions into
a $\Q$-algebra.

We begin in \S\ref{ai41} by defining the multiplication $*$ on
$\CF(\fObj_\A)$, and showing it is associative. Section \ref{ai42}
extends this to locally constructible functions, and \S\ref{ai43}
constructs left or right representations of the algebras
$\CF(\fObj_\A)$, using configuration moduli stacks. Section
\ref{ai44} shows the subspace $\CFi(\fObj_\A)$ of functions {\it
supported on indecomposables} is a Lie subalgebra of
$\CF(\fObj_\A)$, and \S\ref{ai45} that under extra conditions on
$\A$ the subspace $\CF_\fin(\fObj_\A)$ of functions with {\it finite
support\/} is a subalgebra of~$\CF(\fObj_\A)$.

Section \ref{ai46} proves the $\Q$-algebra $\CF_\fin(\fObj_\A)$
is isomorphic to the {\it universal enveloping algebra\/}
$U\bigl(\CF_\fin^\ind(\fObj_\A)\bigr)$. Section \ref{ai47}
defines a commutative comultiplication $\De$ on $\CF_\fin
(\fObj_\A)$ making it into a bialgebra, and \S\ref{ai48} defines
multilinear operations $P_\sIp$ on $\CF(\fObj_\A)$ for all finite
posets $(I,\pr)$, which satisfy an analogue of associativity.
Finally, \S\ref{ai49} gives some examples from quiver
representations~$\modKQ$.

Throughout this section we fix an algebraically closed field $\K$ of
{\it characteristic zero}, so that we can apply the constructible
functions theory of~\S\ref{ai22}.

\subsection{An associative algebra structure on $\CF(\fObj_\A)$}
\label{ai41}

We now extend the Ringel--Hall algebra idea to the stacks
set up of \S\ref{ai32}. First we define the {\it identity}
$\de_{[0]}$ and {\it multiplication} $*$ on~$\CF(\fObj_\A)$.

\begin{dfn} Suppose Assumption \ref{ai3ass} holds. Define
$\de_{[0]}\in\CF(\fObj_\A)$ to be the characteristic
function of the point $[0]\in\fObj_\A(\K)$, so that
$\de_{[0]}\bigl([X]\bigr)=1$ if $X\cong 0$, and
$\de_{[0]}\bigl([X]\bigr)=0$ otherwise. For $f,g\in\CF
(\fObj_\A)$ we define $f\ot g\in\CF(\fObj_\A\t\fObj_\A)$
by $(f\ot g)\bigl([X],[Y]\bigr)=f([X])g([Y])$ for all~$\bigl([X],
[Y]\bigr)\in(\fObj_\A\t\fObj_\A)(\K)=\fObj_\A(\K)\t\fObj_\A(\K)$.

Using the diagrams of 1-morphisms of stacks and pullbacks,
pushforwards of constructible functions
\begin{equation*}
\text{
\begin{footnotesize}
$\displaystyle
\begin{gathered}
\xymatrix@C=54pt@R=15pt{
\fObj_\A\t\fObj_\A &
\fM(\{1,2\},\le)_\A
\ar[l]_{\bs\si(\{1\})\t\bs\si(\{2\})}
\ar[r]^{\bs\si(\{1,2\})}
& \fObj_\A,\\
\CF\bigl(\fObj_\A\bigr)\t\CF\bigl(\fObj_\A\bigr)
\!\!\!\!\!\!\!\!\!\!\!\!\!\!\!
\ar@<.5ex>[dr]^(0.6){\qquad\qquad(\bs\si(\{1\}))^*\cdot(\bs\si(\{2\}))^*}
\ar[d]_{\ot}
\\
\CF\bigl(\fObj_\A\!\t\!\fObj_\A\bigr)
\ar[r]^(.47){(\bs\si(\{1\})\t\bs\si(\{2\}))^*}
& \CF\bigl(\fM(\{1,2\},\le)_\A\bigr)
\ar[r]^(.56){\CF^\stk(\bs\si(\{1,2\}))}
& \CF\bigl(\fObj_\A\bigr),
}
\end{gathered}
$
\end{footnotesize}}
\end{equation*}
define a bilinear operation $*:\CF(\fObj_\A)\t\CF(\fObj_\A)
\ra\CF(\fObj_\A)$ by
\e
\begin{split}
f*g&=\CF^\stk(\bs\si(\{1,2\}))\bigl[\bs\si(\{1\})^*(f)\cdot
\bs\si(\{2\})^*(g)\bigr]\\
&=\CF^\stk(\bs\si(\{1,2\}))\bigl[(\bs\si(\{1\})\t\bs\si(\{2\}))^*
(f\ot g)\bigr].
\end{split}
\label{ai4eq1}
\e
This is well-defined as $\bs\si(\{1,2\})$ is representable and
$\bs\si(\{1\})\t\bs\si(\{2\})$ of finite type by Theorem
\ref{ai3thm}(b),(c), so $\CF^\stk(\bs\si(\{1,2\}))$ and
$(\bs\si(\{1\})\t\bs\si(\{2\}))^*$ are well-defined maps of
constructible functions as in~\S\ref{ai22}.
\label{ai4def1}
\end{dfn}

\begin{rem} Our convention on the order of multiplication,
with $(f*g)([X])$ an integral of $f(Y)g(X/Y)$ over subobjects
$Y\subset X$, agrees with Frenkel et al.\ \cite{FMV} and
Riedtmann \cite[\S 2]{Ried}. However, Lusztig \cite[\S\S 3.1,
10.19, 12.10]{Lusz} and Ringel use the opposite convention,
with $(f*g)([X])$ an integral of~$f(X/Y)g(Y)$.
\label{ai4rem1}
\end{rem}

Here is the basic result, saying $\CF(\fObj_\A)$ is a
$\Q$-algebra. The proof is related to Ringel \cite{Ring2}
and Lusztig~\cite[\S 10.19]{Lusz}.

\begin{thm} In the situation above, $\de_{[0]}*f=f*\de_{[0]}=f$
for all\/ $f$ in $\CF(\fObj_\A),$ and\/ $*$ is associative. Thus\/
$\CF(\fObj_\A)$ is a $\Q$-algebra, with identity\/ $\de_{[0]}$ and
multiplication\/~$*$.
\label{ai4thm1}
\end{thm}

\begin{proof} By considering the $(\{1,2\},\le)$-configuration
$(\si,\io,\pi)$ with $\si(\{1\})=0$ and $\si(\{2\})=X$ for
$X\in\A$, which is unique up to isomorphism, we find from
Definition \ref{ai4def1} that $(\de_{[0]}*f)\bigl([X]\bigr)
=f\bigl([X]\bigr)$, so $\de_{[0]}*f=f$ as this holds for all
$X\in\A$, and similarly~$f*\de_{[0]}=f$.

Define $\al,\be:\{1,2,3\}\ra\{1,2\}$ by $\al(1)\!=\!\al(2)\!=\!1$,
$\al(3)\!=\!2$, $\be(1)\!=\!1$, $\be(2)\!=\!\be(3)\!=\!2$. Consider
the commutative diagram of 1-morphisms, and the corresponding diagram
of pullbacks and pushforwards:
\begin{gather}
\begin{gathered}
\hskip -.3in
\xymatrix@C=25pt@R=30pt{
\fObj_\A\!\t\!\fObj_\A\!\t\!\fObj_\A
&
\fObj_\A\!\t\fM(\{2,3\},\le)_\A
\ar[l]^(0.4){\substack{\id_{\fObj_\A}\t\\
\bs\si(\{2\})\t\bs\si(\{3\})}}
\ar[r]
\ar@<-1ex>@{}[r]_(0.6){\id_{\fObj_\A}\t\bs\si(\{2,3\})}
& \fObj_\A\!\t\!\fObj_\A
\\
\fM(\{1,2\},\le)_\A\!\t\!\fObj_\A
\ar@<5ex>[u]_(0.4){\substack{\bs\si(\{1\})\t\bs\si(\{2\})\\ \t\id_{\fObj_\A}}}
\ar@<-5ex>[d]^(0.4){\bs\si(\{1,2\})\t\id_{\fObj_\A}}
& \fM(\{1,2,3\},\le)_\A
\ar[l]
\ar@<1ex>@{}[l]^(0.4){\substack{S(\{1,2,3\},\le,\{1,2\})\\ \t\bs\si(\{3\})}}
\ar[d]^{Q(\{1,2,3\},\le,\{1,2\},\le,\al)}
\ar[r]
\ar@<-1ex>@{}[r]_{Q(\{1,2,3\},\le,\{1,2\},\le,\be)}
\ar[u]_(0.4){\substack{\bs\si(\{1\})\t\\
S(\{1,2,3\},\le,\{2,3\})}}
&\fM(\{1,2\},\le)_\A
\ar@<-2ex>[u]^(0.4){\substack{\bs\si(\{1\})\t\\
\bs\si(\{2\})}}
\ar@<2ex>[d]_(0.7){\bs\si(\{1,2\})}
\\
\fObj_\A\!\t\!\fObj_\A
&\fM(\{1,2\},\le)_\A
\ar[l]_{\bs\si(\{1\})\t\bs\si(\{2\})}
\ar[r]^(0.6)
{\bs\si(\{1,2\})}
& \fObj_\A,
}
\hskip -.2in
\end{gathered}
\label{ai4eq2}
\allowdisplaybreaks \\[5pt]
\hskip -.5in
\text{
\begin{footnotesize}
$\displaystyle
\begin{gathered}
\xymatrix@C=10pt@R=30pt{
{}\quad\CF\bigl(\fObj_\A\!\t\!\fObj_\A\!\t\!\fObj_\A\bigr)
\ar@<-9ex>[d]^{\substack{(\bs\si(\{1\})\!\t\!\bs\si(\{2\})\\
\t\id_{\fObj_\A})^*}}
\ar[r]_(0.55){\substack{(\id_{\fObj_\A}\t\\
\bs\si(\{2\})\t\bs\si(\{3\}))^*}}
&\CF\bigl(\fObj_\A\!\t\fM(\{2,3\},\le)_\A\bigr)
\ar[r]
\ar@<-1ex>@{}[r]_(0.6){\CF^\stk(\id_{\fObj_\A}\t\bs\si(\{2,3\}))}
\ar[d]^(0.6){\substack{(\bs\si(\{1\})\t\\
S(\{1,2,3\},\le,\{2,3\}))^*}}
& \CF\bigl(\fObj_\A\!\t\!\fObj_\A\bigr)
\ar@<4ex>[d]_(0.6){\substack{(\bs\si(\{1\})\t\\
\bs\si(\{2\}))^*}}
\\
{}\qquad\CF\bigl(\fM(\{1,2\},\le)_\A\!\t\!\fObj_\A\bigr)
\ar@<-9ex>[d]^(0.35){\CF^\stk(\bs\si(\{1,2\})\t\id_{\fObj_\A})}
\ar[r]^(0.58){\substack{(S(\{1,2,3\},\le,\{1,2\})\\ \t\bs\si(\{3\}))^*}}
&
\CF\bigl(\fM(\{1,2,3\},\le)_\A\bigr)
\ar[d]_(0.6){\CF^\stk(Q(\{1,2,3\},\le,\{1,2\},\le,\al))}
\ar[r]
\ar@<-1ex>@{}[r]_(0.6){\CF^\stk(Q(\{1,2,3\},\le,\{1,2\},\le,\be))}
&\CF\bigl(\fM(\{1,2\},\le)_\A\bigr)\qquad{}
\ar@<4ex>[d]_(0.55){\CF^\stk(\bs\si(\{1,2\}))}
\\
\CF\bigl(\fObj_\A\!\t\!\fObj_\A\bigr)
\ar[r]^{(\bs\si(\{1\})\t\bs\si(\{2\}))^*}
& \CF\bigl(\fM(\{1,2\},\le)_\A\bigr)
\ar[r]^(0.7){\CF^\stk(\bs\si(\{1,2\}))}
&
\CF\bigl(\fObj_\A\bigr).
\!\!\!\!\!\!\!\!\!\!\!\!\!\!\!\!\!\!\!\!\!
}
\end{gathered}
$
\end{footnotesize}
}
\hskip -.4in
\label{ai4eq3}
\end{gather}
To show the maps in \eq{ai4eq3} are well-defined we use
Theorem \ref{ai3thm}(a),(b),(c) to show the corresponding
1-morphisms are representable or finite type, and
note that $S(\{1,2,3\},\le,\{1,2\})\t\bs\si(\{3\})$,
$\bs\si(\{1\})\t S(\{1,2,3\},\le,\{2,3\})$ are finite
type by a similar proof to Theorem~\ref{ai3thm}(c).

The top left square in \eq{ai4eq3} commutes by \eq{ai2eq3}, and the
bottom right by \eq{ai2eq2}. Now \cite[Th.~7.10]{Joyc3} implies that
\begin{equation*}
\xymatrix@C=190pt@R=10pt{
*+[r]{\fM(\{1,2\},\le)_\A}
\ar[d]^{\,\bs\si(\{1,2\})}
& *+[l]{\fM(\{1,2,3\},\le)_\A}
\ar[l]^{S(\{1,2,3\},\le,\{1,2\})}
\ar[d]_{Q(\{1,2,3\},\le,\{1,2\},\le,\al)\,}
\\
*+[r]{\fObj_\A} & *+[l]{\fM(\{1,2\},\le)_\A}
\ar[l]_{\bs\si(\{1\})}
}
\end{equation*}
is a Cartesian square. Taking fibre products with
$\bs\si(\{2\}):\fM(\{1,2\})_\A\ra\fObj_\A$ then shows the bottom
left square in \eq{ai4eq2} is also Cartesian, so the bottom left
square in \eq{ai4eq3} commutes by \eq{ai2eq4}. Similarly the top
right square in \eq{ai4eq3} commutes. Therefore \eq{ai4eq3}
commutes. Now let $f,g,h\in\CF(\fObj_\A)$, so that $f\ot g\ot h\in
\CF(\fObj_\A\t\fObj_\A\t\fObj_\A)$. As \eq{ai4eq3} commutes,
applying the two routes round the outside of the square to $f\ot
g\ot h$ shows that $(f*g)*h=f*(g*h)$. Thus $*$ is associative, and
$\CF(\fObj_\A)$ is an algebra.
\end{proof}

Because of the use of Cartesian squares and Theorem \ref{ai2thm1} in
the proof of Theorem \ref{ai4thm1}, to make $*$ associative we must
use the {\it stack pushforward\/} $\CF^\stk$ in \eq{ai4eq1}, and
other pushforwards such as $\CF^\na$ will in general give
nonassociative multiplications. In particular, as $\CF^\stk$ depends
on the stabilizer groups $\Iso_\K(x)$, we cannot afford to forget
this information by passing to coarse moduli schemes, if they
existed. This is an important reason for working with Artin stacks,
rather than some simpler class of spaces.

Define the {\it composition algebra} $\cC\!\subseteq\!\CF(\fObj_\A)$
to be the $\Q$-subalgebra generated by functions $f$ supported on
$[X]\in\fObj_\A(\K)$ with $X$ a {\it simple} object in $\A$. In
examples, the composition algebra $\cC$ is usually more interesting
than the Ringel--Hall algebra $\CF(\fObj_\A)$. When there are
only finitely many simple objects up to isomorphism $\cC$ is
{\it finitely generated}.

\subsection{Extension to locally constructible functions}
\label{ai42}

Next we observe that the associative multiplication $*$ in
\S\ref{ai41} extends to a large subspace $\dLCF(\fObj_\A)$
of the {\it locally constructible functions}~$\LCF(\fObj_\A)$.

\begin{dfn} Suppose Assumption \ref{ai3ass} holds. Define
$\dLCF(\fObj_\A)$ to be the $\Q$-vector subspace of $\LCF
(\fObj_\A)$ consisting of functions $f$ supported on subsets
$\coprod_{\al\in S}\fObj_\A^\al(\K)$ in $\fObj_\A(\K)$ for
$S\subset\bar C(\A)$ a {\it finite} subset. Following
\eq{ai4eq1}, define $*:\dLCF(\fObj_\A)\t\dLCF(\fObj_\A)
\ra\dLCF(\fObj_\A)$ by
\e
f*g=\LCF^\stk(\bs\si(\{1,2\}))\bigl[\bs\si(\{1\})^*(f)\cdot
\bs\si(\{2\})^*(g)\bigr].
\label{ai4eq4}
\e

To see this is well-defined, recall the disjoint union of stacks
\cite[Th.~7.5]{Joyc3}
\begin{equation*}
\fM(\{1,2\},\le)_\A=\ts\coprod_{\ka:\{1,2\}\ra\bar C(\A)}
\fM(\{1,2\},\le,\ka)_\A.
\end{equation*}
Let $f,g\in\dLCF(\fObj_\A)$ be supported on $\coprod_{\al\in S}
\fObj_\A^\al(\K)$, $\coprod_{\al\in T}\fObj_\A^\al(\K)$ respectively
for finite $S,T\subseteq\bar C(\A)$. Then $\bs\si(\{1\})^*(f)\cdot
\bs\si(\{2\})^*(g)$ in \eq{ai4eq4} is locally constructible and
supported on the finite number of $\M(\{1,2\},\le,\ka)_\A$ for
which $\ka(1)\in S$ and $\ka(2)\in T$. By Theorem \ref{ai3thm}(b)
$\bs\si(\{1,2\}):\fM(\{1,2\},\le,\ka)_\A\ra\fObj_\A^{\ka(\{1,2\})}$
is representable and of {\it finite type}. Thus
$\LCF^\stk(\bs\si(\{1,2\}))[\cdots]$ is well-defined in
\eq{ai4eq4}, and lies in $\LCF(\fObj_\A)$. But clearly
$f*g$ is supported on $\bigcup_{\al\in S,\;\be\in T}
\fObj_\A^{\al+\be}(\K)$, so $f*g\in\dLCF(\fObj_\A)$.
\label{ai4def2}
\end{dfn}

Actually, $*$ often makes sense on even larger subspaces of
$\LCF(\fObj_\A)$. For $f*g$ to be well-defined, all we need is that
for each $\ga\in\bar C(\A)$, there should exist only finitely many
pairs $\al,\be\in\bar C(\A)$ with $\ga=\al+\be$ and
$f\vert_{\fObj_\A^\al(\K)}$, $g\vert_{\fObj_\A^\be(\K)}$ both
nonzero. If it happens that for all $\ga\in\bar C(\A)$ there are
only finitely many pairs $\al,\be\in\bar C(\A)$ with $\al+\be=\ga$
then this holds automatically, and $f*g$ is well-defined for all
$f,g\in\LCF(\fObj_\A)$. In particular, this holds for all the quiver
examples of \cite[\S 10]{Joyc3}, so in these examples
$\LCF(\fObj_\A)$ will be a $\Q$-algebra.

We shall deal only with $\dLCF(\fObj_\A)$, though, as it is
sufficient for the applications in \cite{Joyc4,Joyc5}, where it is
useful, for instance, that $\de_{\fObj_\A^\al(\K)}$ lies in
$\dLCF(\fObj_\A)$. Here is the analogue of Theorem \ref{ai4thm1}.
The proof follows that of Theorem \ref{ai4thm1}, replacing
$\CF(\cdots)$ by $\LCF(\cdots)$ and $\CF^\stk(\cdots)$ by
$\LCF^\stk(\cdots)$, and arguing as in Definition \ref{ai4def2} to
show the operators $\LCF^\stk(\cdots)$ are well-defined.

\begin{thm} In the situation above $\dLCF(\fObj_\A)$ is a
$\Q$-algebra, with identity\/ $\de_{[0]}$ and associative
multiplication\/ $*$, and\/ $\CF(\fObj_\A)$ is a $\Q$-subalgebra.
\label{ai4thm2}
\end{thm}

In the rest of the section we give many results for $\CF(\fObj_\A)$.
Mostly these have straightforward generalizations to $\dLCF(\fObj_\A)$,
which we leave as exercises for the reader, just making the occasional
comment. Here are two other remarks:
\begin{itemize}
\setlength{\itemsep}{0pt}
\setlength{\parsep}{0pt}
\item The $\Q$-subalgebra of $\dLCF(\fObj_\A)$ generated by the
characteristic functions $\de_{\fObj_\A^\al(\K)}$ of $\fObj_\A^\al(\K)$
for $\al\in\bar C(\A)$ may be an interesting algebra.
\item One can also consider {\it infinite sums} in $\dLCF(\fObj_\A)$
or $\LCF(\fObj_\A)$. We call an infinite sum $\sum_{i\in I}f_i$ with
$f_i\in\LCF(\fObj_\A)$ {\it convergent\/} if for all constructible
$C\subseteq\fObj_\A(\K)$, only finitely many $f_i\vert_C$ are
nonzero. Then $\sum_{i\in I}f_i$ makes sense, and lies in
$\LCF(\fObj_\A)$. In \cite{Joyc5} we will prove identities which are
convergent infinite sums of products in~$\dLCF(\fObj_\A)$.
\end{itemize}

\subsection{Representations of Ringel--Hall algebras}
\label{ai43}

Here is a way to construct representations of the
algebra $\CF(\fObj_\A)$ of \S\ref{ai41}. Although it is
very simple, I did not find this method used explicitly
or implicitly in the Ringel--Hall algebra literature.

\begin{dfn} Let Assumption \ref{ai3ass} hold. Define
$\al,\be:\{1,2,3\}\ra\{1,2\}$ by
\e
\begin{gathered}
\al(1)=\al(2)=1,\;\> \al(3)=2,\;\> \be(1)=1,\;\> \be(2)=\be(3)=2.
\end{gathered}
\label{ai4eq5}
\e
Using the diagram of pullbacks, pushforwards of constructible functions
\begin{equation*}
\text{
\begin{footnotesize}
$\displaystyle
\xymatrix@C=54pt@R=15pt{
\CF\bigl(\fObj_\A\bigr)\!\t\!\CF\bigl(\fM(\{1,2\},\le)_\A\bigr)
\ar[dr]_(0.2){\begin{subarray}{l}
(\bs\si(\{2\}))^*\cdot\\
(Q(\{1,2,3\},\le,\{1,2\},\le,\al))^*\!\!\!\end{subarray}}
\ar@<-10ex>[d]_{\ot}
& \CF\bigl(\fM(\{1,2\},\le)_\A\bigr) \\
\CF\bigl(\fObj_\A\!\t\fM(\{1,2\},\le)_\A\bigr)
\ar[r]
\ar@<-1ex>@{}[r]_(0.6){(\bs\si(\{2\})\t Q(\{1,2,3\},\le,\{1,2\},\le,\al))^*}
& \CF\bigl(\fM(\{1,2,3\},\le)_\A\bigr)
\ar[u]_{\substack{\CF^\stk(Q(\{1,2,3\},\le,\\
{}\qquad\qquad\{1,2\},\le,\be))}}
}
$
\end{footnotesize}}
\end{equation*}
define $*_L:\CF(\fObj_\A)\t\CF(\fM(\{1,2\},\le)_\A)
\ra\CF(\fM(\{1,2\},\le)_\A)$ by
\e
\begin{split}
f*_Lr&=\CF^\stk(Q(\{1,2,3\},\le,\{1,2\},\le,\be))\\
&\qquad \bigl[\bs\si(\{2\})^*(f)\cdot
(Q(\{1,2,3\},\le,\{1,2\},\le,\al))^*(r)\bigr].
\end{split}
\label{ai4eq6}
\e
This is well-defined as $Q(\{1,2,3\},\le,\{1,2\},\le,\be)$ is
representable by Theorem \ref{ai3thm}(a), and one can show
$\bs\si(\{2\})\t Q(\cdots,\al)$ is finite type. In the same way,
define $*_R:\CF\bigl(\fM(\{1,2\},\le)_\A\bigr)\!\t\!
\CF(\fObj_\A)\!\ra\!\CF\bigl(\fM(\{1,2\},\le)_\A
\bigr)$ by
\begin{align*}
r*_Rf&=\CF^\stk(Q(\{1,2,3\},\le,\{1,2\},\le,\al))\\
&\qquad \bigl[
(Q(\{1,2,3\},\le,\{1,2\},\le,\be))^*(r)
\cdot\bs\si(\{2\})^*(f)\bigr].
\end{align*}
For $X\in\A$, write $V^{\sst[X]}$ for the vector subspace of $f\in
\CF\bigl(\fM(\{1,2\},\le)_\A\bigr)$ supported on points
$[(\si,\io,\pi)]$ with~$\si(\{1,2\})\cong X$.
\label{ai4def3}
\end{dfn}

The proof of the next theorem is modelled on that of
Theorem \ref{ai4thm1}. Note that we do {\it not\/} claim
that $\CF\bigl(\fM(\{1,2\},\le)_\A\bigr),V^{\sst[X]}$
are {\it two-sided\/} representations of $\CF(\fObj_\A)$,
only that $*_L$ and $*_R$ {\it separately} define left and
right representations. That is, we do not claim that
$(f*_Lr)*_Rg=f*_L(r*_Rg)$ for $f,g\in\CF(\fObj_\A)$ and
$r\in\CF\bigl(\fM(\{1,2\},\le)_\A\bigr)$, and in general
this is false.

\begin{thm} Above, if\/ $f,g\in\CF(\fObj_\A)$ and\/
$r\in\CF\bigl(\fM(\{1,2\},\le)_\A\bigr)$ then
$\de_{[0]}*_Lr=r$, $(f*g)*_Lr=f*_L(g*_Lr)$ and\/
$r*_R\de_{[0]}=r$, $r*_R(f*g)=(r*_Rg)*_Rf$. Thus,
$*_L,*_R$ define left and right representations
of the algebra $\CF(\fObj_\A)$ on the vector
space~$\CF\bigl(\fM(\{1,2\},\le)_\A)$.

Furthermore, for $X\in\A$ the subspace $V^{\sst[X]}$ is
closed under both actions $*_L,*_R$ of\/ $\CF(\fObj_\A)$.
Hence $V^{\sst[X]}$ is a left and right representation
of\/~$\CF(\fObj_\A)$.
\label{ai4thm3}
\end{thm}

\begin{proof} The first part of Theorem \ref{ai4thm1} generalizes
easily to show that $\de_{[0]}*_Lr=r$. Define $\ga,\de:\{1,2,3,4\}
\ra\{1,2\}$ and $\ep,\ze,\eta:\{1,2,3,4\}\ra\{1,2,3\}$ by
\begin{gather*}
\ga\!:\!1,2,3\!\mapsto\!1,\; \ga\!:\!4\!\mapsto\!2,\;
\de\!:\!1\!\mapsto\!1,\; \de\!:\!2,3,4\!\mapsto\!2,\;
\ep\!:\!1,2\!\mapsto\!1,\; \ep\!:\!3\!\mapsto\!2,\\
\ep\!:\!4\!\mapsto\!3,\;
\ze\!:\!1\!\mapsto\!1,\; \ze\!:\!2,3\!\mapsto\!2,
\ze\!:\!4\!\mapsto\!3,\; \eta\!:\!1\!\mapsto\!1,\;
\eta\!:\!2\!\mapsto\!2,\; \eta\!:\!3,4\!\mapsto\!3,
\end{gather*}
and $\al,\be$ as \eq{ai4eq5}. Consider the diagram of 1-morphisms
\e
\hskip -.1in
\text{
\begin{small}
$\displaystyle
\begin{gathered}
\xymatrix@C=45pt{
{\begin{subarray}{l}\ts\fObj_\A\t\\
\ts\fObj_\A\t\\
\ts\fM(\{1,2\},\le)_\A\end{subarray}}
&
{\begin{subarray}{l}\ts\fObj_\A\t\\
\ts\fM(\{1,2,3\},\le)_\A\end{subarray}}
\ar[l]
\ar@<1ex>@{}[l]^(0.43){\begin{subarray}{l}
\sst\id_{\fObj_\A}\t\bs\si(\{2\})\t \\
\sst Q(\{1,2,3\},\le,\{1,2\},\le,\al)\end{subarray}}
\ar[r]
\ar@<-1ex>@{}[r]_(0.65){\begin{subarray}{l}\sst\id_{\fObj_\A}\t\\
\sst Q(\{1,2,3\},\le,\{1,2\},\le,\be)\end{subarray}}
&
{\begin{subarray}{l}\ts\fObj_\A\t\\
\ts\fM(\{1,2\},\le)_\A\end{subarray}}
\\
{\begin{subarray}{l}\ts\fM(\{2,3\},\le)_\A\t\\
\ts\fM(\{1,2\},\le)_\A\end{subarray}}
\ar@<3ex>[u]_(0.4){\begin{subarray}{l}(\bs\si(\{2\})\t\bs\si(\{3\}))\\
\t\id_{\fM(\{1,2\},\le)_\A}\end{subarray}}
\ar@<-3ex>[d]^(0.55){\bs\si(\{2,3\})\t\id_{\fM(\{1,2\},\le)_\A}}
& \fM(\{1,2,3,4\},\le)_\A
\ar[l]
\ar@<1ex>@{}[l]^(0.4){\substack{
\sst S(\{1,2,3,4\},\le,\{2,3\})\t \\
\sst Q(\{1,2,3,4\},\le,\{1,2\},\le,\ga)\\ {}}}
\ar[r]
\ar@<-1ex>@{}[r]_{\sst Q(\{1,2,3,4\},\le,\{1,2,3\},\le,\eta)}
\ar[u]_(0.4){\sst\bs\si(\{2\})\t Q(\{1,2,3,4\},\le,\{1,2,3\},\le,\ep)}
\ar[d]^(0.6){
\sst Q(\{1,2,3,4\},\le,\{1,2,3\},\le,\ze)}
&\fM(\{1,2,3\},\le)_\A
\ar@<-4ex>[u]^(0.2){\sst\bs\si(\{2\})\t Q(\{1,2,3\},\le,\{1,2\},\le,\al)}
\ar@<4ex>[d]_(0.4){\sst Q(\{1,2,3\},\le,\{1,2\},\le,\be)}
\\
{\begin{subarray}{l}\ts\fObj_\A\t\\
\ts\fM(\{1,2\},\le)_\A\end{subarray}}
&\fM(\{1,2,3\},\le)_\A
\ar[l]_(0.6){\substack{\sst\bs\si(\{2\})\t\\
\sst Q(\{1,2,3\},\le,\{1,2\},\le,\al)}}
\ar[r]
\ar@<1ex>@{}[r]^(0.65){\sst Q(\{1,2,3\},\le,\{1,2\},\le,\be)}
& \fM(\{1,2\},\le)_\A,\!\!\!\!\!\!\!\!\!\!\!\!
}
\end{gathered}
$
\end{small}
}
\hskip -.3in
\label{ai4eq7}
\e
analogous to \eq{ai4eq2}. It is not difficult to show
\eq{ai4eq7} commutes, and the top right and bottom
left squares are Cartesian. Thus there is a commutative
diagram of spaces $\CF(\cdots)$ and pullbacks/pushforwards
analogous to \eq{ai4eq3}. Applying this to $f\ot g\ot r$ in
$\CF\bigl(\fObj_\A\!\t\!\fObj_\A\!\t\fM(\{1,2\},\le)_\A\bigr)
{}$ gives $(f*g)*_Lr=f*_L(g*_Lr)$, and $*_L$ is a left
representation. The proof for $*_R$ is similar.

Finally, if $r\in V^{\sst[X]}$ then $r$ is supported on points
$[(\si,\io,\pi)]$ with $\si(\{1,2\})\cong X$. If we have a diagram
\begin{equation*}
\xymatrix@C=95pt{
[(\si,\io,\pi)] &
[(\si',\io',\pi')]
\ar@{|->}[l]_{Q(\{1,2,3\},\le,\{1,2\},\le,\al)_*}
\ar@{|->}[r]^{Q(\{1,2,3\},\le,\{1,2\},\le,\be)_*}
&[(\ti\si,\ti\io,\ti\pi)],
}
\end{equation*}
then $X\cong\si(\{1,2\})\cong\si'(\{1,2,3\})\cong\ti\si(\{1,2\})$
as $\al(\{1,2,3\})=\{1,2\}=\be(\{1,2,3\})$. Hence in \eq{ai4eq6},
$(Q(\{1,2,3\},\le,\{1,2\},\le,\al))^*(r)$ is supported on points
$[(\si',\io',\pi')]$ with $\si'(\{1,2,3\})\cong X$, and thus
$f*_Lr$ is supported on points $[(\ti\si,\ti\io,\ti\pi)]$ with
$\ti\si(\{1,2\})\cong X$. That is, $f*_Lr$ lies in $V^{\sst[X]}$, so
$V^{\sst[X]}$ is closed under $*_L$. The proof for $*_R$ is similar.
\end{proof}

This shows that the big representation $\CF\bigl(\fM(\{1,2\},\le)_\A)$
contains many smaller subrepresentations $V^{\sst[X]}$. In examples,
this may be a useful tool for constructing {\it finite-dimensional
representations} of interesting infinite-dimensional algebras, such
as universal enveloping algebras.

The ideas above extend easily to representations $*_L,*_R$ of
$\dLCF(\fObj_\A)$ in \S\ref{ai42} on $\dLCF(\fM(\{1,2\},\le)_\A)$,
the subspace of $f\in\LCF(\fM(\{1,2\},\le)_\A)$ supported on the
union of $\M(\{1,2\},\le,\ka)_\A$ over finitely many~$\ka:\{1,2\}
\ra\bar C(\A)$.

\subsection{Indecomposables and Lie algebras}
\label{ai44}

We now give an analogue of ideas of Ringel \cite{Ring1}
and Riedtmann \cite{Ried}, who both define a Lie algebra
structure on spaces of functions on isomorphism classes
of {\it indecomposable} objects in an abelian category.

\begin{dfn} Suppose Assumption \ref{ai3ass} holds. An
object $0\not\cong X\in\A$ is called {\it decomposable} if $X\cong
Y\op Z$ for $0\not\cong Y,Z\in\A$. Otherwise $X$ is {\it
indecomposable}. Write $\CFi(\fObj_\A)$ for the subspace of $f$ in
$\CF(\fObj_\A)$ {\it supported on indecomposables}, that is,
$f\bigl([X]\bigr)\ne 0$ implies $0\not\cong X$ is indecomposable.

Decomposability can be characterized in terms of the
finite-dimensional $\K$-algebra $\End(X)=\Hom(X,X)$:
$X$ is decomposable if and only if there exist
$0\ne e_1,e_2\in\End(X)$ with $1=e_1+e_2$, $e_1^2=e_1$,
$e_2^2=e_2$ and $e_1e_2=e_2e_1=0$. Then $e_1,e_2$ are called
{\it orthogonal idempotents}. Given such $e_1,e_2$ we can
define nonzero objects $Y=\Im e_1$ and $Z=\Im e_2$ in
$\A$, and there exists an isomorphism $X\cong Y\op Z$
identifying $e_1,e_2$ with $\id_Y,\id_Z$. By choosing a
set of {\it primitive orthogonal idempotents} in $\End(X)$
one can show that each $0\not\cong X\in\A$ may be written
$X\cong V_1\op\cdots\op V_n$ for indecomposable
$V_1,\ldots,V_n$, unique up to order and isomorphism.

Define a bilinear bracket $[\,,\,]:\CF(\fObj_\A)\t\CF(\fObj_\A)
\ra\CF(\fObj_\A)$ by
\e
[f,g]=f*g-g*f,
\label{ai4eq8}
\e
for $*$ defined in \eq{ai4eq1}. Since $*$ is associative by
Theorem \ref{ai4thm1}, $[\,,\,]$ satisfies the {\it Jacobi identity},
and makes $\CF(\fObj_\A)$ into a {\it Lie algebra} over~$\Q$.
\label{ai4def4}
\end{dfn}

The following result is related to Riedtmann \cite[\S 2]{Ried},
and~\cite[Prop.~2.2.8]{FMV}.

\begin{thm} In the situation above, $\CFi(\fObj_\A)$ is
closed under the Lie bracket\/ $[\,,\,],$ and is a Lie algebra
over\/~$\Q$.
\label{ai4thm4}
\end{thm}

\begin{proof} Let $f,g\in\CFi(\fObj_\A)$ and $Y\in\A$
with $(f*g)([Y])\ne 0$. By a long but elementary calculation
involving properties of the Euler characteristic, we can show
that either (i) $Y$ is indecomposable;
\begin{itemize}
\setlength{\itemsep}{0pt}
\setlength{\parsep}{0pt}
\item[(ii)] $Y\cong X\op Z$ for $X,Z\in\A$ indecomposable
with $X\not\cong Z$ and
\begin{align*}
(f*g)([Y])&=\bigl(f([X])g([Z])+f([Z])g([X])\bigr)\cdot\\
&\quad \chi\bigl(\Aut(X\op Z)/\Aut(X)\t\Aut(Z)\bigr);\text{ or}
\end{align*}
\item[(iii)] $Y\cong X\op X$ for $X\in\A$ indecomposable and
\end{itemize}
\begin{equation*}
(f*g)([Y])=f([X])g([X])\cdot\chi\bigl(\Aut(X\op X)/\Aut(X)\t\Aut(X)\bigr).
\end{equation*}
In (ii),(iii) we have $(f*g)([Y])=(g*f)([Y])$, so that
$[f,g]([Y])=0$ by \eq{ai4eq8}. Hence the only possibility in which
$[f,g]([Y])\ne 0$ is (i), when $Y$ is indecomposable. Thus
$[f,g]\in\CFi(\fObj_\A)$, as we have to prove.
\end{proof}

In the same way we find $\dLCF^\ind(\fObj_\A)$ is a Lie
subalgebra of~$\dLCF(\fObj_\A)$.

\subsection{Constructible functions with finite support}
\label{ai45}

Here is some more notation. As functions with finite support are
always constructible, we do not generalize this to~$\dLCF(\fObj_\A)$.

\begin{dfn} Write $\CF_\fin(\fObj_\A)$ for the subspace of
$f\in\CF(\fObj_\A)$ with {\it finite support}, that is, $f$ is
nonzero on only finitely many points in $\fObj_\A(\K)$. Define
$\CFi_\fin(\fObj_\A)=\CFi(\fObj_\A) \cap\CF_\fin(\fObj_\A)$. For
each $[X]\in\fObj_\A(\K)$, write $\de_{[X]}:\fObj_\A(\K)\ra\{0,1\}$
for the characteristic function of $[X]$. Then the
$\de_{\smash{[X]}}$ form a basis for $\CF_\fin(\fObj_\A)$, and the
$\de_{\smash{[X]}}$ for indecomposable $X$ form a basis
for~$\CF_\fin^\ind(\fObj_\A)$.
\label{ai4def5}
\end{dfn}

We want $\CF_\fin(\fObj_\A)$ to be a subalgebra of $\CF(\fObj_\A)$.
To prove this we need an extra assumption, which holds in Example
\ref{ai4ex1} below.

\begin{ass} For all $X,Z\!\in\!\A$, there are only finitely
isomorphism classes of $Y\!\in\!\A$ for which there exists an
exact sequence $0\!\ra\!X\!\ra\!Y\!\ra\!Z\!\ra\!0$ in~$\A$.
\label{ai4ass}
\end{ass}

If {\it all constructible sets in $\fObj_\A$ are finite} then
$\CF_\fin(\fObj_\A)\!=\!\CF(\fObj_\A)$ and $\CF_\fin^\ind
(\fObj_\A)\!=\!\CFi(\fObj_\A)$, and Assumption \ref{ai4ass} holds
automatically.

\begin{prop} If Assumptions \ref{ai3ass} and \ref{ai4ass} hold
then $\CF_\fin(\fObj_\A)$ is closed under $*$, and\/
$\CF_\fin^\ind(\fObj_\A)$ is closed under~$[\,,\,]$.
\label{ai4prop1}
\end{prop}

\begin{proof} For $X,Z\in\A$, $\de_{[X]}*\de_{[Z]}$ is supported
on the set of $[Y]\in\fObj_\A(\K)$ for which there exists an
exact sequence $0\ra X\ra Y\ra Z\ra 0$. By Assumption
\ref{ai4ass} this set is finite, so $\de_{[X]}*\de_{[Z]}$ lies
in $\CF_\fin(\fObj_\A)$. As the $\de_{[X]}$ form a basis
of $\CF_\fin(\fObj_\A)$, it is closed under~$*$.
\end{proof}

\subsection{Universal enveloping algebras}
\label{ai46}

We study the universal enveloping algebras of~$\CFi(\fObj_\A),
\CFi_\fin(\fObj_\A)$.

\begin{dfn} Let $\g$ be a Lie algebra over $\Q$. The {\it
universal enveloping algebra} $U(\g)$ is the $\Q$-algebra
generated by $\g$ with the relations $xy-yx=[x,y]$ for
all $x,y\in\g$. Multiplication in $U(\g)$ will be written
as juxtaposition, $(x,y)\mapsto xy$. Each Lie algebra
representation of $\g$ extends uniquely to an algebra
representation of $U(\g)$, so $U(\g)$ is a powerful tool for
studying the representation theory of $\g$. See Humphreys
\cite{Hump} for an introduction to these ideas.

In Theorem \ref{ai4thm4} the embedding of the Lie algebra
$\CFi(\fObj_\A)$ with bracket $[\,,\,]$ in the algebra
$\CF(\fObj_\A)$ with multiplication $*$ satisfying \eq{ai4eq8}
implies there is a unique $\Q$-algebra homomorphism
\e
\Phi:U\bigl(\CFi(\fObj_\A)\bigr)\ra\CF(\fObj_\A)
\label{ai4eq9}
\e
with $\Phi(1)=\de_{[0]}$ and $\Phi(f_1\cdots f_n)= f_1*\cdots*f_n$
for $f_1,\ldots,f_n\in\CFi(\fObj_\A)$. Similarly, in Proposition
\ref{ai4prop1} there is a unique homomorphism
\e
\Phi_\fin:U\bigl(\CFi_\fin(\fObj_\A)\bigr)\ra\CF_\fin (\fObj_\A).
\label{ai4eq10}
\e
\label{ai4def6}
\end{dfn}

The next two results are similar to Riedtmann~\cite[\S 3]{Ried}.

\begin{prop} $\Phi$ and\/ $\Phi_\fin$ above are injective. Hence,
the $\Q$-subalgebras of\/ $\CF(\fObj_\A),\CF_\fin(\fObj_\A)$
generated by $\CFi(\fObj_\A),\CFi_\fin(\fObj_\A)$ are isomorphic to
$U\bigl(\CFi(\fObj_\A)\bigr),\ab U\bigl(\CFi_\fin(\fObj_\A)\bigr)$
respectively.
\label{ai4prop2}
\end{prop}

\begin{proof} We do the case of $\Phi_\fin$ first, so suppose
Assumptions \ref{ai3ass} and \ref{ai4ass} hold. Let $V_1,V_2\in\A$
be indecomposable. Applying (i)--(iii) in the proof of Theorem
\ref{ai4thm4} to $f=\de_{[V_1]}$, $g=\de_{[V_2]}$ we find
that $\de_{[V_1]}*\de_{[V_2]}$ is supported on points
$[V_1\op V_2]$ and $[Y]$ for $Y\in\A$ indecomposable, and
$\bigl(\de_{[V_1]}*\de_{[V_2]}\bigr)\bigl([V_1\op V_2]\bigr)=
\chi\bigl(\Aut(V_1\op V_2)/\Aut(V_1)\t\Aut(V_2)\bigr)$.
That is, $\chi\bigl(\Aut(V_1\op V_2)/\Aut(V_1)\t\Aut(V_2)\bigr)
\cdot\de_{[V_1\op V_2]}-\de_{[V_1]}*\de_{[V_2]}$ is supported
on points $[Y]$ for indecomposable $Y$. It is not difficult to
generalize this to show that if $V_1,\ldots,V_m\in\A$ are
indecomposable then
\e
\begin{gathered}
\chi\bigl(\Aut(V_1\op\cdots\op V_m)/\Aut(V_1)\t\cdots\t\Aut(V_m)\bigr)
\cdot\de_{[V_1\op\cdots\op V_m]}-\\
\quad\text{$\de_{[V_1]}*\de_{[V_2]}*\cdots*\de_{[V_m]}$ is supported on
points $[W_1\op\cdots\op W_k]$}\\
\text{for indecomposable $W_1,\ldots,W_k\in\A$ and $1\le k<m$.}
\end{gathered}
\label{ai4eq11}
\e

Let the $V_1,\ldots,V_m$ have $a$ equivalence classes under isomorphism
with sizes $m_1,\ldots,m_a$, so that $m=m_1+\cdots+m_a$. Using facts
about the finite-dimensional algebras $\End(V_1\op\cdots\op V_m),\ab
\End(V_1),\ldots,\End(V_m)$ and the {\it Jacobson radical\/}
from Benson \cite[\S 1]{Bens} we find there is an isomorphism of
$\K$-varieties
\begin{equation*}
\Aut(V_1\!\op\!\cdots\!\op\!V_m)/\Aut(V_1)\!\t\!\cdots\!\t\!\Aut(V_m)
\!\cong\!\K^l\!\t\!\ts\prod_{i=1}^a\bigl(\GL(m_i,\K)/(\K^\t)^{m_i}\bigr),
\end{equation*}
which allows us to compute the Euler characteristic
\e
\chi\bigl(\Aut(V_1\op\cdots\op V_m)/\Aut(V_1)\t\cdots\t\Aut(V_m)\bigr)
=\ts\prod_{i=1}^am_i!.
\label{ai4eq12}
\e

To show $\Phi_\fin$ is injective, write $I=\bigl\{[X]
\in\fObj_\A(\K):X$ is indecomposable$\bigr\}$, and let $\le$ be any
arbitrary {\it total order} on $I$. Then $\bigl\{ \de_{[X]}:[X]\in
I\bigr\}$ is a basis for $\CFi_\fin (\fObj_\A)$, so the {\it
Poincar\'e--Birkhoff--Witt Theorem} \cite[Cor.~C, p.~92]{Hump} shows
\e
\{1\}\cup\bigl\{\de_{[X_1]}\de_{[X_2]}\cdots\de_{[X_n]}:n\ge 1,\;\>
[X_1],\ldots,[X_n]\in I,\;\> [X_1]\le\cdots\le[X_n]\bigr\}
\label{ai4eq13}
\e
is a basis for $U\bigl(\CFi_\fin(\fObj_\A)\bigr)$. Suppose $u\in
U\bigl(\CFi_\fin(\fObj_\A)\bigr)$ is nonzero. If $u=c\cdot 1$ for
$c\in\Q\sm\{0\}$ then $\Phi_\fin(u)=c\,\de_{[0]}\ne 0$. Otherwise,
there exists a basis element
$\de_{[X_1]}\de_{[X_2]}\cdots\de_{[X_n]}$ from \eq{ai4eq13} with $n$
greatest such that the coefficient $u_{[X_1]\cdots[X_n]}$ of
$\de_{[X_1]}\de_{[X_2]}\cdots\de_{[X_n]}$ is nonzero.

We shall evaluate the constructible function $\Phi_\fin(u)$ at
the point $[X_1\op\cdots\op X_n]$. Let $\de_{[V_1]}\cdots
\de_{[V_m]}$ be any basis element from \eq{ai4eq13} with
nonzero coefficient in $u$. Then $m\le n$, by choice of $n$. We have
\begin{align*}
\Phi_\fin(\de_{[V_1]}\cdots\de_{[V_m]})\bigl([X_1\!\op\!\cdots
\!\op\!X_n]\bigr)&=(\de_{[V_1]}*\cdots*\de_{[V_m]})\bigl([X_1
\op\cdots\op X_n]\bigr)\\
=c\,\de_{[V_1\op\cdots\op V_m]}\bigl([X_1\!\op\!\cdots\!\op\!X_n]\bigr)
&=\begin{cases}
c, & [V_1\!\op\!\cdots\!\op\!V_m]\!=\![X_1\!\op\!\cdots\!\op\!X_n], \\
0, & [V_1\!\op\!\cdots\!\op\!V_m]\!\ne\![X_1\!\op\!\cdots\!\op\!X_n],
\end{cases}
\end{align*}
where $c$ is the nonzero integer \eq{ai4eq12} and in the second
line we use \eq{ai4eq11} and $m\le n$ to see that $c\,\de_{[V_1
\op\cdots\op V_m]}-\de_{[V_1]}*\cdots*\de_{[V_m]}$ is zero
at~$[X_1\op\cdots\op X_n]$.

As $V_i,X_j$ are indecomposable, $[V_1\op\cdots\op V_m]=
[X_1\op\cdots\op X_n]$ if and only if $m=n$ and $[V_1],
\ldots,[V_m]$ and $[X_1],\ldots,[X_n]$ are the same up to a
permutation of $1,\ldots,n$. But by assumption
$[V_1]\le\cdots\le[V_m]$ and $[X_1]\le\cdots\le[X_n]$ in the total
order $\le$ on $I$, so $[V_1\op\cdots\op V_m]=[X_1\op\cdots\op X_n]$
only if $[V_i]=[X_i]$ for all $i$. Therefore
$\Phi_\fin(u)\bigl([X_1\op\cdots\op X_n]
\bigr)=c\,u_{[X_1]\cdots[X_n]}\ne 0$, as $\de_{[X_1]}\cdots
\de_{[X_n]}$ is the only basis element making a nonzero
contribution. Hence $\Phi_\fin(u)\ne 0$ for all $0\ne u\in
U\bigl(\CFi_\fin(\fObj_\A)\bigr)$, and $\Phi_\fin$ is injective.

Showing $\Phi$ is injective uses essentially the same ideas, but is
a little more tricky as we cannot choose a basis for
$\CFi(\fObj_\A)$ consisting of functions with {\it disjoint
support}. We leave it as an exercise for the reader. For the last
part, as $\Phi,\Phi_\fin$ are injective they are isomorphisms with
their images, which are the $\Q$-subalgebras generated
by~$\CFi(\fObj_\A),\CFi_\fin(\fObj_\A)$.
\end{proof}

We shall show $\Phi_\fin$ is an isomorphism. This will enable
us to identify the algebra $\CF_\fin(\fObj_\A)$ in examples.

\begin{prop} Let Assumptions \ref{ai3ass} and \ref{ai4ass} hold.
Then $\Phi_\fin$ in \eq{ai4eq10} is an isomorphism.
\label{ai4prop3}
\end{prop}

\begin{proof} As $\Phi_\fin$ is injective by Proposition
\ref{ai4prop3}, we need only show it is surjective. Suppose
by induction that for some $m\ge 1$, $\Im\Phi_\fin$ contains
$\de_{[V_1\op\cdots\op V_n]}$ whenever $1\le n<m$ and
$V_i\in\A$ are indecomposable. This is trivial for $m=1$.
Let $V_1,\ldots,V_m\in\A$ be indecomposable. Then by
\eq{ai4eq11}, $c\,\de_{[V_1\op\cdots\op V_m]}-\de_{[V_1]}*
\cdots*\de_{[V_m]}$ lies in the span of functions $\de_{[W_1\op
\cdots\op W_k]}$ for indecomposable $W_i\in\A$ and $1\le k<m$,
where $c$ is the nonzero integer~\eq{ai4eq12}.

By induction $\de_{[W_1\op\cdots\op W_k]}\in\Im\Phi_\fin$, so
$c\,\de_{[V_1\op\cdots\op V_m]}-\de_{[V_1]}*\cdots*\de_{[V_m]}$
lies in $\Im\Phi_\fin$. But $\de_{[V_1]}*\cdots*\de_{[V_m]}=
\Phi_\fin(\de_{[V_1]}\cdots\de_{[V_m]})\in\Im\Phi_\fin$, so
$\de_{[V_1\op\cdots\op V_m]}$ lies in $\Im\Phi_\fin$. Thus
the inductive hypothesis holds for $m+1$, and by induction
$\de_{[V_1\op\cdots\op V_n]}\in\Im\Phi_\fin$ whenever $n\ge 1$
and $V_i\in\A$ are indecomposable.

But each $0\not\cong X\in\A$ may be written $X\cong V_1\op
\cdots\op V_n$ for $V_i\in\A$ indecomposable, by Definition
\ref{ai4def4}, so $\de_{[X]}\in\Im\Phi_\fin$ from above. Also
$\de_{[0]}=\Phi_\fin(1)\in\Im\Phi_\fin$, so $\de_{[X]}\in
\Im\Phi_\fin$ for all $X\in\A$. As the $\de_{[X]}$ are a
basis for $\CF_\fin(\fObj_\A)$, this proves
$\Im\Phi_\fin=\CF_\fin(\fObj_\A)$, and $\Phi_\fin$ is surjective.
\end{proof}

We shall see in \S\ref{ai47} that with extra structures on
$\CF_\fin(\fObj_\A)$, $\Phi_\fin$ is actually an isomorphism
of $\Q$-bialgebras and of Hopf algebras.

\subsection{Comultiplication and bialgebras}
\label{ai47}

Next we explain how to define a {\it cocommutative
comultiplication} on the $\Q$-algebra $\CF_\fin(\fObj_\A)$,
making it into a {\it bialgebra}. Our treatment is based on
Ringel \cite{Ring4}, who defines a similar comultiplication
on degenerate Ringel--Hall algebras at $q=1$. For an
introduction to bialgebras, see Joseph~\cite[\S 1]{Jose}.

\begin{dfn} Let Assumptions \ref{ai3ass} and \ref{ai4ass} hold,
so that $\CF_\fin(\fObj_\A)$ is a $\Q$-algebra by
Proposition \ref{ai4prop1}. Let
\e
\Psi:\CF_\fin(\fObj_\A)\ot\CF_\fin(\fObj_\A)\ra
\CF_\fin(\fObj_\A\t\fObj_\A)
\label{ai4eq14}
\e
be the unique linear map with $\Psi(f\ot g)=f\ot g$
for $f,g\in\CF_\fin(\fObj_\A)$, in the notation of
Definition \ref{ai4def1}. Since $\de_{([X],[Y])}$ for
$X,Y\in\A$ form a basis for $\CF_\fin(\fObj_\A\t\fObj_\A)$
and $\Psi(\de_{[X]}\ot\de_{[Y]})=\de_{[X]}\ot\de_{[Y]}=\de_{([X],
[Y])}$, we see that $\Psi$ is an isomorphism of $\Q$-vector spaces.

For $I$ any finite set, let $\bu$ be the partial order on
$\{1,2\}$ with $i\bu j$ if and only if $i=j$. Consider
the diagram of 1-morphisms
\begin{equation*}
\xymatrix@C=70pt{
\fObj_\A &
\fM(\{1,2\},\bu)_\A
\ar[l]_{\bs\si(\{1,2\})}
\ar[r]^{\bs\si(\{1\})\t\bs\si(\{2\})}
& \fObj_\A\t\fObj_\A.
}
\end{equation*}
By \cite[Prop.~7.9]{Joyc3} $\bs\si(\{1\})\t\bs\si(\{2\})$ is a
1-isomorphism, and so is representable. Now
$\bs\si(\{1,2\})_*\ci(\bs\si(\{1\})\t\bs\si(\{2\}))_*^{-1}:
(\fObj_\A\t\fObj_\A)(\K)\ra\fObj_\A(\K)$ maps $([Y],[Z])
\mapsto[Y\op Z]$. Since for any $X\in\A$ there are only finitely
many pairs $Y,Z\in\A$ up to isomorphism with $Y\op Z\cong X$, and
$(\bs\si(\{1\})\t\bs\si(\{2\}))_*$ is a bijection, we see that
$\bs\si(\{1,2\})_*$ takes only finitely many points to each point in
$\fObj_\A(\K)$. Thus the following maps are well-defined:
\begin{equation*}
\xymatrix@C=22pt{
\CF_\fin(\fObj_\A)
\ar[r]
\ar@<1ex>@{}[r]^(0.4){(\bs\si(\{1,2\}))^*}
& \CF_\fin\bigl(\fM(\{1,2\},\bu)_\A\bigr)
\ar[rr]
\ar@{}@<1ex>[rr]^{\CF^\stk(\bs\si(\{1\})\!\t\bs\si(\{2\}))}
&& \CF_\fin(\fObj_\A\!\t\!\fObj_\A).
}
\end{equation*}
As $\Psi$ is an isomorphism $\Psi^{-1}$ exists, so we may
define the {\it comultiplication}
\e
\begin{split}
&\De:\CF_\fin(\fObj_\A)\ra\CF_\fin(\fObj_\A)
\ot\CF_\fin(\fObj_\A)\\
\text{by}\quad
&\De=\Psi^{-1}\ci\CF^\stk
(\bs\si(\{1\})\t\bs\si(\{2\}))\ci(\bs\si(\{1,2\}))^*.
\end{split}
\label{ai4eq15}
\e
Define the {\it counit\/} $\ep:\CF_\fin(\fObj_\A)\ra\Q$
by~$\ep:f\mapsto f([0])$.
\label{ai4def7}
\end{dfn}

\begin{thm} Let Assumptions \ref{ai3ass} and \ref{ai4ass} hold.
Then $\CF_\fin(\fObj_\A)$ with\/ $*,\De,\de_{[0]},\ep$ is
a cocommutative bialgebra.
\label{ai4thm5}
\end{thm}

\begin{proof} For the axioms of a bialgebra, see Joseph
\cite[\S 1.1]{Jose}. First we show $\CF_\fin(\fObj_\A),
\De,\ep$ form a cocommutative coalgebra. As $(\{1,2\},\bu)$
is preserved by exchanging 1,2, equation \eq{ai4eq15} is
unchanged by exchanging the factors in
$\CF_\fin(\fObj_\A)\ot\CF_\fin(\fObj_\A)$, so
$\De$ is {\it cocommutative}. To show $\De$ is
{\it coassociative} we must show the following commutes:
\begin{equation*}
\text{
\begin{small}
$\displaystyle
\xymatrix@C=160pt@R=10pt{
*+[r]{\CF_\fin(\fObj_\A)}
\ar[r]_\De
\ar[d]^{\;\De}
& *+[l]{\CF_\fin(\fObj_\A)\ot\CF_\fin(\fObj_\A)}
\ar[d]_{\De\ot\id\;}
\\
*+[r]{\CF_\fin(\fObj_\A)\!\ot\!\CF_\fin(\fObj_\A)}
\ar[r]^(0.4){\id\ot\De}
& *+[l]{\CF_\fin(\fObj_\A)
\!\ot\!\CF_\fin(\fObj_\A)
\!\ot\!\CF_\fin(\fObj_\A).}
}
$
\end{small}
}
\end{equation*}
Using \cite[Prop.~7.9]{Joyc3} and the fact $\Psi$ is an isomorphism,
this follows provided
\begin{equation*}
\text{
\begin{small}
$\displaystyle
\xymatrix@C=190pt@R=10pt{
*+[r]{\CF_\fin(\fObj_\A)}
\ar[r]_{\bs\si(\{1,2\})^*}
\ar[d]^{\;\bs\si(\{1,2\})^*}
& *+[l]{\CF_\fin\bigl(\fM(\{1,2\},\bu)_\A\bigr)}
\ar[d]_{Q(\{1,2,3\},\bu,\{1,2\},\bu,\al)^*\;}
\\
*+[r]{\CF_\fin\bigl(\fM(\{1,2\},\bu)_\A\bigr)}
\ar[r]^{Q(\{1,2,3\},\bu,\{1,2\},\bu,\be)^*}
& *+[l]{\CF_\fin\bigl(\fM(\{1,2,3\},\bu)_\A\bigr)}
}
$
\end{small}
}
\end{equation*}
commutes, where $\al,\be$ are as in \eq{ai4eq5}. But this
is immediate as $\bs\si(\{1,2\})\ci Q(\{1,2,3\},\bu,\{1,2\},
\bu,\al)=\bs\si(\{1,2\})\ci Q(\{1,2,3\},\bu,\{1,2\},\bu,\be)$.

Since $\CF_\fin(\fObj_\A)$ is cocommutative, to show $\ep$
is a {\it counit\/} we need
\begin{equation*}
\text{
\begin{small}
$\displaystyle
\xymatrix@R=10pt{
\CF_\fin(\fObj_\A)
\ar@<.5ex>[dr]^{\quad f\mapsto 1\ot f}
\ar[d]^{\;\De}
\\
\CF_\fin(\fObj_\A)\!\ot\!\CF_\fin(\fObj_\A)
\ar[r]^(0.55){\ep\ot\id}
& \Q\!\ot\!\CF_\fin(\fObj_\A)
}
$
\end{small}
}
\end{equation*}
to commute. This holds if $(\De f)\bigl([0],[X]\bigr)=
f\bigl([X]\bigr)$ for all $f\in\CF_\fin(\fObj_\A)$ and
$[X]\in\fObj_\A(\K)$, which is clear from \eq{ai4eq15}, as
there is just one point $[(\si,\io,\pi)]\in\M(\{1,2\},\bu)_\A$
with $\si(\{1\})\cong 0$ and $\si(\{2\})\cong X$, and it
has~$\si(\{1,2\})\cong X$.

Next we prove $\De$ is {\it multiplicative}, that is,
$\De$ is an algebra homomorphism from $\CF_\fin(\fObj_\A)$
with multiplication $*$ to $\CF_\fin(\fObj_\A)\!\ot\!\CF_\fin
(\fObj_\A)$ with multiplication $*\ot *$. Define $I=\{a,b,c,d\}$
and a partial order $\pr$ by $a\pr b$, $c\pr d$, and $i\pr i$ for
$i\in I$. Define maps $\mu,\nu:I\ra\{1,2\}$ by
\begin{equation*}
\mu(a)=\mu(b)=1,\; \mu(c)=\mu(d)=2,\;
\nu(a)=\nu(c)=1,\; \nu(b)=\nu(d)=2.
\end{equation*}
Then using \cite[Prop.~7.9]{Joyc3} and the fact that $\Psi$ is an
isomorphism, calculation shows $\De$ is multiplicative provided the
unbroken arrows `$\longra$' commute in
\e
{}\hskip -.6in
\text{
\begin{footnotesize}
$\displaystyle
\begin{gathered}
\xymatrix@!0@C=147pt@R=15pt{
*+[r]{\CF_\fin(\fObj_\A\!\t\!\fObj_\A)}
\ar[rr]_(0.4){(\bs\si(\{a,c\})\t\bs\si(\{b,d\}))^*}
\ar[dd]^{(\bs\si(\{1\})\t\bs\si(\{2\}))^*}
&&
*+[l]{\CF_\fin\bigl(\fM(\{a,c\},\bu)_\A\t\fM(\{b,d\},\bu)_\A\bigr)}
\ar[dd]_{(S(I,\pr,\{a,c\})\t S(I,\pr,\{b,d\}))^*}
\\
\\
*+[r]{\CF_\fin\bigl(\fM(\{1,2\},\le)_\A\bigr)}
\ar@<.5ex>[rr]^{Q(I,\pr,\{1,2\},\le,\nu)^*}
\ar@{-->}[dr]_(0.75){\pi_2^*}
\ar[dd]^{\CF^\stk(\bs\si(\{1,2\}))}
&&
*+[l]{\CF\text{ or }\CF_\fin\bigl(\fM(I,\pr)_\A\bigr)}
\ar[dd]_{\CF^\stk(Q(I,\pr,\{1,2\},\le,\mu))}
\\
&\CF(\fF)
\ar@{-->}[dr]_(0.1){\CF^\stk(\pi_1)}
\ar@{-->}[ur]^(0.2){\phi^*}
\\
*+[r]{\CF_\fin(\fObj_\A)}
\ar@<-.5ex>[rr]^(0.28){(\bs\si(\{1,2\}))^*}
&& *+[l]{\CF\text{ or }\CF_\fin\bigl(\fM(\{1,2\},\bu)_\A\bigr).}
}
\end{gathered}
$
\end{footnotesize}
}
\hskip -.4in{}
\label{ai4eq16}
\e
Here we write `$\CF$ or $\CF_\fin(\cdots)$' as the arrows
`$\longra$' map to $\CF_\fin(\cdots)$, but the arrows
`$\dashra$' defined below may map to~$\CF(\cdots)$.

The top square of \eq{ai4eq16} commutes as the corresponding
1-morphisms do. The bottom square is more tricky, since
although in
\e
\begin{gathered}
\xymatrix@!0@C=140pt@R=15pt{
*+[r]{\fM(\{1,2\},\le)_\A}
\ar[dd]^{\,\bs\si(\{1,2\})}
&& *+[l]{\fM(I,\pr)_\A}
\ar[dd]_{Q(I,\pr,\{1,2\},\le,\mu)}
\ar@<-.5ex>[ll]^(0.5){Q(I,\pr,\{1,2\},\le,\nu)}
\ar@{-->}[dl]_(0.7)\phi
\\
& \fF
\ar@{-->}[ul]^(0.2){\pi_2}
\ar@{-->}[dr]_(0.2){\pi_1}
\\
*+[r]{\fObj_\A}
&&
*+[l]{\fM(\{1,2\},\bu)_\A}
\ar@<.5ex>[ll]_(0.75){\bs\si(\{1,2\})}
}
\end{gathered}
\label{ai4eq17}
\e
the square of 1-morphisms `$\longra$' commutes, it is {\it not\/} a
Cartesian square, and so Theorem \ref{ai2thm1} does not apply.

To get round this, define $\fF$ to be the {\it fibre product
stack\/} of the bottom left corner of \eq{ai4eq17}. Since the outer
square of \eq{ai4eq17} commutes, there exist 1-morphisms
$\pi_1,\pi_2,\phi$ in \eq{ai4eq17} unique up to 2-isomorphism, such
that \eq{ai4eq17} commutes, and the bottom left quadrilateral is a
Cartesian square. We can then add maps `$\dashra$' in \eq{ai4eq16}.
The bottom left quadrilateral of \eq{ai4eq16} commutes by
\eq{ai2eq4}, and the central triangle by \eq{ai2eq3}. It remains
only to show the right hand triangle commutes.

We may justify this as follows. Points of $\fF(\K)$
may be naturally identified with isomorphism classes
of quadruples $(X,S,T,U)$, where $X\in\A$ and
$S,T,U\subset X$ are subobjects of $X$ with
$S,\ldots,X/U\in\A$, such that $X=S\op T$. An
{\it isomorphism} $\phi:(X,S,T,U)\ra(X',S',T',U')$
is an isomorphism $\phi:X\ra X'$ in $\A$ such that
$\phi(S)=S'$, $\phi(T)=T'$, $\phi(U)=U'$. Then
$\pi_1,\pi_2$ act on $\fF(\K)$ by
\begin{equation*}
(\pi_1)_*:[(X,S,T,U)]\mapsto[(X,S,T)]
\quad\text{and}\quad
(\pi_2)_*:[(X,S,T,U)]\mapsto[(X,U)],
\end{equation*}
where we identify $[(\si,\io,\pi)]\!\in\!\M(\{1,2\},\bu)_\A$
and $[(\si',\io',\pi')]\!\in\!\M(\{1,2\},\le)_\A$~with
\begin{small}
\begin{gather*}
\bigl[\bigl(\si(\{1,2\}),\io(\{1\},\{1,2\}):\si(\{1\})\!\ra\!
\si(\{1,2\}),\io(\{2\},\{1,2\}):\si(\{2\})\!\ra\!
\si(\{1,2\})\bigr)\bigr],
\\
\text{\begin{normalsize}and\end{normalsize}}\quad
\bigl[\bigl(\si'(\{1,2\}),\io'(\{1\},\{1,2\}):\si'(\{1\})\ra
\si'(\{1,2\})\bigr)\bigr]
\quad\text{\begin{normalsize}respectively.\end{normalsize}}
\end{gather*}
\end{small}

There is a closed substack $\fG$ of $\fF$ such that $\phi:
\fM(I,\pr)_\A\ra\fG$ is a 1-isomorphism. A point
$[(X,S,T,U)]\in\fF(\K)$ lies in $\fG(\K)$ if and only if
$U=(S\cap U)\op(T\cap U)$. Here, since $X=S\op T$ we have
$(S\cap U)\op(T\cap U)\subset U\subset X$, but it can
happen that $(S\cap U)\op(T\cap U)\ne U$. To understand
this, consider the case in which $S,T,U$ are vector
subspaces of a vector space $X$.

Write $\CF(\fG),\CF(\fF\sm\fG)$ for
the subspaces of $\CF(\fF)$ supported on
$\fG(\K),\ab
\fF(\K)\sm\fG(\K)$ respectively. Then
$\CF(\fF)=\CF(\fG)\op\CF(\fF\sm\fG)$, so
it suffices to show the two triangles
\e
{}\hskip -.5in
\text{
\begin{small}
$\displaystyle
\begin{gathered}
\xymatrix@C=6pt@R=4pt{
&\CF\bigl(\fM(I,\pr)_\A\bigr)
\ar@<-2ex>[dd]^{\substack{\CF^\stk(Q(I,\pr,\\
\{1,2\},\le,\mu))}}
&&\CF\bigl(\fM(I,\pr)_\A\bigr)
\ar@<-2ex>[dd]^{\substack{\CF^\stk(Q(I,\pr,\\
\{1,2\},\le,\mu))}}
\\
\CF(\fG)
\ar@{-->}[dr]_(0.3){\CF^\stk(\pi_1)\,\,{}}
\ar@{-->}[ur]^(0.3){\phi^*}
&&
\qquad \CF(\fF\sm\fG)
\ar@{-->}[dr]_(0.3){\CF^\stk(\pi_1)}
\ar@{-->}[ur]^(0.3){\phi^*}
\\
& \CF\bigl(\fM(\{1,2\},\bu)_\A\bigr)
&& \CF\bigl(\fM(\{1,2\},\bu)_\A\bigr)
}
\end{gathered}
$
\end{small}
}
\hskip -.2in{}
\label{ai4eq18}
\e
commute. As $\phi:\fM(I,\pr)_\A\ra\fG$ is a 1-isomorphism we have
$\phi^*=\CF^\stk(\phi^{-1})$ on $\CF(\fG)$, so the left triangle
commutes by~\eq{ai2eq2}.

In the right triangle, $\phi^*\!=\!0$ as each $f\in\CF(\fF\sm\fG)$
is zero on $\fG(\K)\!=\! \phi_*(\M(I,\pr)_\A)$. We shall show
$\CF^\stk(\pi_1)\!=\!0$ too. Let
$[(X,S,T,U)]\!\in\!\fF(\K)\sm\fG(\K)$, so that
$(\pi_1)_*:[(X,S,T,U)]\!\mapsto\![(X,S,T)]$. We have stabilizer
groups $\Iso_\K([(X,S,T,U)])\!=\!\Aut(X,S,T,U)$ and $\Iso_\K
([(X,S,T)])\!=\!\Aut(X,S,T)$ in $\Aut(X)$. Thus Definition
\ref{ai2def4} yields
\e
m_{\pi_1}\bigl([(X,S,T,U)]\bigr)=\chi\bigl(\Aut(X,S,T)/
\Aut(X,S,T,U)\bigr).
\label{ai4eq19}
\e

As $X=S\op T$ there is a subgroup $\bigl\{\id_S+\al\id_T:\al
\in\K\sm\{0\}\bigr\}\cong\K^\t$ in the centre of $\Aut(X,S,T)$.
Since $U\ne(S\cap U)\op(T\cap U)$, it is easy to see this group
intersects $\Aut(X,S,T,U)$ in the identity. Thus $\K^\t$ acts {\it
freely} on the left on $\Aut(X,S,T)/\Aut(X,S,T,U)$, fibring it by
$\K^\t$ orbits. But $\chi(\K^\t)=0$, so properties of $\chi$ show
that $m_{\smash{\pi_1}}\bigl([(X,S,T,U)]\bigr)=0$ in \eq{ai4eq19}.
Definition \ref{ai2def4} then shows that $\CF^\stk(\pi_1)f=0$ for
all $f\in\CF(\fF\sm\fG)$. Hence the right triangle in \eq{ai4eq18}
commutes, and thus \eq{ai4eq16} commutes, and $\De$ is {\it
multiplicative}.

By Proposition \ref{ai4prop1}, to show $\CF_\fin(\fObj_\A)$
is a {\it bialgebra} it remains only to verify some compatibilities
\cite[\S 1.1.3]{Jose} between the unit $\de_{[0]}$ and counit $\ep$,
which follow from the easy identities $\ep(\de_{[0]})=1$,
$\De\,\de_{[0]}=\de_{[0]}\ot\de_{[0]}$ and $(f*g)([0])=f([0])g([0])$
for all $f,g\in\CF_\fin(\fObj_\A)$. This completes the proof.
\end{proof}

We can determine $\De,\ep$ on the subspace~$\CFi_\fin(\fObj_\A)$.

\begin{lem} If\/ $f\!\in\!\CFi_\fin(\fObj_\A)$ then $\De f
\!=\!f\!\ot\de_{[0]}\!+\de_{[0]}\!\ot\!f$ and\/~$\ep(f)\!=\!0$.
\label{ai4lem1}
\end{lem}

\begin{proof} If $X\in\A$ is indecomposable then
$\bs\si(\{1,2\})^*\bigl([X]\bigr)\subseteq\M(\{1,2\},\bu)_\A$ is two
points $[(\si,\io,\pi)],[(\si',\io',\pi')]$, where
$\si(\{1\})\!=\!\si(\{1,2\})\!=\!X$, $\si(\{2\})\!=\!0$ and
$\si'(\{1\})\!=\!0$, $\si'(\{2\})\!=\!\si'(\{1,2\})\!=\!X$. Thus
$\De\,\de_{[X]}=\de_{[X]}\ot\de_{[0]}+\de_{[0]}\ot\de_{[X]}$ by
\eq{ai4eq15}, proving the first equation as such $\de_{[X]}$ form a
basis for $\CFi_\fin(\fObj_\A)$. Also $\ep(f)=f([0])=0$ as $f$ is
supported on $[X]$ for $X$ indecomposable.
\end{proof}

Let $\g$ be a Lie algebra. Then as in Joseph \cite[\S 1.2.6]{Jose},
the universal enveloping algebra $U(\g)$ has the structure of a
bialgebra with comultiplication $\De$ and counit $\ep$ satisfying
$\De x=1\ot x+x\ot 1$, $\ep(x)=0$ for all $x\in\g\subset U(\g)$.
Since $\g$ generates $U(\g)$ as an algebra, $\De,\ep$ are
determined on the whole of $U(\g)$ by their values on $\g$.
So we deduce:

\begin{cor} In Proposition \ref{ai4prop3}, $\Phi_\fin$
is an isomorphism of\/~$\Q$-bialgebras.
\label{ai4cor1}
\end{cor}

When we try to make $\CF(\fObj_\A)$ into a bialgebra
in the same way, without assuming {\it finite support},
we run into the following problem: if we use spaces
$\CF(\cdots)$ rather than $\CF_\fin(\cdots)$ in
\eq{ai4eq14} then $\Psi$ is injective, but generally
{\it not\/} surjective. Thus $\Psi^{-1}$ does not
exist, and $\De$ in \eq{ai4eq15} is not well-defined.

There are two natural solutions to this. The first is
to omit $\Psi^{-1}$ in \eq{ai4eq15}, so $\De$ maps
$\CF(\fObj_\A)\!\ra\!\CF(\fObj_\A\!\t\!\fObj_\A)$,
where we regard $\CF(\fObj_\A\!\t\!\fObj_\A)$ as a {\it
topological completion} of $\CF(\fObj_\A)\ot\CF
(\fObj_\A)$. Then the proof of Theorem \ref{ai4thm5}
works with few changes, but what we get is not strictly
a bialgebra.

The second is to restrict to a subalgebra
$\H$ of $\CF(\fObj_\A)$ such that \eq{ai4eq15}
yields a well-defined comultiplication
$\De:\H\ra\H\ot\H$. That is, as $\Psi$ is injective,
$\Psi^{-1}$ is well-defined on $\Im\Psi$, so
\eq{ai4eq15} makes sense if $\CF^\stk(\bs\si
(\{1\})\t\bs\si(\{2\}))\ci(\bs\si(\{1,2\}))^*$ maps
$\H\ra\Im\Psi$. We take this approach in our next theorem.

\begin{thm} Let Assumption \ref{ai3ass} hold, $\L\subseteq\CFi
(\fObj_\A)$ be a Lie subalgebra, and\/ $\H_\L$ the subalgebra of\/
$\CF(\fObj_\A)$ generated by $\L$. In particular, $\L$ can be the
Lie subalgebra generated by functions supported on points $[X]$ for
$X\in\A$ simple, and then $\H_\L=\cC$, the composition algebra
of\/~\S\ref{ai41}.

Then \eq{ai4eq15} yields a well-defined comultiplication
$\De:\H_\L\!\ra\!\H_\L\!\ot\!\H_\L$, where
\e
\Psi:\CF(\fObj_\A)\ot\CF(\fObj_\A)\ra
\CF(\fObj_\A\t\fObj_\A)
\label{ai4eq20}
\e
is injective and\/ $\Psi^{-1}$ defined on $\Im\Psi$, and\/
$\ep:f\!\mapsto\!f([0])$ defines a counit\/ $\ep:\H_\L\!\ra\!\Q$,
which make $\H_\L$ into a cocommutative bialgebra isomorphic
to~$U(\L)$.
\label{ai4thm6}
\end{thm}

\begin{proof} First we show $\De$ is well-defined. Changing
our point of view, omit $\Psi^{-1}$ from \eq{ai4eq15} so
that $\De$ maps $\CF(\fObj_\A)\!\ra\!\CF(\fObj_\A\!\t\!
\fObj_\A)$ and is well-defined, and regard $\Psi$ in
\eq{ai4eq20} as an {\it identification}, so that
$\CF(\fObj_\A)\ot\CF(\fObj_\A)$ becomes a {\it
vector subspace} of $\CF(\fObj_\A\t\fObj_\A)$. Then
we must prove that
\e
\De(\H_\L)\subseteq \H_\L\ot\H_\L
\subseteq\CF(\fObj_\A)\ot\CF(\fObj_\A)\subseteq
\CF(\fObj_\A\t\fObj_\A).
\label{ai4eq21}
\e

The proof of Theorem \ref{ai4thm5} still shows that $\De$
is {\it multiplicative} with respect to the natural
product on $\CF(\fObj_\A\!\t\!\fObj_\A)$, which
is essentially the right hand column of \eq{ai4eq16}.
Furthermore, the subspaces $\H_\L\ot\H_\L$ and
$\CF(\fObj_\A)\ot\CF(\fObj_\A)$ are closed under
this product, which equals $*\ot *$ upon them. Let $f_1,
\ldots,f_n\in\L$. The proof of Lemma \ref{ai4lem1} shows
that $\De f_i=f_i\!\ot\de_{[0]}\!+\de_{[0]}\!\ot\!f_i$
in $\H_\L\ot\H_\L$, so by multiplicativity of $\De$ we have
\begin{equation*}
\De(f_1*\cdots*f_n)=
(f_1\!\ot\de_{[0]}\!+\de_{[0]}\!\ot\!f_1)(*\ot *)\cdots
(*\ot *)(f_n\!\ot\de_{[0]}\!+\de_{[0]}\!\ot\!f_n).
\end{equation*}

Thus $\De(f_1*\cdots*f_n)$ lies in $\H_\L\ot\H_\L$. As
$\H_\L$ is spanned by such $f_1*\cdots*f_n$ and $\de_{[0]}$,
and $\De\,\de_{[0]}=\de_{[0]}\ot\de_{[0]}\in\H_\L\ot\H_\L$
as in Theorem \ref{ai4thm5}, we have proved \eq{ai4eq21},
and $\De:\H_\L\ra\H_\L\ot\H_\L$ is well-defined. The
proof of Theorem \ref{ai4thm5} now shows $\H_\L$ is a
cocommutative bialgebra. Finally, Proposition
\ref{ai4prop2} shows that $\Phi$ in \eq{ai4eq9} restricts
to an injective morphism $U(\L)\ra\CF(\fObj_\A)$,
which has image $\H_\L$. Hence $\H_\L\cong U(\L)$ as an
algebra, and the isomorphism of bialgebras follows as in
Corollary~\ref{ai4cor1}.
\end{proof}

As in \cite[\S 1.1.7]{Jose}, a {\it Hopf algebra}
is a bialgebra $A$ equipped with an {\it antipode}
$S:A\ra A$ satisfying certain conditions. If a
bialgebra $A$ admits an antipode $S$, then $S$ is
unique. Now for $\g$ a Lie algebra, $U(\g)$ is actually
a Hopf algebra \cite[\S 1.2.6]{Jose}. Therefore in the
situations of Proposition \ref{ai4prop3} and Theorem
\ref{ai4thm6}, there must exist unique antipodes
$S:\CF_\fin(\fObj_\A)\ra\CF_\fin(\fObj_\A)$ and
$S:\H_\L\ra\H_\L$ making the bialgebras into Hopf algebras.

However, there does not appear to be a simple formula
for $S$ in terms of constructible functions. (The most
obvious answer, that $(Sf)([X])=(-1)^kf([X])$ if $X\in\A$
has $k$ indecomposable factors, does not work.) So we
shall not try to determine the antipodes~$S$.

\subsection{Other algebraic operations from finite posets}
\label{ai48}

We define a family of {\it multilinear operations}
$P_\sIp$ on~$\CF(\fObj_\A)$.

\begin{dfn} Let Assumption \ref{ai3ass} hold and $(I,\pr)$
be a finite poset. Using
\begin{equation*}
\xymatrix@C=60pt{{\prod_{i\in I}\fObj_\A} &
\fM(I,\pr)_\A
\ar[l]_{\prod_{i\in I}\bs\si(\{i\})}
\ar[r]^{\bs\si(I)}
& \fObj_\A,
}
\end{equation*}
define a multilinear operation $P_\sIp:\prod_{i\in I}\CF(\fObj_\A)
\ra\CF(\fObj_\A)$ by
\begin{equation*}
P_\sIp(f_i:i\in I)=\ts\CF^\stk(\bs\si(I))\bigl[\prod_{i\in I}
\bs\si(\{i\})^*(f_i)\bigr],
\end{equation*}
which exists as $\bs\si(I)$ is representable, $\prod_{i\in I}
\bs\si(\{i\})$ finite type by Theorem~\ref{ai3thm}.
\label{ai4def8}
\end{dfn}

This generalizes Definition \ref{ai4def1}, as $P_{\sst(\{1,2\},\le)}
(f_1,f_2)=f_1*f_2$. In this notation, Theorem \ref{ai4thm1} shows
$*$ is associative by proving that
\e
\begin{split}
P_{\sst(\{1,2\},\le)}\bigl(P_{\sst(\{1,2\},\le)}(f_1,f_2),f_3\bigr)&=
P_{\sst(\{1,2,3\},\le)}(f_1,f_2,f_3)\\
&=P_{\sst(\{1,2\},\le)}\bigl(f_1,P_{\sst(\{1,2\},\le)}(f_2,f_3)\bigr).
\end{split}
\label{ai4eq22}
\e
Here is a generalization, which shows that if we substitute one
operation $P_\sJl$ into another $P_\sKt$, we get a third $P_\sIp$.
It is a constructible functions version of the notion
\cite[Def.~5.7]{Joyc3} of {\it substitution} of configurations.

\begin{thm} Let Assumption \ref{ai3ass} hold, $(J,\ls),(K,\tl)$
be nonempty finite posets with\/ $J\cap K=\emptyset$, and\/
$l\in K$. Set\/ $I=J\cup(K\sm\{l\})$, and define $\pr$ on\/ $I$ by
\begin{equation*}
i\pr j\quad\text{for $i,j\in I$ if}\quad
\begin{cases}
i\ls j, & i,j\in J, \\
i\tl j, & i,j\in K\sm\{l\}, \\
l\tl j, & i\in J,\quad j\in K\sm\{l\}, \\
i\tl l, & i\in K\sm\{l\},\quad j\in J.
\end{cases}
\end{equation*}
Let\/ $f_j:j\in J$ and\/ $g_k:k\in K\sm\{l\}$ lie in
$\CF(\fObj_\A)$. Then
\e
\begin{split}
&P_\sIp\bigl(f_j:j\in J,\; g_k:k\in K\sm\{l\}\bigr)=\\
&P_\sKt\bigl(P_\sJl(f_j:j\in J)_l,\; g_k:k\in K\sm\{l\}\bigr).
\end{split}
\label{ai4eq23}
\e
\label{ai4thm7}
\end{thm}

\begin{proof} Define $\phi:I\ra K$ by $\phi(i)=l$ if $i\in J$,
and $\phi(i)=i$ if $i\in K\sm\{l\}$. Consider the commutative
diagram of 1-morphisms, and the corresponding diagram of
pullbacks and pushforwards:
\begin{gather}
\begin{gathered}
\hskip -.2in
\xymatrix@C=60pt@R=20pt{
{\prod_{i\in I}\!\fObj_\A} \\
{\fM(J,\ls)_\A\!\t\!\!\prod\limits_{k\in I\sm J}\!\!\fObj_\A}
\ar[u]^{\substack{(\prod_{j\in J}\bs\si(\{j\}))\t \\
\prod_{k\in I\sm J}\id_{\fObj_\A}}}
\ar[d]^(0.4){\bs\si(J)\t\prod_{k\in I\sm J}\id_{\fObj_\A}}
& \fM(I,\pr)_\A
\ar@(u,r)[ul]_(0.3){\quad\prod_{i\in I}\bs\si(\{i\})}
\ar[l]_(0.4){\substack{S(I,\pr,J)\t \\ \prod_{k\in I\sm J}\bs\si(\{k\})}}
\ar[d]_(0.7){Q(I,\pr,K,\tl,\phi)}
\ar[dr]^{\quad\bs\si(I)}\\
{\fObj_\A\!\t\!\!\prod\limits_{k\in I\sm J}\!\!\fObj_\A}
&\fM(K,\tl)_\A
\ar[l]^(0.48){\bs\si(\{l\})\t
\prod_{k\in I\sm J}\bs\si(\{k\})}
\ar[r]_(0.6){\bs\si(K)}
& \fObj_\A,\qquad
}
\hskip -.5in
\end{gathered}
\label{ai4eq24}
\allowdisplaybreaks \\[5pt]
\text{
\begin{small}
$\displaystyle
\begin{gathered}
\hskip -.1in
\xymatrix@C=30pt@R=22pt{
\CF\bigl(\prod_{i\in I}\!\fObj_\A\bigr)
\ar[d]_{\substack{((\prod_{j\in J}\bs\si(\{j\}))\t \\
\prod_{k\in I\sm J}\id_{\fObj_\A})^*}}
\ar@(r,u)[dr]^(0.7){\qquad
(\prod_{i\in I}\bs\si(\{i\}))^*}
\\
\CF\bigl(\fM(J,\ls)_\A\!\t\!\prod_{k\in I\sm J}\!\fObj_\A\bigr)
\ar[d]^(0.35){\CF^\stk(\bs\si(J)\t\prod_{k\in I\sm J}\id_{\fObj_\A})}
\ar[r]
\ar@<1ex>@{}[r]^(0.57){\substack{(S(I,\pr,J)\t \\
\prod_{k\in I\sm J}\bs\si(\{k\}))^*}}
&
\CF\bigl(\fM(I,\pr)_\A\bigr)
\ar[d]_(0.65){\CF^\stk(Q(I,\pr,K,\tl,\phi))}
\ar[dr]^(0.5){\quad\CF^\stk(\bs\si(I))}
\\
\CF\bigl(\fObj_\A\!\t\!\prod_{k\in I\sm J}\!\fObj_\A\bigr)
\ar[r]
\ar@<-1ex>@{}[r]_(0.6){(\bs\si(\{l\})\t\prod_{k\in I\sm J}
\bs\si(\{k\}))^*}
& \CF\bigl(\fM(K,\tl)_\A\bigr)
\ar[r]_(0.48){\CF^\stk(\bs\si(K))}
&
\CF\bigl(\fObj_\A\bigr).\qquad\qquad
}
\hskip -.8in
\end{gathered}
$
\end{small}
}
\label{ai4eq25}
\end{gather}
Calculation shows that \eq{ai4eq23} holds provided \eq{ai4eq25}
commutes. The proof of this is similar to that for \eq{ai4eq3} in
Theorem \ref{ai4thm1}. By \cite[Th.~7.10]{Joyc3}
\begin{equation*}
\xymatrix@C=130pt@R=8pt{
*+[r]{\fM(J,\ls)_\A}
\ar[d]^{\,\bs\si(J)}
& *+[l]{\fM(I,\pr)_\A}
\ar[l]^{S(I,\pr,J)}
\ar[d]_{Q(I,\pr,K,\tl,\phi)\,}
\\
*+[r]{\fObj_\A}  & *+[l]{\fM(K,\tl)_\A}
\ar[l]_{\bs\si(\{l\})}
}
\end{equation*}
is a Cartesian square, so the square in \eq{ai4eq24} is Cartesian,
and the square in \eq{ai4eq25} commutes by \eq{ai2eq4}. We leave the
details to the reader.
\end{proof}

Applying the theorem and induction shows that {\it any}
multilinear operation on $\CF(\fObj_\A)$ obtained by
combining operations $P_\sJl$ is of the form $P_\sIp$
for some poset $(I,\pr)$. For instance, the posets
$(\{1,2\},\bu)$, $(\{1,2,3\},\bu)$ of Definition
\ref{ai4def7} and Theorem \ref{ai4thm7} give an
analogue of~\eq{ai4eq22}:
\begin{align*}
P_{\sst(\{1,2\},\bu)}\bigl(P_{\sst(\{1,2\},\bu)}(f_1,f_2),f_3\bigr)&=
P_{\sst(\{1,2,3\},\bu)}(f_1,f_2,f_3)\\
&=P_{\sst(\{1,2\},\bu)}\bigl(f_1,P_{\sst(\{1,2\},\bu)}(f_2,f_3)\bigr).
\end{align*}
This shows that $P_{\sst(\{1,2\},\bu)}$ gives an {\it associative,
commutative multiplication} on $\CF(\fObj_\A)$, different from~$*$.

Often the operations $P_\sIp$ map $\prod_{i\in I}\CF_\fin
(\fObj_\A)\ra\CF_\fin(\fObj_\A)$, which holds for $(I,\bu)$,
for example. Now Theorem \ref{ai4thm7} implies that Ringel--Hall
algebras $\H$ admit many {\it extra algebraic operations}
$P_\sIp$, which generalize multiplication $*$, and satisfy many
compatibilities. It is an interesting question whether these
operations may be useful tools in studying algebras which occur
as Ringel--Hall algebras, such as certain $U(\g)$. See Remark
\ref{ai4rem2} below on this.

Combining the ideas of Definitions \ref{ai4def2} and \ref{ai4def8}
we may also define $P_\sIp:\prod_{i\in I}\dLCF(\fObj_\A)\ra\ab
\dLCF(\fObj_\A)$, satisfying the analogue of Theorem~\ref{ai4thm7}.

\subsection{Examples from quivers}
\label{ai49}

Let $\Ga$ be a Dynkin diagram which is a disjoint union of
diagrams of type $A,D$ or $E$, and $\g$ be the corresponding
a finite-dimensional semisimple Lie algebra over $\Q$. Then
we have a decomposition $\g=\n_+\op\h\op\n_-$, where $\h$ is
a Cartan subalgebra, and the nilpotent Lie subalgebra $\n_+$
is a direct sum of one-dimensional subspaces indexed by the
set $\Phi_+$ of {\it positive roots} $\al$ of~$\g$.

Gabriel showed that if $Q=(Q_0,Q_1,b,e)$ is a quiver with
underlying graph $\Ga$, then isomorphism classes $[V]$ of
indecomposable representations $V$ of $Q$ are in 1-1
correspondence with $\al\in\Phi_+$. Later, Ringel
\cite{Ring2} used Ringel--Hall algebras over finite
fields to recover the Lie bracket on $\n_+$ on the vector
space spanned by such $[V]$. Here is a constructible functions
version of this, adapted from Riedtmann \cite{Ried} and
Frenkel et al.~\cite[\S 4]{FMV}.

\begin{ex} Let $\g=\n_+\op\h\op\n_-$, $\Ga$ and $Q$ be as above.
Set $\A=\modKQ$, the abelian category of representations of $Q$ over
$\K$, and define $K(\A)=\Z^{Q_0},\fF_\A$ as in
\cite[Ex.~10.5]{Joyc3}. Then Assumption \ref{ai3ass} holds by
\cite[\S 10]{Joyc3}. Also Gabriel's result implies {\it all
constructible sets in $\fObj_\A$ are finite}, so
$\CF(\fObj_\modKQ)=\CF_\fin(\fObj_\modKQ)$ and Assumption
\ref{ai4ass} holds automatically.

Now Riedtmann's Lie algebra $L(\C\vec{Q})\ot_\Z\Q$ of \cite[\S
2]{Ried} coincides exactly with the Lie algebra
$\CFi(\fObj_\modKQ)=\CF_\fin^\ind(\fObj_\modKQ)$ defined in
\S\ref{ai44} above when $\K=\C$. So \cite[Prop., p.~542]{Ried} gives
an isomorphism of Lie algebras
\begin{equation*}
\CFi(\fObj_\modKQ)=\CF_\fin^\ind(\fObj_\modKQ)\cong\n_+.
\end{equation*}
Proposition \ref{ai4prop3} and Corollary \ref{ai4cor1} thus
yield an isomorphism of bialgebras
\e
\CF(\fObj_\modKQ)=\CF_\fin(\fObj_\modKQ)\cong U(\n_+).
\label{ai4eq26}
\e
\label{ai4ex1}
\end{ex}

\begin{rem} We defined $P_\sIp:\prod_{i\in I}\CF(\fObj_\modKQ)
\ra\CF(\fObj_\modKQ)$ for each finite poset $(I,\pr)$ in
\S\ref{ai48}, which by \eq{ai4eq26} yield operations $P_\sIp$
on $U(\n_+)$ generalizing multiplication, and satisfying
various compatibilities. The author does not know if these
$P_\sIp$ on $U(\n_+)$ are good for anything.

However, we can say one thing: in general they depend on the choice
of {\it orientation} on $\Ga$ used to make $Q$. For example, if
$\Ga$ is $A_2$: $\smash{\mathop{\bu} \limits^{\sst
i}-\mathop{\bu}\limits^{\sst j}}$, then calculation shows that the
two possible choices $\smash{ \mathop{\bu}\limits^{\sst
i}\ra\mathop{\bu}\limits^{\sst j}}$ and
$\smash{\mathop{\bu}\limits^{\sst i}\leftarrow\mathop{\bu}
\limits^{\sst j}}$ for $Q$ give different answers for
$P_{\sst(\{1,2\},\bu)}$ on $U(\n_+)$. So the $P_\sIp$ on $U(\n_+)$
do not seem to be canonical.
\label{ai4rem2}
\end{rem}

We generalize Example \ref{ai4ex1} to {\it affine Lie algebras} and
{\it Kac--Moody Lie algebras}, based on Lusztig \cite[\S 12]{Lusz}
and Frenkel et al.~\cite[\S 5.6]{FMV}.

\begin{ex} Let $\Ga$ be a finite undirected graph with vertex set
$Q_0$ and no edge joining a vertex with itself, and let $\g$ be
the corresponding {\it Kac--Moody algebra} over $\Q$, as in Kac
\cite{Kac}. Then $\g$ is a $\Q$-Lie algebra with generators
$e_i,f_i,h_i$ for $i\in Q_0$, which satisfy certain relations.
It has a decomposition $\g=\n_+\op\h\op\n_-$, where $\n_+,\n_-,\h$
are the Lie subalgebras of $\g$ generated by the $e_i$, the $f_i$ and
the $h_i$ respectively, and $\h$ is abelian, the Cartan subalgebra.

Let $Q$ be a quiver with underlying graph $\Ga$, and {\it without
oriented cycles}. Take $\A=\modKQ$, and define $K(\A)=\Z^{Q_0},
\fF_\A$ as in \cite[Ex.~10.5]{Joyc3}. Then Assumption \ref{ai3ass}
holds by \cite[\S 10]{Joyc3}, but in general Assumption \ref{ai4ass}
does not. There is up to isomorphism one simple object
$V_i\in\modKQ$ for each $i\in Q_0$. Write $\L$ for the Lie
subalgebra of $\CF(\fObj_\modKQ)$ generated by the isomorphism
classes of simples $\de_{[V_i]}$ for $i\in Q_0$. Then the subalgebra
of $\CF(\fObj_\modKQ)$ generated by $\L$ is the composition algebra
$\cC$, and there is an isomorphism of bialgebras $U(\L)\cong\cC$ by
Theorem~\ref{ai4thm6}.

There are now natural identifications between the algebras
$\CF(\fObj_\modKQ)$ above and $\M(\Om)$ in Lusztig
\cite[\S 10.19]{Lusz}, and between their complexification
and $L(Q)$ in Frenkel et al.\ \cite[\S 2.2]{FMV}; and
between the algebras $\cC$ above and $\M_0(\Om)$ in Lusztig
\cite[\S 10.19]{Lusz}; and between the Lie algebras $\L\ot_\Q\C$
above and $\n^*(Q)$ in Frenkel et al.\ \cite[\S 2.2.5]{FMV};
{\it except that\/} Lusztig uses the opposite order of
multiplication to us, as in Remark~\ref{ai4rem1}.

From Lusztig \cite[Prop.~10.20]{Lusz} (which he attributes to
Schofield), there is a unique algebra isomorphism $\cC\cong
U(\n_+)$ identifying the generators $\de_{[V_i]}$ of $\cC$
with the generators $e_i$ of $\n_+$ for $i\in Q_0$. This implies
an isomorphism of the Lie subalgebras $\L$ and $\n_+$
of $\cC$ and~$U(\n_+)$.

When $\Ga$ is an {\it affine Dynkin diagram}, this isomorphism
$\L\cong\n_+$ also follows from Frenkel et al.\
\cite[Cor.~5.6.27]{FMV}. They also describe \cite[\S 5.6]{FMV} the
{\it isomorphism classes of indecomposables} in $\modKQ$, and
\cite[Cor.~5.6.30]{FMV} they define $\L$ explicitly as a subspace
of~$\CFi_\fin (\fObj_\modKQ)$.
\label{ai4ex2}
\end{ex}

We return to these examples in~\S\ref{ai54}.

\section{Stack algebras}
\label{ai5}

We now develop analogues of the Ringel--Hall algebras
$\CF(\fObj_\A)$ of \S\ref{ai4} in which we replace constructible
functions by {\it stack functions}, as in \S\ref{ai23} and
\S\ref{ai25}. We call the corresponding algebras {\it stack
algebras}. Throughout this section $\K$ is an algebraically closed
field of {\it arbitrary characteristic}, except when we specify
characteristic zero for results comparing stack and constructible
functions.

\subsection{Different kinds of stack algebras}
\label{ai51}

Following Definition \ref{ai4def1}, we define a {\it multiplication}
$*$ on~$\uSF(\fObj_\A)$.

\begin{dfn} Suppose Assumption \ref{ai3ass} holds. Using the
1-morphism diagram
\begin{equation*}
\xymatrix@C=70pt{ \fObj_\A\t\fObj_\A & \fM(\{1,2\},\le)_\A
\ar[l]_{\bs\si(\{1\})\t\bs\si(\{2\})} \ar[r]^(0.55){\bs\si(\{1,2\})}
& \fObj_\A, }
\end{equation*}
define a bilinear operation $*:\uSF(\fObj_\A)\t\uSF(\fObj_\A)\ra
\uSF(\fObj_\A)$~by
\e
f*g=\bs\si(\{1,2\})_*\bigl[(\bs\si(\{1\})\t\bs\si(\{2\}))^*(f\ot
g)\bigr],
\label{ai5eq1}
\e
using operations $\phi_*,\phi^*,\ot$ from \S\ref{ai23}. Here
$(\bs\si(\{1\})\t \bs\si(\{2\}))^*$ is well-defined as
$\bs\si(\{1\})\t\bs\si(\{2\})$ is of finite type by Theorem
\ref{ai3thm}(c). As $\bs\si(\{1,2\})$ is representable by Theorem
\ref{ai3thm}(b), the restriction maps $*:\SF(\fObj_\A)\!\t\!
\SF(\fObj_\A)\!\ra\!\SF(\fObj_\A)$, that is, $\SF(\fObj_\A)$ is
closed under~$*$.

If Assumption \ref{ai2ass} holds, \eq{ai5eq1} defines
multiplications $*$ on $\uSF(\fObj_\A,\Up,\La)$, $\uoSF,\oSF
(\fObj_\A,\Up,\La)$, $\uoSF,\oSF(\fObj_\A,\Up,\La^\ci)$ and
$\uoSF,\oSF(\fObj_\A,\Th,\Om)$. Write
$\bde_{[0]}\!\in\!\SF(\fObj_\A),\ldots,\oSF(\fObj_\A,\Th,\Om)$ for
$\bde_C$ in Definition \ref{ai2def6} with $C\!=\!\{[0]\}$. Then
$\bde_{[0]}=[(\Spec\K,0)]$, where $0:\Spec\K\!\ra\!\fObj_\A$
corresponds to $0\!\in\!\A$.
\label{ai5def1}
\end{dfn}

Here is the analogue of Theorem~\ref{ai4thm1}.

\begin{thm} If Assumptions \ref{ai2ass}, \ref{ai3ass} hold\/
$\uSF,\SF(\fObj_\A),\uSF(\fObj_\A,\Up,\La),\ab\ldots,\ab
\oSF(\fObj_\A,\Th,\Om)$ are algebras with identity $\bde_{[0]}$ and
multiplication $*$. Also
\e
\begin{gathered}
\pi_{\fObj_\A}^\stk:\SF(\fObj_\A)\!\ra\!\CF(\fObj_\A),\;\>
\bar\pi_{\fObj_\A}^\stk:\oSF(\fObj_\A,\Up,\La^\ci)\!\ra\!
\CF(\fObj_\A),\\
\text{and}\qquad
\bar\pi_{\fObj_\A}^\stk:\oSF(\fObj_\A,\Th,\Om)\ra\CF(\fObj_\A)
\end{gathered}
\label{ai5eq2}
\e
are morphisms of\/ $\K$-algebras when $\cha\K=0$ and\/
$\bar\pi_{\fObj_\A}^\stk$ is defined.
\label{ai5thm1}
\end{thm}

\begin{proof} To show $*$ is {\it associative} on $\uSF(\fObj_\A)$
we follow the proof of Theorem \ref{ai4thm1} replacing $\CF(\cdots)$
by $\uSF(\cdots)$, using Theorem \ref{ai2thm2} to show the analogue
of \eq{ai4eq3} commutes. To show $\bde_{[0]}$ is an identity on
$\uSF(\fObj_\A)$ we can copy Theorem \ref{ai4thm1}, or note that
$\bde_{[0]}=[(\Spec\K,0)]$ and show directly that $[(\Spec\K,0)]*
[(\fR,\rho)]=[(\fR,\rho)]$ by constructing an explicit 1-isomorphism
\begin{equation*}
i:\fR\longra(\Spec\K\t\fR)_{0\t\rho,\fObj_\A\t\fObj_\A,\bs\si(\{1\})\t
\bs\si(\{2\})}\fM(\{1,2\},\le)_\A
\end{equation*}
with $\bs\si(\{1,2\})\ci\pi_{\fM(\{1,2\},\le)_\A}\ci i$ 2-isomorphic
to $\rho$. Thus $\uSF(\fObj_\A)$ is a $\Q$-algebra. Also
$\SF(\fObj_\A)$ is a $\Q$-subalgebra as it is closed under $*$. The
arguments for $\uSF(\fObj_\A,\Up,\La),\ldots,\oSF(\fObj_\A,\Th,\Om)$
are the same.

To prove \eq{ai5eq2} are $\K$-algebra morphisms, note that
$\bde_{[0]}=\io_{\fObj_\A}(\de_{[0]})$ and $\pi_{\fObj_\A}^\stk
\ci\io_{\fObj_\A}$ is the identity by Definition \ref{ai2def6}, so
$\pi_{\fObj_\A}^\stk(\bde_{[0]})=\de_{[0]}$. For the
$\pi_{\fObj_\A}^\stk$ case, $\pi_{\fObj_\A}^\stk(f*g)=
\pi_{\fObj_\A}^\stk(f)*\pi_{\fObj_\A}^\stk(g)$ follows from
equations \eq{ai4eq1}, \eq{ai5eq1} and Theorem \ref{ai2thm3}(b),(c).
The $\bar\pi_{\fObj_\A}^\stk$ cases are the same.
\end{proof}

It is easy to show by example that $\io_{\fObj_\A}:\CF(\fObj_\A)\ra
\SF(\fObj_\A)$ is {\it not\/} in general an algebra morphism, as the
$\io_\fF$ do not commute with pushforwards. We will need the
analogues of the $P_\sIp$ of \S\ref{ai48} in \S\ref{ai52}, so we
define them now.

\begin{dfn} Let Assumptions \ref{ai2ass} and \ref{ai3ass} hold
and $(I,\pr)$ be a finite poset. Define a multilinear operation
$P_\sIp\!:\!\prod_{i\in I}\uSF(\fObj_\A)\!\ra\!\uSF(\fObj_\A)$ on
$\uSF(\fObj_\A)$ (or on $\uSF(\fObj_\A,\Up,\La),\uoSF(\fObj_\A,\Up,
\La),\uoSF(\fObj_\A,\Up,\La^\ci)$, or $\uoSF(\fObj_\A,\Th,\Om)$) by
\e
P_\sIp(f_i:i\in I)=\ts\bs\si(I)_*\bigl[\bigl(\prod_{i\in I}
\bs\si(\{i\})\bigr)^*\bigl(\bigot_{i\in I}f_i\bigr)\bigr],
\label{ai5eq3}
\e
using 1-morphisms from $\fM(I,\pr)_\A$. Since $\bs\si(I)$ is
representable, $\SF(\fObj_\A)$ and $\oSF(\fObj_\A,*,*)$ are closed
under $P_\sIp$. Also $P_{\sst(\{1,2\},\le)}$ is $*$ in
Definition~\ref{ai5def1}.
\label{ai5def2}
\end{dfn}

Here is the analogue of Theorem \ref{ai4thm7}, proved as for
Theorem~\ref{ai5thm1}.

\begin{thm} Let Assumptions \ref{ai2ass} and \ref{ai3ass} hold and\/
$(J,\ls),(K,\tl),l$ and\/ $(I,\pr)$ be as in Theorem \ref{ai4thm7}.
Let\/ $f_j:j\in J$ and\/ $g_k:k\in K\sm\{l\}$ lie in
$\uSF(\fObj_\A)$ or $\uSF,\uoSF(\fObj_\A,*,*)$. Then \eq{ai4eq23}
holds in the same space.
\label{ai5thm2}
\end{thm}

Now apart from $\uSF(\fObj_\A,\Up,\La)$ these algebras are all
inconveniently large for our later work, and we will find it useful
to define subalgebras $\SFa(\fObj_\A)$, $\oSFa(\fObj_\A,*,*)$ using
the algebra structure on stabilizer groups in~$\fObj_\A$.

\begin{dfn} Suppose Assumptions \ref{ai2ass} and \ref{ai3ass} hold,
and $[(\fR,\rho)]$ be a generator of $\SF(\fObj_\A)$. Let
$r\in\fR(\K)$ with $\rho_*(r)=[X]\in\fObj_\A(\K)$, for $X\in\A$.
Then $\rho$ induces a morphism of stabilizer $\K$-groups
$\rho_*:\Iso_\K(r)\ra\Iso_\K([X]) \cong\Aut(X)$. As $\rho$ is
representable this is {\it injective}, and induces an isomorphism of
$\Iso_\K(r)$ with a $\K$-subgroup of~$\Aut(X)$.

Now $\Aut(X)=\End(X)^\t$ is the $\K$-group of invertible elements in
a {\it finite-dimensional\/ $\K$-algebra} $\End(X)=\Hom(X,X)$. We
say the $[(\fR,\rho)]$ {\it has algebra stabilizers} if whenever
$r\in\fR(\K)$ with $\rho_*(r)=[X]$, the $\K$-subgroup
$\rho_*\bigl(\Iso_\K(r)\bigr)$ in $\Aut(X)$ is the $\K$-group $A^\t$
of invertible elements in a $\K$-subalgebra $A$ in $\End(X)$. (Here
when $X\cong 0$ we allow $\End(X)=\{0\}$ as a $\K$-algebra.)

Write $\SFa(\fObj_\A),\oSFa(\fObj_\A,\!\Up,\La),
\oSFa(\fObj_\A,\!\Up,\La^\ci),\oSFa(\fObj_\A,\Th,\Om)$ for the
subspaces of $\SF(\fObj_\A)$, $\oSF(\fObj_\A,\Up,\La)$,
$\oSF(\fObj_\A,\Up,\La^\ci)$, $\oSF(\fObj_\A,\Th,\Om)$ respectively
spanned over $\Q,\La,\La^\ci$ or $\Om$ by $[(\fR,\rho)]$ with
algebra stabilizers.

One of the relations \cite[Def.~5.17(iii)]{Joyc2} defining
$\oSF(\fObj_\A,*,*)$ mixes $[(\fR,\rho)]$ with and without algebra
stabilizers, so $\oSF(\fObj_\A,*,*)$ can contain expressions
$\sum_ic_i[(\fR_i,\rho_i)]$ in which $[(\fR_i,\rho_i)]$ does not
have algebra stabilizers. Propositions \ref{ai5prop2} and
\ref{ai5prop3} are tools for getting round the problems this causes.
Also, $\io_{\fObj_\A}$ in Definition \ref{ai2def6} maps
$\CF(\fObj_\A)\ra\SFa(\fObj_\A)$, as it involves $[(\fR,\rho)]$ with
$\rho$ an inclusion of substacks, so the condition above holds
with~$A=\End(X)$.
\label{ai5def3}
\end{dfn}

\begin{prop} $\SFa(\fObj_\A)$ and\/ $\oSFa(\fObj_\A,*,*)$ are closed
under $*$ and contain $\bde_{[0]}$, and so are subalgebras. They are
also closed under the~$P_\sIp$.
\label{ai5prop1}
\end{prop}

\begin{proof} When $[(\si,\io,\pi)]\in\M(\{1,2\},\le)_\A$ the
stabilizer group maps
\begin{align*}
&\bs\si(\{1,2\})_*:\Aut\bigl((\si,\io,\pi)\bigr)\ra
\Aut(\si(\{1,2\})),\\
&(\bs\si(\{1\}\t\bs\si(\{2\}))_*:\Aut\bigl((\si,\io,\pi)\bigr)
\ra\Aut(\si(\{1\}))\t\Aut(\si(\{2\}))
\end{align*}
are the restrictions to $\K$-groups of invertible elements of
$\K$-algebra morphisms $\End\bigl((\si,\io,\pi)\bigr)\!\ra\!
\End(\si(\{1,2\}))$, $\End\bigl((\si,\io,\pi)\bigr)\!\ra\!
\End(\si(\{1\}))\t\End(\si(\{2\}))$. Using this and fibre products
of $\K$-algebras, we find $\SFa(\fObj_\A)$ is closed under $*$.
Replacing $(\{1,2\},\le)$ by $(I,\pr)$ we also see $\SFa(\fObj_\A)$
is closed under $P_\sIp$. It contains $\bde_{[0]}$ as
$\bde_{[0]}=\io_{\fObj_\A}(\de_{[0]})$ and $\io_{\fObj_\A}$ maps to
$\SFa(\fObj_\A)$. The arguments for $\oSFa(\fObj_\A,*,*)$ are the
same.
\end{proof}

Now much of \cite{Joyc2} involves taking $[(\fR,\rho)]$ with
$\fR\cong[V/G]$ and then making linear combinations of
$[([W/H],\rho\ci\io^{W,H})]$ for certain $\K$-subgroups $H\subseteq
G$ and $H$-invariant $\K$-subvarieties $W\subseteq V$. If
$[(\fR,\rho)]$ has algebra stabilizers, for a general 1-isomorphism
$\fR\cong[V/G]$ these $[([W/H],\rho\ci\io^{W,H})]$ may {\it not\/}
have algebra stabilizers, which is why \cite[Def.~5.17(iii)]{Joyc2}
mixes $[(\fR,\rho)]$ with and without algebra stabilizers, as we
said above.

The next two propositions construct special 1-isomorphisms
$\fR\cong[V/A^\t]$ such that the $[([W/H],\rho\ci\io^{W,H})]$
automatically have algebra stabilizers. Using these we can do the
operations of \cite{Joyc2} within $\SFa(\fObj_\A)$
and~$\oSFa(\fObj_\A,*,*)$.

\begin{prop} Let\/ $S\subset\fObj_\A(\K)$ be constructible. Then
there exists a finite decomposition $S=\coprod_{l\in L}\fF_l(\K)$,
where $\fF_l$ is a finite type $\K$-substack of\/ $\fObj_\A$, and\/
$1$-isomorphisms $\fF_l\cong[U_l/A_l^\t]$, for $A_l$ a
finite-dimensional\/ $\K$-algebra and\/ $U_l$ a quasiprojective
$\K$-variety acted on by $A_l^\t$, such that if\/ $u\in U_l(\K)$
projects to $[X]\in\fF_l(\K)\subset\fObj_\A(\K)$ then there exists a
subalgebra $B_u$ of\/ $A_l$ with\/ $\Stab_{A_l^\t}(u)=B_u^\t$ and an
isomorphism $B_u\cong\End(X)$ compatible with the
isomorphisms~$\Stab_{A_l^\t}(u)\cong\Iso_\K([X])\cong\Aut(X)$.
\label{ai5prop2}
\end{prop}

\begin{proof} A result of Kresch \cite[Prop.~3.5.9]{Kres} used
in \cite{Joyc1,Joyc2} shows an algebraic $\K$-stack can be {\it
stratified by global quotient stacks} if and only if it has {\it
affine geometric stabilizers}. Applied to $S\subset\fObj_\A(\K)$,
this yields a finite decomposition $S=\coprod_{l\in L}\fF_l(\K)$ for
$\K$-substacks $\fF_l$ of $\fObj_\A$ with $\fF_l\cong[U_l/A_l^\t]$
for $A_l$ a finite-dimensional\/ $\K$-algebra which we may take to
be $\GL(m_l,\K)$, and $U_l$ a quasiprojective $\K$-variety acted on
by~$A_l^\t$.

Kresch \cite[p.~520]{Kres} uses affine geometric stabilizers to find
a vector bundle acted on by a group scheme. In our case we know
more: the stabilizers are of the form $\Aut(X)=\End(X)^\t$ for
finite-dimensional algebras $\End(X)$. Putting this into the
argument we can take the group scheme action to be the restriction
of an algebra representation of an algebra scheme. The proposition
follows.
\end{proof}

\begin{prop} Let\/ $[(\fR,\rho)]\in\SF(\fObj_\A)$ and use the
notation of Proposition \ref{ai5prop2} with\/ $S=\rho_*(\fR(\K))$.
Then $[(\fR,\rho)]=\sum_{l\in L}[([Y_l/A_l^\t],\rho_l)]$ where $Y_l$
is a quasiprojective $\K$-variety and\/ $\rho_l:[Y_l/A_l^\t]\ra\fF_l
\subseteq\fObj_\A$ is induced by an $A_l^\t$-equivariant morphism
$\phi_l:Y_l\ra U_l$, under the $1$-isomorphism
$\fF_l\cong[U_l/A_l^\t]$. Moreover $[(\fR,\rho)]$ has algebra
stabilizers if and only if\/ $\Stab_{A_l^\t}(y)=C_y^\t$ for some
subalgebra $C_y\subseteq A_l$, for all\/ $l\in L$ and\/~$y\in
Y_l(\K)$.
\label{ai5prop3}
\end{prop}

\begin{proof} Let $\pi_l:U_l\ra\fF_l\subset\fObj_\A$ be the natural
projection from $\fF_l\cong[U_l/A_l^\t]$. As $\fR,U_l$ are finite
type and $\rho$ is representable,
$Z_l=\fR\t_{\rho,\fObj_\A,\pi_l}U_l$ is a {\it finite type algebraic
$\K$-space}, with an action of $A_l^\t$. Thus there exists a finite
decomposition $Z_l=\coprod_{i\in I_l}Z_l^i$ for $A_l^\t$-invariant
quasiprojective $\K$-subvarieties $Z_l^i$. Define $Y_l$ to be the
scheme-theoretic disjoint union of the $Z_l^i$ for~$i\in I_l$.

By our convention in \S\ref{ai21} that $\K$-varieties need not be
irreducible, nor connected, $Y_l$ is a quasiprojective $\K$-variety,
acted on by $A_l^\t$. Let $\pi_{Z_l}:Y_l\ra Z_l$ be the obvious
morphism and $\pi_\fR:Z_l\ra\fR$, $\pi_{U_l}:Z_l\ra U_l$ the
projections from the fibre product. Then $\pi_\fR$ and
$\pi_\fR\ci\pi_{Z_l}$ are $A_l^\t$-invariant and so push down to
1-morphisms $\pi'_\fR:[Z_l/A_l^\t]\ra\fR$ and
$\pi''_\fR:[Y_l/A_l^\t] \ra\fR$. Define $\rho_l=\rho\ci\pi''_\fR$
and $\phi_l=\pi_{U_l}\ci\pi_{Z_l}$. It is now easy to see that
\begin{equation*}
[(\fR,\rho)]=\ts\sum_{l\in L}[([Z_l/A_l^\t],
\rho\ci\pi'_\fR)]=\sum_{l\in L}[([Y_l/A_l^\t],\rho_l)],
\end{equation*}
since $\pi'_\fR$ embeds $[Z_l/A_l^\t]$ as the $\K$-substack of $\fR$
over $\fF_l$, and $[Z_l/A_l^\t]$ and $[Y_l/A_l^\t]$ split into the
same pieces $[Z_l^i/A_l^\t]$ for $i\in I_l$. The first part follows.

For the second part, if $r\in\fR(\K)$ with $\rho_*(r)=[X]$ then
$[X]\in\fF_l(\K)$ for unique $l\in L$, and $r=(\pi_\fR)_*(z)$ for
$z\in Z_l(\K)$. As $(\pi_{Z_l})_*:Y_l(\K)\ra Z_l(\K)$ is a bijection
$z=(\pi_{Z_l})_*(y)$ for unique $y\in Y_l(\K)$. Let
$u=(\phi_l)_*(y)\in U_l(\K)$. Then $B_u^\t=\Stab_{\smash{A_l^\t}}(u)
\cong\Aut(X)$ as $u$ projects to $[X]$, for some subalgebra
$B_u\subseteq A_l$ with compatible isomorphism $B_u\cong\End(X)$.
Now $\Iso_\K(r)\cong\Stab_{\smash{A_l^\t}}(y)$, and
$\rho_*:\Iso_\K(r)\ra \Iso_\K([X])$ is just inclusion
$\Stab_{\smash{A_l^\t}}(y)\subseteq B_u^\t$ as $\K$-subgroups of
$A_l^\t$. Thus, $[(\fR,\rho)]$ has algebra stabilizers if and only
if $\Stab_{\smash{A_l^\t}}(y)=C_y^\t$ for some subalgebra
$C_y\subseteq B_u$, for all $l,u,y$. But a subalgebra of $B_u$ is
the same as a subalgebra of $A_l$ lying in $B_u$, and the proof is
complete.
\end{proof}

\begin{cor} $\SFa(\fObj_\A)$ and\/ $\oSFa(\fObj_\A,*,*)$ are closed
under the operators $\Pi^\mu,\Pi^\vi_n,\hat\Pi^\nu_\fF$
of\/~{\rm\cite[\S 5.2]{Joyc2}}.
\label{ai5cor1}
\end{cor}

\begin{proof} Let $[(\fR,\rho)]$ have algebra stabilizers, and
use the notation of Proposition \ref{ai5prop3}. Then
$\Pi^\mu([(\fR,\rho)])=\sum_{l\in L}\Pi^\mu([(
[Y_l/A_l^\t],\rho_l)])$. The definition \cite[Def.~5.10]{Joyc2} of
$\Pi^\mu([([Y_l/A_l^\t],\rho_l)])$ gives a linear combination of
$[([W_l/H_l],\rho_l\ci\io^{W_l,H_l})]$ for certain $\K$-subgroups
$H_l\subseteq A_l^\t$ and $H_l$-invariant $\K$-subvarieties
$W_l\subseteq Y_l$. We may take $A_l=\End(\K^{m_l})$, and then
\cite[Ex.~5.7]{Joyc2} implies the $H_l$ appearing in the sum are of
the form $B_l^\t$ for subalgebras $B_l\subseteq A_l$. It easily
follows that $[([W_l/H_l],\rho_l\ci\io^{W_l,H_l})]$ has algebra
stabilizers, which proves what we want. The
$\Pi^\vi_n,\hat\Pi^\nu_\fF$ cases are the same.
\end{proof}

Combining Proposition \ref{ai5prop3} with the proof of the first
part of Proposition \ref{ai2prop2} we find $\oSFa(\fObj_\A,*,*)$ are
generated by $[(U\t[\Spec\K/T],\rho)]$ as in Proposition
\ref{ai2prop2} with algebra stabilizers. But $T\cong(\K^\t)^k\t K$
for $K$ finite abelian can only be of the form $B^\t$ for a
$\K$-algebra $B$ if $K$ is trivial, giving:

\begin{cor} $\oSFa(\fObj_\A,\Up,\La)$, $\oSFa(\fObj_\A,\Up,\La^\ci)$
and\/ $\oSFa(\fObj_\A,\Th,\Om)$ are generated over $\La,\La^\ci$
and\/ $\Om$ respectively by elements $[(U\t[\Spec\K/T],\rho)]$ with
algebra stabilizers, for\/ $U$ a quasiprojective $\K$-variety
and\/~$T\cong(\K^\t)^k$.
\label{ai5cor2}
\end{cor}

We define projections $\Pi_{[I,\ka]}$
on~$\oSFa(\fObj_\A),\oSFa(\fObj_\A,*,*)$.

\begin{dfn} Let Assumptions \ref{ai2ass} and \ref{ai3ass} hold. Write
$C(\A)=\bar C(\A)\sm\{0\}$, as in Definition \ref{ai3def4}. Consider
pairs $(I,\ka)$ with $I$ a finite set and $\ka:I\ra C(\A)$ a map.
Define an equivalence relation `$\approx$' on such $(I,\ka)$ by
$(I,\ka)\approx(I',\ka')$ if there exists a bijection $i:I\ra I'$
with $\ka'\ci i=\ka$. Write $[I,\ka]$ for the $\approx$-equivalence
class of $(I,\ka)$. For such an $[I,\ka]$ we will define projections
\e
\begin{split}
&\Pi_{[I,\ka]}:\SFa(\fObj_\A)\ra\SFa(\fObj_\A),\\
&\Pi_{[I,\ka]}:\oSFa(\fObj_\A,*,*)\ra\oSFa(\fObj_\A,*,*),
\end{split}
\label{ai5eq4}
\e
for $*,*=\Up,\La$ or $\Up,\La^\ci$ or $\Th,\Om$, using the operators
$\hat\Pi^\nu_\fF$ of \cite[Def.~5.15]{Joyc2}. Define
\begin{align*}
\nu:\bigl\{&(T,[X],\phi):\text{$T$ a $\K$-group isomorphic to
$(\K^\t)^k\!\t\!K$, $K$ finite abelian,}\\
&\text{$[X]\!\in\!\fObj_\A(\K)$,
$\phi:T\!\ra\!\Iso_\K([X])\!=\!\Aut(X)$ a $\K$-group
morphism}\bigr\}\!\ra\!\Q
\end{align*}
by $\nu(T,[X],\phi)=1$ if $T\cong(\K^\t)^{\md{I}}$, $\phi$ is
injective, and there exists a splitting $X\cong\bigop_{i\in I}X_i$
in $\A$ with $[X_i]=\ka(i)$ for all $i\in I$, such that $\phi(T)$ is
the $\K$-subgroup $\{\sum_{i\in I}\la_i\id_{X_i}:\la_i\in\K^\t\}$ in
$\Aut(X)$, and $\nu(T,[X],\phi)=0$ otherwise. This depends only on
the equivalence class $[I,\ka]$ of $(I,\ka)$. Note too that the
$\K$-subgroup of $\Aut(X)$ above is $A^\t$, where $A$ is the
subalgebra $\{\sum_{i\in I}\la_i\id_{X_i}:\la_i\in\K\}$ in
$\End(X)$. Then $\nu$ is an $\fObj_\A$-{\it weight function} in the
sense of \cite[Def.~5.15]{Joyc2}. Define projections $\Pi_{[I,\ka]}$
on $\SFa(\fObj_\A)$ and $\oSFa(\fObj_\A,*,*)$ to be the operators
$\hat\Pi^\nu_{\fObj_\A}$ of \cite[Def.~5.15]{Joyc2}. Corollary
\ref{ai5cor1} implies these map as in~\eq{ai5eq4}.
\label{ai5def4}
\end{dfn}

Here is what this means. The $\Pi_{[I,\ka]}$ are refinements of the
$\Pi^\vi_n$ of \cite[\S 5.2]{Joyc2}, which project to components
with {\it virtual rank\/} $n$. Now if $[X]\in\fObj_\A(\K)$ then
$\Iso_\K([X])\cong\Aut(X)$. Maximal tori in $\Aut(X)$ are of the
form $\{\sum_{i\in I}\la_i\id_{X_i}:\la_i\in\K^\t\}$, where
$X=\bigop_{i\in I}X_i$ with $0\not\cong X_i$ {\it indecomposable}.
Thus the rank of $\Iso_\K([X])$ is the {\it number of indecomposable
factors of\/}~$X$.

Now $[(\fR,\rho)]\in\SFa(\fObj_\A)$ also has stabilizers of the form
$A^\t$ for subalgebras $A\subseteq\End(X)$, so we can treat $A$ as
like $\End(Y)$ for some `object' $Y$, and $\rk\,A^\t$ as the `number
of indecomposables' in $Y$. As the $\Pi^\vi_n$ project to components
with `virtual rank' $n$, and in $\SFa(\fObj_\A)$ we equate rank with
number of indecomposables, so on $\SFa(\fObj_\A)$ we should think of
$\Pi^\vi_n$ as projecting to stack functions {\it supported on
objects with\/ $n$ virtual indecomposable factors}. That is, the
idea of `virtual rank' in $\SF(\fF)$ translates to `number of
virtual indecomposables' in~$\SFa(\fObj_\A)$.

Each $X\in\A$ splits as $X=\bigop_{i\in I}X_i$ with $0\not\cong X_i$
indecomposable, uniquely up to bijective changes of indexing set $I$
and isomorphisms of $X_i$. Defining $\ka:I\ra C(\A)$ by
$\ka(i)=[X_i]$, we see that the equivalence class $[I,\ka]$ depends
only on $[X]$, and $X$ has $\md{I}$ indecomposables. In the same
way, the $\Pi_{[I,\ka]}$ project to components whose virtual
indecomposables are of equivalence class~$[I,\ka]$.

\begin{prop} In the situation above $\Pi_{[I,\ka]}^2=\Pi_{[I,\ka]}$,
and\/ $\Pi_{[I,\ka]}\Pi_{[J,\la]}=0$ if\/ $[I,\ka]\ne[J,\la]$. Also,
if\/ $f\in\SFa(\fObj_\A)$ or $\oSFa(\fObj_\A,*,*)$ then
\e
f=\!\!\!
\sum_{\text{eq. classes $[I,\ka]$}}\!\!\!\!\!\!\!\!
\Pi_{[I,\ka]}(f) \;\>\text{and\/}\;\> \Pi^\vi_n(f)=\!\!\!
\sum_{\text{eq. classes $[I,\ka]:\md{I}=n$}} \!\!\!\!\!\!\!\!
\Pi_{[I,\ka]}(f),
\label{ai5eq5}
\e
where the sums make sense as only finitely many $\Pi_{[I,\ka]}(f)$
are nonzero.
\label{ai5prop4}
\end{prop}

\begin{proof} The analogue of \cite[Th.~5.12(c)]{Joyc2} for the
$\hat\Pi^\nu_\fF$ says that $\hat\Pi^{\nu_1}_\fF\ci
\hat\Pi^{\nu_2}_\fF=\hat\Pi^{\nu_1\nu_2}_\fF$. If $\nu_1,\nu_2$ are
the $\fObj_\A$-weight functions defined in Definition \ref{ai5def4}
using $[I,\ka]$ and $[J,\la]$ then $\nu_1\nu_2$ is $\nu_1$ if
$[I,\ka]=[J,\la]$ and 0 otherwise, so the first part follows. For
the second, let $[(\fR,\rho)]$ have algebra stabilizers, use the
notation of Proposition \ref{ai5prop3}, and define
$\Pi_{[I,\ka]}([([Y_l/A_l^\t],\rho_l)])$ using
\cite[Def.~5.15]{Joyc2}. This gives a finite sum over subtori
$P,Q,R$ in $A_l^\t$ with $(\K^\t)^{\md{I}}\cong R\subseteq P\cap Q$
of a term involving the subset $(Y_l)^{P,R}_{\nu,1}$ of $y\in
Y_l^P(\K)$ such that $R$ induces a decomposition of type $[I,\ka]$
of the image point in~$\fObj_\A(\K)$.

Each $y\in Y_l^P(\K)$ induces such a decomposition for a unique
$[I,\ka]$, and the map $y\mapsto[I,\ka]$ is constructible, and so
realizes only finitely many $[I,\ka]$ on $Y_l^P(\K)$. Hence
$\Pi_{[I,\ka]}([([Y_l/A_l^\t],\rho_l)])\ne 0$ for only finitely many
$[I,\ka]$, proving the last line. Summing over all $[I,\ka]$ yields
a sum over $P,Q,R$ of a term involving the whole of $Y_l^P(\K)$.
Comparing with \cite[Def.~5.10]{Joyc2} and using
\cite[Th.~5.12(a)]{Joyc2} gives
$\sum_{[I,\ka]}\Pi_{[I,\ka]}([([Y_l/A_l^\t],\rho_l)])=
\Pi^1([([Y_l/A_l^\t],\rho_l)])=[([Y_l/A_l^\t],\rho_l)]$. Restricting
to $\md{I}=n$ fixes $\dim R=n$, and the sum reduces to
$\Pi^\vi_n([([Y_l/A_l^\t],\rho_l)])$. Equation \eq{ai5eq5} follows.
\end{proof}

The $\Pi_{[I,\ka]}$ are also defined on $\uSF,\SF(\fObj_\A)$ and
$\uoSF,\oSF(\fObj_\A,*,*)$ with $\Pi_{[I,\ka]}^2\!=\! \Pi_{[I,\ka]}$
and $\Pi_{[I,\ka]}\Pi_{[J,\la]}\!=\!0$ if $[I,\ka]\!\ne\![J,\la]$,
but on these larger spaces \eq{ai5eq5} does not hold. Using
Proposition \ref{ai2prop2} and Corollary \ref{ai5cor2} we can show
\begin{equation*}
\ts\sum_{\text{eq. classes $[I,\ka]$}}
\Pi_{[I,\ka]}:\oSF(\fObj_\A,*,*)\longra\oSFa(\fObj_\A,*,*)
\end{equation*}
is a surjective projection. But the same is not true
for~$\SF(\fObj_\A),\SFa(\fObj_\A)$.

\subsection{The Lie algebras $\SFai(\fObj_\A)$ and
$\oSFai(\fObj_\A,*,*)$}
\label{ai52}

Next we study stack function analogues of the Lie algebra
$\CFi(\fObj_\A)$ of~\S\ref{ai44}.

\begin{dfn} Let Assumptions \ref{ai2ass} and \ref{ai3ass} hold.
Define $\SFai(\fObj_\A)$, $\oSFai(\fObj_\A,\Up,\La)$,
$\oSFai(\fObj_\A,\Up,\La^\ci)$ and $\oSFai(\fObj_\A,\Th,\Om)$ to be
the subspaces of $f\in\SFa(\fObj_\A)$ or $\oSFa(\fObj_\A,*,*)$
satisfying~$\Pi^\vi_1(f)=f$.
\label{ai5def5}
\end{dfn}

We interpreted $\Pi^\vi_n$ above as projecting to $f$ `supported on
objects with $n$ virtual indecomposable factors'. So
$\SFai(\fObj_\A)$, $\oSFai(\fObj_\A,*,*)$ should be thought of as
stack functions {\it supported on virtual indecomposables}, and are
good analogues of $\CFi(\fObj_\A)$. Our goal is to prove that
$\SFai(\fObj_\A)$, $\oSFai(\fObj_\A,*,*)$ are {\it Lie subalgebras}
of $\SFa(\fObj_\A)$, $\oSFa(\fObj_\A,*,*)$. To do this we must study
the relationship between multiplication $*$ and projections
$\Pi^\vi_n$, that is, express $\Pi^\vi_n(f*g)$ in terms of
$\Pi^\vi_l(f)$ and~$\Pi^\vi_m(g)$.

\begin{prop} Let\/ $T\subseteq\fObj_\A(\K)\t\fObj_\A(\K)$ be
constructible. Then there exist a finite decomposition
$T=\coprod_{m\in M}\fG_m(\K)$, where $\fG_m$ is a finite type
$\K$-substack of\/ $\fObj_\A\t\fObj_\A$, $1$-isomorphisms
$\fG_m\cong[V_m/G_m]$ for $G_m$ a special\/ $\K$-group and\/ $V_m$ a
quasiprojective $\K$-variety, finite-dimensional representations
$E_m^0,E_m^1$ of\/ $G_m$, and morphisms $\jmath_m:(\K^\t)^2\!\ra\!
G_m$ for all\/ $m\in M$, satisfying:
\begin{itemize}
\setlength{\itemsep}{0pt}
\setlength{\parsep}{0pt}
\item[{\rm(a)}] Let\/ $v\in V_m(\K)$ project to $([X],[Y])\!\in\!
\fG_m(\K)\!\subset\!\fObj_\A(\K)\!\t\!\fObj_\A(\K)$,~so
\e
\Stab_{G_m}(v)\cong\Iso_\K([X],[Y])\cong\Aut(X)\t\Aut(Y).
\label{ai5eq6}
\e
Then there are isomorphisms $E^0_m\cong\Hom(Y,X)$ and\/
$E^1_m\cong\Ext^1(Y,X)$ such that\/ \eq{ai5eq6} identifies the
action of\/ $\Stab_{G_m}(v)$ on $E^i_m$ with the action of\/
$\Aut(X)\!\t\!\Aut(Y)$ on $\Hom(Y,X),\Ext^1(Y,X)$ given by
$(\al,\be)\cdot e\!=\!\al\ci e\ci\be^{-1}$.
\item[{\rm(b)}] $\jmath_m$ maps into the centre of\/ $G_m$, and\/
$\jmath_m\bigl((\K^\t)^2\bigr)$ acts freely on $V_m$. Thus, in
{\rm(a)} $\jmath_m$ maps $(\K^\t)^2\ra\Stab_{G_m}(v)$. Composing
this with\/ \eq{ai5eq6} gives the map
$(\de,\ep)\mapsto(\de\id_X,\ep\id_Y)$, for~$\de,\ep\in\K^\t$.
\item[{\rm(c)}] Write $i_m:\fG_m\ra\fObj_\A\t\fObj_\A$ for
the inclusion $1$-morphism. Then there is a $1$-isomorphism
\e
{}\!\!\!\!\!\!\!\!\!\!\!\!\!\!
\fG_m\!\t_{i_m,\fObj_\A\!\t\!\fObj_\A,\bs\si(\{1\})\!\t\bs\si(\{2\})}
\!\fM(\{1,2\},\le)_\A\!\cong\!\bigl[V_m\!\t\!E^1_m/G_m\!\lt\!
E^0_m\bigr].
\label{ai5eq7}
\e
Here multiplication on $G_m\lt E^0_m$ is
$(\ga,e)\cdot(\ga',e')=(\ga\ga',e+\ga\cdot e')$, and\/ $E^0_m$ acts
trivially on $V_m\t E^1_m$, and\/ $G_m$ acts in the given way on
$V_m,E^1_m$.
\item[{\rm(d)}] Equation \eq{ai5eq7} is a substack of\/
$\fM(\{1,2\},\le)_\A$, so its $\K$-points are $[(\si,\io,\pi)]$ for
$(\{1,2\},\le)$-configurations $(\si,\io,\pi)$. Let\/ $X,Y,v$ be as
in {\rm(a)}, and\/ $(v,e)\in(V_m\t E^1_m)(\K)$ project to
$[(\si,\io,\pi)]$. Then we can choose $\si(\{1\})=X$,
$\si(\{2\})=Y$, and\/ $e\in\Ext^1(Y,X)$ corresponds to the exact
sequence
\e
\xymatrix@C=22pt{ 0 \ar[r] &\si(\{1\}) \ar[rr]^{\io(\{1\},\{1,2\})}
&&\si(\{1,2\}) \ar[rr]^{\pi(\{1,2\},\{2\})} &&\si(\{2\}) \ar[r] & 0.
}
\label{ai5eq8}
\e
\end{itemize}
\label{ai5prop5}
\end{prop}

\begin{proof} The proof of Proposition \ref{ai5prop2} easily
generalizes to give a finite decomposition $T=\coprod_{m\in
M}\fG_m(\K)$ and 1-isomorphisms $\fG_m\cong[U_m/A_m^\t\t B_m^\t]$,
for $A_m,B_m$ finite-dimensional $\K$-algebras, such that if $u\in
U_m(\K)$ projects to $([X],[Y])$ then $\Stab_{A_m^\t\t B_m^\t}(u)=
C_u^\t\t D_u^\t$ for subalgebras $C_u\subseteq A_m$, $D_u\subseteq
B_m$ with isomorphisms $C_u\cong\End(X)$, $D_u\cong\End(Y)$ inducing
the isomorphism~$\Stab_{A_m^\t\t B_m^\t}(u)\cong\Aut(X)\t\Aut(Y)$.

The functions $([X],[Y])\mapsto\dim\Hom(Y,X)$ or $\dim\Ext^1(Y,X)$
are locally constructible on $(\fObj_\A\t\fObj_\A)(\K)$, and so take
finitely many values on $T$. Refining $T=\coprod_{m\in M}\fG_m(\K)$,
we can make $\dim\Hom(Y,X),\dim\Ext^1(Y,X)$ constant on each
$\fG_m(\K)$. Refining further, $\Hom(Y,X),\Ext^1(Y,X)$ are the
fibres over $([X],[Y])\in\fG_m(\K)$ of {\it vector bundles} over
$\fG_m$, in the stack sense. These pull back under the projection
$U_m\ra\fG_m$ to vector bundles $\mathcal{E}^0_m,\mathcal{E}^1_m$
over $U_m$, with fibres $\K$-vector spaces $E^0_m,E^1_m$, and the
$A_m^\t\t B_m^\t$-action on $U_m$ lifts to actions on
$\mathcal{E}^0_m,\mathcal{E}^1_m$ preserving the vector bundle
structure.

Define $V_m$ to be the quasiprojective $\K$-variety of triples
$(u,\al^0,\al^1)$, for $u\in U_m(\K)$ and
$\al^i:(\mathcal{E}^i_m)_v\ra E^i_m$ vector space isomorphisms
between the fibre of $\mathcal{E}^i_m$ over $v$ for $i=0,1$ and
$E^i_m$. Then $V_m$ is a principal bundle over $\ti V_m$ with
structure group $\Aut(E^0_m)\t\Aut(E^1_m)$, and the $A_m^\t\t
B_m^\t$-action lifts to $V_m$ and commutes with the
$\Aut(E^0_m)\t\Aut(E^1_m)$-action. Define $G_m=A_m^\t\t
B_m^\t\t\Aut(E^0_m)\t\Aut(E^1_m)$, which is special as it is a
product of groups of the form $A^\t$ for finite-dimensional algebras
$A$. It acts on $V_m$ with $[V_m/G_m]\cong[U_m/A_m^\t\t
B_m^\t]\cong\fG_m$. Define actions of $G_m$ on $E_m^i$ for $i=0,1$
by $(a,b,\be^0,\be^1)\cdot e^i=\be^ie^i$. It is now easy to see that
(a) holds for~$V_m,G_m$.

For (b), define $\jmath_m(\de,\ep)=(\de\id_{A_m},\ep\id_{B_m},
\de\ep^{-1}\id_{E^0_m},\de\ep^{-1}\id_{E^1_m})$. This is clearly a
$\K$-group morphism to the centre of $G_m$. If $v=(u,\al^0,\al^1)$
in $V_m(\K)$ then $\Stab_{A_m^\t\t B_m^\t}(u)=C_u^\t\t D_u^\t$ as
above, and $\de\id_{A_m}\in C_u^\t$, $\ep\id_{B_m}\in D_u^\t$ as
$C_u,D_u$ are subalgebras, so $(\de\id_{A_m},\ep\id_{B_m})$ fixes
$u$. The identification of actions in (a) then shows
$(\de\id_{A_m},\ldots,\de\ep^{-1} \id_{E^1_m})$ fixes $v$. Thus,
$\jmath_m((\K^\t)^2)$ fixes each $v\in V_m(\K)$, and acts freely on
$V_m$. Composing with \eq{ai5eq6} gives
$(\de,\ep)\mapsto(\de\id_X,\ep\id_Y)$, as
$\de\id_{A_m}=\de\id_{C_u}$, $\ep\id_{B_m}=\ep\id_{D_u}$ are
identified with $\de\id_X,\ep\id_Y$ under $C_u\cong\End(X)$,
$D_u\cong\End(Y)$.

For (c) and (d), note that the fibre of
$\bs\si(\{1\})\t\bs\si(\{2\}):\fM(\{1,2\},\le)_\A
\ra\fObj_\A\t\fObj_\A$ over $([X],[Y])\in(\fObj_\A\t\fObj_\A)(\K)$
is the family of $[(\si,\io,\pi)]$ for $(\si,\io,\pi)$ a
$(\{1,2\},\le)$-configuration with $\si(\{1\})\cong X$ and
$\si(\{2\})\cong Y$. This is equivalent to a short exact sequence
\eq{ai5eq8}, which are classified by $\Ext^1(Y,X)$. But as the
isomorphisms $X\!\cong\!\si(\{1\})$, $Y\!\cong\!\si(\{2\})$ are not
prescribed we divide by $\Aut(X)\!\t\!\Aut(Y)$, so as sets we have
\begin{equation*}
\M(\{1,2\},\le)_\A\supset\bigl(\bs\si(\{1\})\t\bs\si(\{2\})
\bigr)_*^{-1}\bigl(([X],[Y])\bigr)\cong
\frac{\Ext^1(Y,X)}{\Aut(X)\t\Aut(Y)}.
\end{equation*}

To describe this fibre as a stack we must take stabilizer groups
into account. One can show that if \eq{ai5eq8} corresponds to
$e\in\Ext^1(Y,X)$ then $\Aut\bigl((\si,\io,\pi)\bigr)\cong H_e
\lt\Hom(Y,X)$, where $H_e$ is the $\K$-subgroup of
$\Aut(X)\t\Aut(Y)$ fixing $e$. If $([X],[Y])$ is the image of $v\in
V_m(\K)$, we have 1-isomorphisms
\begin{gather*}
\Spec\K\t_{X\t Y,\fObj_\A\t\fObj_\A,\bs\si(\{1\})\t\bs\si(\{2\})}
\fM(\{1,2\},\le)_\A\cong\\
\bigl[\Ext^1(Y,X)/(\Aut(X)\t\Aut(Y))\lt\Hom(Y,X)\bigr]
\cong\bigl[E^1_m/\Stab_{G_m\lt E^0_m}(v)\bigr].
\end{gather*}
Here $\Hom(Y,X)$ acts trivially on $\Ext^1(Y,X)$, but contributes to
the stabilizers, and similarly $E^0_m$ acts trivially on $V_m$.
Parts (c) and (d) follow by a families version of this argument.
\end{proof}

\begin{cor} Let\/ $f,g$ lie in $\SFa(\fObj_\A)$ or $\oSFa(\fObj_\A,*,*)$,
choose $T\subseteq\fObj_\A(\K)\t\fObj_\A(\K)$ constructible with\/
$f\ot g$ supported on $T$, and use the notation of Proposition
\ref{ai5prop5}. Then arguing as in Proposition \ref{ai5prop3}, we
may write
\e
f\ot g=\ts\sum_{m\in M,\; n\in N_m}c_{mn}\bigl[\bigl([W_{mn}/
G_m],\tau_{mn}\bigr)\bigr],
\label{ai5eq9}
\e
where $N_m$ is finite, $c_{mn}\in\Q,\La,\La^\ci$ or $\Om$, $W_{mn}$
is a quasiprojective $\K$-variety acted on by $G_m$, and\/
$\tau_{mn}:[W_{mn}/G_m]\ra\fG_m\subseteq\fObj_\A\t\fObj_\A$ is
induced by a $G_m$-equivariant morphism $\phi_{mn}:W_{mn}\ra V_m$.
Moreover
\e
\bigl(\bs\si(\{1\})\!\t\!\bs\si(\{2\})\bigr)^*
\bigl([([W_{mn}/G_m],\tau_{mn})]\bigr)\!=\!\bigl[\bigl([W_{mn}
\!\t\!E^1_m/G_m\!\lt\!E^0_m],\xi_{mn}\bigr)\bigr]
\label{ai5eq10}
\e
in $\SF(\fM(\{1,2\},\le)_\A)$ or $\oSF(\fM(\{1,2\},\le)_\A,*,*)$,
where $\phi_{mn}$ induces
\e
\xi_{mn}:[W_{mn}\t E^1_m/G_m\lt E^0_m]\longra\bigl[V_m\t
E^1_m/G_m\lt E^0_m\bigr],
\label{ai5eq11}
\e
using \eq{ai5eq7} to regard the right hand side of\/ \eq{ai5eq11} as
a substack of\/ $\fM(\{1,2\},\le)_\A$. Combining \eq{ai5eq1},
\eq{ai5eq9} and\/ \eq{ai5eq10} gives
\e
f*g=\ts\sum_{m\in M,\; n\in N_m}c_{mn}\bigl[\bigl([W_{mn}\t
E^1_m/G_m\lt E^0_m],\bs\si(\{1,2\})\ci\xi_{mn}\bigr)\bigr].
\label{ai5eq12}
\e
\label{ai5cor3}
\end{cor}

Here we have used a formula for the fibre product of quotient stacks
from the proof of \cite[Th.~4.12]{Joyc2} to deduce \eq{ai5eq10}. Our
next theorem, which will be important in \cite[\S 8]{Joyc4}, proves
a relationship between the operators $\Pi^\vi_k$ and~$P_\sIp$.

\begin{thm} Let Assumptions \ref{ai2ass} and \ref{ai3ass} hold,
$(I,\tl)$ be a finite poset, $k\ge 0$, and\/ $f_i$ for $i\in I$ lie
in $\SFai(\fObj_\A)$ or $\oSFai(\fObj_\A,*,*)$. Then
\e
\begin{split}
&\Pi^\vi_k\bigl[P_\sIt(f_i:i\in I)\bigr]=\\
&\sum_{\substack{\text{iso. classes}\\ \text{of finite sets}\\
\text{$K$, $k\le\md{K}\le\md{I}$}}}\,\,
\sum_{\substack{\text{surjective $\phi:I\ra K$.}\\
\text{Define $\pr$ on $I$ by $i\pr j$}\\
\text{if $i\tl j$ and $\phi(i)=\phi(j)$}}}\!\!\! N_{I,K,\phi,k}\cdot
P_\sIp(f_i:i\in I),
\end{split}
\label{ai5eq13}
\e
where $N_{I,K,\phi,k}\in\Q$ depends only on $I,K,\phi$ up to
isomorphism and\/~$k$.
\label{ai5thm3}
\end{thm}

\begin{proof} It is easy to partially generalize Proposition
\ref{ai5prop5}(a),(b) from $(\{1,2\},\le)$ to $(I,\pr)$ as follows.
Choose $T\subseteq\prod_{i\in I}\fObj_\A(\K)$ constructible with
$\bigot_{i\in I}f_i$ supported on $T$. Then there exist a finite
decomposition $T=\coprod_{m\in M}\fG_m(\K)$, where $\fG_m$ is a
finite type $\K$-substack of $\prod_{i\in I}\fObj_\A$,
1-isomorphisms $\fG_m\cong[V_m/G_m]$ for $G_m$ a special $\K$-group
and $V_m$ a quasiprojective $\K$-variety, and morphisms
$\jmath_m:(\K^\t)^I\ra G_m$ for all $m\in M$, satisfying:
\begin{itemize}
\setlength{\itemsep}{0pt}
\setlength{\parsep}{0pt}
\item[(a)] Let $v\in V_m(\K)$ project to $\prod_{i\in I}[X_i]
\in\fG_m(\K)\subset\prod_{i\in I}\fObj_\A(\K)$, so that
\e
\ts\Stab_{G_m}(v)\cong\Iso_\K\bigl(\prod_{i\in
I}[X_i]\bigr)\cong\prod_{i\in I}\Aut(X_i).
\label{ai5eq14}
\e
\item[(b)] $\jmath_m$ maps into the centre of $G_m$, and
$\jmath_m\bigl((\K^\t)^I\bigr)$ acts freely on $V_m$. Thus, in (a)
$\jmath_m$ maps $(\K^\t)^I\ra\Stab_{G_m}(v)$. Composing this with
\eq{ai5eq14} gives the map $\de\mapsto\prod_{i\in
I}\de(i)\id_{X_i}$, for~$\de\in(\K^\t)^I$.
\end{itemize}

Generalizing part (c) is more tricky. Write $i_m:\fG_m\ra
\prod_{i\in I}\fObj_\A$ for the inclusion 1-morphism. Then we can
form $\fG_m\t_{i_m,\prod_{i\in I}\fObj_\A,\prod_{i\in
I}\bs\si(\{i\})}\fM(I,\tl)_\A$, regarded as a $\K$-substack of
$\fM(I,\tl)_\A$. It is of {\it finite type}, as $\fG_m$ is and
$\prod_{i\in I}\bs\si(\{i\})$ is by Theorem \ref{ai3thm}(c). Also
$\fM(I,\tl)_\A$ has affine geometric stabilizers. So using
\cite[Prop.~3.5.9]{Kres} as in Proposition \ref{ai5prop2}, we may
write
\begin{equation*}
\bigl(\fG_m\t_{i_m,\prod_{i\in I}\fObj_\A,\prod_{i\in
I}\bs\si(\{i\})}\fM(I,\tl)_\A\bigr)(\K)=\ts\coprod_{p\in
P_m}\fH_{mp}(\K),
\end{equation*}
where $P_m$ is finite and $\fH_{mp}$ a $\K$-substack of
$\fM(I,\tl)_\A$ with a 1-isomorphism
\e
\fH_{mp}\cong[Y_{mp}/G_m\lt K_{mp}],
\label{ai5eq15}
\e
where $Y_{mp}$ is a quasiprojective $\K$-variety and $K_{mp}$ a {\it
nilpotent\/} $\K$-group acted on by $G_m$, such that $\prod_{i\in
I}\bs\si(\{i\}):\fH_{mp}\ra\fG_m$ is induced by a morphism
$\psi_{mp}:Y_{mp}\ra V_m$ equivariant w.r.t.\ the natural
projection~$G_m\lt K_{mp}\ra G_m$.

The only nontrivial claim here is that we can take the quotient
group in \eq{ai5eq15} to be $G_m\lt K_{mp}$ with $K_{mp}$ {\it
nilpotent}, rather than an arbitrary $\K$-group with a morphism to
$G_m$. When $(I,\tl)=(\{1,2\},\le)$ this follows from Proposition
\ref{ai5prop5}(c) with $K_{mp}=E^0_m$. The general case can be
proved by the inductive argument on $\md{I}$ in the proof of Theorem
\ref{ai6thm2} below, which builds up $\fM(I,\tl)_\A$ by repeated
fibre products with $\bs\si(\{1\})\t\bs\si(\{2\}):\fM(\{1,2\})_\A
\ra\fObj_\A\t\fObj_\A$. The point of this is that as $K_{mp}$ is
nilpotent we can use the same maximal torus for $G_m$ and $G_m\lt
K_{mp}$, which will be important when we come to apply~$\Pi^\vi_k$.

Now we can generalize Corollary \ref{ai5cor3} to write
\e
\ts\bigot_{i\in I}f_i=\ts\sum_{m\in M,\; n\in N_m}c_{mn}
\bigl[\bigl([W_{mn}/G_m],\tau_{mn}\bigr)\bigr],
\label{ai5eq16}
\e
where $N_m$ is finite, $c_{mn}\in\Q,\La,\La^\ci$ or $\Om$, $W_{mn}$
is a quasiprojective $\K$-variety acted on by $G_m$, and
$\tau_{mn}:[W_{mn}/G_m]\ra\fG_m\subseteq\fObj_\A\t\fObj_\A$ is
induced by a $G_m$-equivariant morphism $\phi_{mn}:W_{mn}\ra V_m$.
Moreover
\e
\begin{split}
&\ts\bigl(\prod_{i\in I}\bs\si(\{i\})\bigr)^*
\bigl([([W_{mn}/G_m],\tau_{mn})]\bigr)=\\
&\ts\sum_{p\in P_m} \bigl[\bigl([W_{mn}
\t_{\phi_{mn},V_m,\psi_{mp}}Y_{mp}/G_m\lt
K_{mp}],\xi_{mnp}\bigr)\bigr]
\end{split}
\label{ai5eq17}
\e
in $\SF(\fM(I,\tl)_\A)$ or $\oSF(\fM(I,\tl)_\A,*,*)$, where
$\pi_{Y_{mp}}$ induces
\e
\xi_{mnp}:\bigl[W_{mn}\t_{\phi_{mn},V_m,\psi_{mp}}Y_{mp}/G_m\lt
K_{mp}\bigr]\longra\bigl[Y_{mp}/G_m\lt K_{mp}\bigr],
\label{ai5eq18}
\e
using \eq{ai5eq15} to regard the right hand side of \eq{ai5eq18} as
a substack of $\fM(I,\tl)_\A$. Combining \eq{ai5eq3}, \eq{ai5eq16}
and \eq{ai5eq17} gives
\e
P_\sIt(f_i:i\in I)=\!\!\sum_{\substack{m\in M,\\ n\in N_m,\; p\in
P_m}}\!\!\!\! \begin{aligned}[t] c_{mn}\bigl[\bigl([W_{mn}
\t_{\phi_{mn},V_m,\psi_{mp}}Y_{mp}/G_m\lt K_{mp}],&\\
\bs\si(I)\ci\xi_{mnp}\bigr)\bigr].&\end{aligned}
\label{ai5eq19}
\e

Since $\Pi^\vi_1(f_i)=f_i$, we deduce from
\cite[Prop.~5.14(iv)]{Joyc2} that
\e
\Pi^\vi_k\bigl[\ts\bigot_{i\in I}f_i\bigr]=\begin{cases}
\bigot_{i\in I}f_i, & \text{$k=\md{I}$,} \\ 0, & \text{otherwise.}
\end{cases}
\label{ai5eq20}
\e
Let $T_m$ be a maximal torus in $G_m$, so that $T_m\t\{0\}=T_m$ is a
maximal torus in $G_m\lt K_{mp}$. Applying \cite[Def.s 5.10 \&
5.13]{Joyc2} to \eq{ai5eq16} we find that
\e
\Pi^\vi_k\bigl[{\ts\bigot_{i\in I}}f_i\bigr]\!=
\!\!\!\!\!\!\!\!\!\!\!\!\!\!\!\!\!\!\!\!
\sum_{\begin{subarray}{l}
\text{$m\in M$, $n\in N_m$, $P\!\in\!\cP(W_{mn},T_m)$,}\\
\text{$Q\!\in\!\cQ(G_m,T_m)$, $R\!\in\!\cR(W_{mn},G_m,T_m):$}\\
\text{$R\!\subseteq\!P\cap Q$, $M^{W_{mn}}_{G_m}(P,Q,R)\!\ne\!0$,
$\dim R\!=\!k$}
\end{subarray}}
\!\!\!\!\!\!\!\!\!\!\!\!\!\!\!\!\!\!
\begin{aligned}[t]
&c_{mn}\,M^{W_{mn}}_{G_m}(P,Q,R)\,\cdot\\
&\qquad \bigl[\bigl([W_{mn}^P/C_{G_m}(Q)],
\tau_{mn}\!\ci\!\io^{P\cap Q}\bigr)\bigr].
\end{aligned}
\label{ai5eq21}
\e

Let $m'\in M$ and $R'\subseteq T_{m'}$ lie in $\cR(W_{m'n},G_{m'},
T_{m'})$ for at least one $n\in N_{m'}$. By \cite[Def.~5.15]{Joyc2}
we can define $\hat\Pi^\nu_{\Pi_{i\in I}\fObj_\A}$ on $\SF(\Pi_{i\in
I}\fObj_\A)$ or $\oSF(\Pi_{i\in I}\fObj_\A,*,*)$ with the $\Pi_{i\in
I}\fObj_\A$-weight function $\nu$ given by $\nu(T,g,\phi)=1$ if
$\phi$ is injective and $\phi_*(T)\subset\Iso_\K(g)$ is identified
with $R'\subseteq \Stab_{G_{m'}}(w)$ for some $n\in N_{m'}$ and
$w\in W_{m'n}(\K)$ projecting to
$g\in\fG_{m'}(\K)\subseteq(\Pi_{i\in I}\fObj_\A)(\K)$ under
$\fG_{m'}\cong[W_{m'n}/G_{m'}]$, and $\nu(T,g,\phi)=0$ otherwise.

Applying $\hat\Pi^\nu_{\Pi_{i\in I}\fObj_\A}$ to \eq{ai5eq21}
projects to those components in the sum with $m=m'$ and $R$
conjugate to $R'$ under the Weyl group $W$ of $G_m$. By
\eq{ai5eq20}, equation \eq{ai5eq21} is zero for $k\ne\md{I}$. So, if
$\dim R'\ne\md{I}$ then the sum of components in \eq{ai5eq21} with
$m=m'$ and $R$ conjugate to $R'$ is zero. But by symmetry in $W$
each of the conjugates of $R'$ give the same answer. Hence, for any
fixed $m\in M$ and $R$ with $\dim R\ne\md{I}$, the sum of components
in \eq{ai5eq21} with these $m,R$ is zero.

By (b) above, $\jmath_m((\K^\t)^I)$ lies in the centre of $G_m$ and
acts trivially on $V_m$, and we can use algebra stabilizers to show
that it also acts trivially on each $W_{mn}$. Thus
$\jmath_m((\K^\t)^I)\subseteq R$ for each $R\!\in\!\cR(W_{mn},
G_m,T_m)$. Also, as $G_m$ is special and each $f_i$ has algebra
stabilizers one can show that each such $R$ is a torus (rather than
the product of a torus with a finite group). And as each $f_i$ is
supported over {\it nonzero} elements $[X_i]$, as composing
$\jmath_m$ with \eq{ai5eq14} takes $\de\mapsto(\de(i)\id_{X_i}
)_{i\in I}$, we see $\jmath_m$ is {\it injective}. Therefore
$\jmath_m((\K^\t)^I)\cong (\K^\t)^{\md{I}}$ is the minimal element
of $\cR(W_{mn},G_m,T_m)$, and the unique element $R$ with $\dim
R=\md{I}$. By \eq{ai5eq16} and \eq{ai5eq20}, the sum of terms in
\eq{ai5eq21} with these $m,R$ is~$\sum_{n\in N_m}c_{mn}\bigl[
\bigl([W_{mn}/G_m],\tau_{mn}\bigr)\bigr]$.

As for \eq{ai5eq21}, equation \eq{ai5eq19} yields an expression for
$\Pi^\vi_k\bigl[P_\sIt(f_i:i\!\in\!I)\bigr]$, but it is not in the
form we want. Applying $\Pi^\vi_k$ to $[([W_{mn}\t_{V_m}Y_{mp}
/G_m\lt K_{mp}],\bs\si(I)\ci\xi_{mnp})]$ involves summing over the
three finite sets:
\begin{align*}
&\cP(W_{mn}\t_{V_m}Y_{mp},T_m)\subseteq\bigl\{P\cap\dot P:
\text{$P\in\cP(W_{mn},T_m)$, $\dot P\in\cP(Y_{mp},T_m)$}\bigr\}\\
&\cQ(G_m\lt K_{mp},T_m)=\bigl\{Q\cap\dot Q:
\text{$Q\in\cQ(G_m,T_m)$, $\dot Q\in\cQ(G_m\lt K_{mp},T_m)$}\bigr\}\\
&\cR(W_{mn}\t_{V_m}Y_{mp},G_m\lt
K_{mp},T_m)\subseteq\bigl\{R\cap\dot R:
\text{$R\in\cR(W_{mn},G_m,T_m)$,}\\
&\qquad\qquad\qquad\qquad\qquad\qquad\qquad\qquad\quad \text{$\dot
R\in\cR(Y_{mp},G_m\lt K_{mp},T_m)$}\bigr\}.
\end{align*}
It follows from the proof in \cite[\S 5]{Joyc2} that the definition
of $\Pi^\vi_k$ is independent of choices, that if we replace
$\cP(W_{mn}\t_{V_m}Y_{mp},T_m)$ by $\cP(W_{mn},T_m)\t
\cP(Y_{mp},T_m)$, and $\cQ(G_m\lt K_{mp},T_m)$ by $\cQ(G_m,T_m)\t
\cQ(G_m\lt K_{mp},T_m)$, and $\cR(W_{mn}\t_{V_m}Y_{mp},G_m\lt
K_{mp},T_m)$ by $\cR(W_{mn},G_m,T_m)\t\cR(Y_{mp},G_m\lt K_{mp},T_m)$
throughout the definition, taking $(P,\dot P)$ to act as $P\cap\dot
P$ on $W_{mn}\!\t_{V_m}\!Y_{mp}$ and so on, we get the same answer.
But defined using these sets we easily find that
\begin{equation*}
M^{W_{mn}\!\t_{V_m}\!Y_{mp}}_{G_m\!\lt\!K_{mp}}\bigl((P,\dot
P),(Q,\dot Q),(R,\dot R)\bigr)= M^{W_{mn}}_{G_m}(P,Q,R)\cdot
M^{Y_{mp}}_{G_m\!\lt\!K_{mp}}(\dot P,\dot Q,\dot R).
\end{equation*}
This yields:
\e
\begin{gathered}
\Pi^\vi_k\bigl[P_\sIt(f_i:i\!\in\! I)\bigr]\!=
\!\!\!\!\!\!\!\!\!\!\!\!\!\!\!\!\!\!\!\!\!\!\!\!\!
\sum_{\begin{subarray}{l}
\text{$m\!\in\!M$, $p\in P_m$, $\dot P\in\cP(Y_{mp},T_m)$,}\\
\text{$\dot Q\in\cQ(G_m\lt K_{mp},T_m)$, $\dot R\in\cR(Y_{mp},G_m\lt
K_{mp},T_m)$:}\\
\text{$\dot R\!\subseteq\!\dot P\!\cap\!\dot Q$,
$M^{Y_{mp}}_{G_m\!\lt\!K_{mp}}(\dot P,\dot Q,\dot R)\ne 0$}
\end{subarray}}
\!\!\!\!\!\!\!\!\!\!\!\!\!\!\!\!\!\!\!\!\!\!\!\!\!\!\!\!\!\!\!
M^{Y_{mp}}_{G_m\!\lt\!K_{mp}}(\dot P,\dot Q,\dot R)\cdot
\\
\raisebox{-16pt}{\begin{Large}$\displaystyle\Biggl[$\end{Large}}\,
\sum_{\begin{subarray}{l}
\text{$n\in N_{nm}$, $P\!\in\!\cP(W_{mn},T_m)$,}\\
\text{$Q\!\in\!\cQ(G_m,T_m)$,}\\
\text{$R\!\in\!\cR(W_{mn},G_m,T_m)$:}\\
\text{$R\!\subseteq\!P\!\cap\!Q$, $M^{W_{mn}}_{G_m}(P,Q,R)\!\ne\!0$,
$\dim R\!\cap\!\dot R\!=\!k$}
\end{subarray}
\!\!\!\!\!\!\!\!\!\!\!\!\!\!\!\!\!\!\!\!\!\!\!\!\!\!\!\!\! }
\begin{aligned}[t]
&c_{mn}\,M^{W_{mn}}_{G_m}(P,Q,R)\,\cdot\\[-3pt]
&\bigl[\bigl([W_{mn}^P\!\t_{V_m}\!Y_{mp}^{\dot P}/
C_{G_m}(Q)\!\lt\!(K_{mp})^{\dot Q}],\\[-3pt]
&\qquad\qquad\qquad \bs\si(I)\!\ci\!\xi_{mnp}\!\ci\!\io^{P\cap
Q,\dot P\cap\dot Q}\bigr)\bigr]
\end{aligned}
\,\raisebox{-16pt}{\begin{Large}$\displaystyle\Biggr]$\end{Large}}.
\end{gathered}
\label{ai5eq22}
\e

Actually, according to the argument we gave the last term should be
\begin{equation*}
\bigl[\bigl([W_{mn}^{P\cap\dot P}\!\t_{V_m}\!Y_{mp}^{P\cap\dot P}/
C_{G_m}(Q\cap\dot Q)\!\lt\!(K_{mp})^{Q\cap\dot Q}],
\bs\si(I)\!\ci\!\xi_{mnp}\!\ci\!\io^{P\cap\dot P\cap Q\cap\dot
Q}\bigr)\bigr].
\end{equation*}
However, using \cite[Lem.~5.9]{Joyc2} we can show that for fixed
$P,\dot P,Q,\dot Q$ in \eq{ai5eq22}, unless $P,\dot P$ are the
smallest elements of $\cP(W_{mn},T_m)$, $\cP(Y_{mp},T_m)$ containing
$P\cap\dot P$, and $Q,\dot Q$ are the smallest elements of
$\cQ(G_m,T_m),\cQ(G_m\lt K_{mp},T_m)$ containing $Q\cap\dot Q$, then
the sum of $M^{Y_{mp}}_{G_m\!\lt\!K_{mp}}(\dot P,\dot Q,\dot R)\cdot
M^{W_{mn}}_{G_m}(P,Q,R)$ over all $R,\dot R$ with fixed $R\cap\dot
R$ is zero. But if $P,\dot P,Q,\dot Q$ are the smallest elements
containing $P\cap\dot P,Q\cap\dot Q$ then $W_{mn}^P=W_{mn}^{P\cap
\dot P}$, $Y_{mp}^{\dot P}=Y_{mp}^{P\cap\dot P}$, $C_{G_m}(Q)=
C_{G_m}(Q\cap\dot Q)$ and $(K_{mp})^{\dot Q}=(K_{mp})^{Q\cap\dot
Q}$, by \cite[Lem.s~5.4(iii) \& 5.6(iii)]{Joyc2}. So \eq{ai5eq22} is
correct.

The important point about the way we have written \eq{ai5eq22} is
that the last line $[\cdots]$ is a sum over certain $R$ of a linear
operation applied to the terms in \eq{ai5eq21} with fixed $m,R$. But
we have already shown that the sum of terms in \eq{ai5eq21} with
fixed $m,R$ is $\sum_{n\in N_m}c_{mn}\bigl[\bigl([W_{mn}/G_m],
\tau_{mn}\bigr)\bigr]$ if $R=\jmath_m((\K^\t)^I)$, and 0 otherwise.
This proves that
\e
\Pi^\vi_k\bigl[P_\sIt(f_i:i\!\in\! I)\bigr]\!=
\!\!\!\!\!\!\!\!\!\!\!\!\!\!\!\!\!\!\!\!\!\!\!\!\!\!\!\!
\!\!\!\!\!\!\!\!\!\!\!\!\!\!
\sum_{\begin{subarray}{l}
\text{$m\!\in\!M$, $n\in N_{nm}$, $p\in P_m$,}\\
\text{$\dot P\!\in\!\cP(Y_{mp},T_m)$, $\dot Q\!\in\!\cQ(G_m\!\lt\!
K_{mp},T_m)$,}\\
\text{$\dot R\!\in\!\cR(Y_{mp},G_m\!\lt\!K_{mp},T_m):\dot
R\!\subseteq\!\dot P\!\cap\!\dot Q$,}\\
\text{$M^{Y_{mp}}_{G_m\!\lt\!K_{mp}}(\dot P,\dot Q,\dot R)\ne 0$,
$\dim\jmath_m((\K^\t)^I)\!\cap\!\dot R\!=\!k$}
\end{subarray}}
\!\!\!\!\!\!\!\!\!\!\!\!\!
\begin{aligned}[t]
&M^{Y_{mp}}_{G_m\!\lt\!K_{mp}}(\dot P,\dot Q,\dot R)\cdot c_{mn}\cdot\\
&\bigl[\bigl([W_{mn}\!\t_{V_m}\!Y_{mp}^{\dot P}/
G_m\!\lt\!(K_{mp})^{\dot Q}],\!\!\!\!\!\! \\
&\qquad\qquad\qquad \bs\si(I)\!\ci\!\xi_{mnp}\bigr)\bigr].
\end{aligned}
\label{ai5eq23}
\e

Now using the argument of \cite[Lem.~5.9]{Joyc2} again, we find that
unless $\dot P,\dot Q\subseteq\jmath_m((\K^\t)^I)$, the sum of over
$\dot R$ of terms with these $\dot P,\dot Q$ in \eq{ai5eq23} is
zero. So we may restrict to $\dot P,\dot Q,\dot R\subseteq\jmath_m
((\K^\t)^I)$, and then as $\jmath_m((\K^\t)^I)$ lies in the centre
of $G_m$ we may replace $G_m,T_m$ by $\jmath_m((\K^\t)^I)$ in the
definitions of $\dot P,\dot Q,\dot R$ and $M^{Y_{mp}}_{G_m\!\lt
\!K_{mp}}(\dot P,\dot Q,\dot R)$, yielding:
\e
\Pi^\vi_k\bigl[P_\sIt(f_i:i\!\in\! I)\bigr]\!=
\!\!\!\!\!\!\!\!\!\!\!\!\!\!\!\!\!\!\!\!\!\!\!\!\!\!\!\!
\!\!\!\!\!\!\!\!\!\!\!\!\!\!
\sum_{\begin{subarray}{l}
\text{$m\!\in\!M$, $n\in N_{nm}$, $p\in P_m$,}\\
\text{$\dot P\!\in\!\cP(Y_{mp},\jmath_m((\K^\t)^I))$,}\\
\text{$\dot Q\!\in\!\cQ(\jmath_m((\K^\t)^I)\!\lt\!K_{mp},
\jmath_m((\K^\t)^I))$,}\\
\text{$\dot R\!\in\!\cR(Y_{mp},\jmath_m((\K^\t)^I)\!\lt\!K_{mp},
\jmath_m((\K^\t)^I))$:}\\
\text{$\dot R\!\subseteq\!\dot P\!\cap\!\dot Q$,
$M^{Y_{mp}}_{\jmath_m((\K^\t)^I)\!\lt\!K_{mp}}(\dot P,\dot Q,\dot
R)\ne 0$, $\dim\dot R\!=\!k$}
\end{subarray}}
\!\!\!\!\!\!\!\!\!\!\!\!\!\!\!\!\!\!\!\!\!\!\!\!\!\!
\!\!\!\!\!\!\!\!\!
\begin{aligned}[t]
&M^{Y_{mp}}_{\jmath_m((\K^\t)^I)\!\lt\!K_{mp}}
(\dot P,\dot Q,\dot R)\cdot c_{mn}\cdot\\
&\bigl[\bigl([W_{mn}\!\t_{V_m}\!Y_{mp}^{\dot P}/
G_m\!\lt\!(K_{mp})^{\dot Q}],\!\!\!\!\!\! \\
&\qquad\qquad\qquad \bs\si(I)\!\ci\!\xi_{mnp}\bigr)\bigr].
\end{aligned}
\label{ai5eq24}
\e

Suppose that $m\in M$, $p\in P_m$, $\dot P\subseteq\jmath_m
((\K^\t)^I)$ is a $\K$-subgroup, and $y\in Y_{mp}(\K)$ is fixed by
$\dot P$. Then $y$ projects to $[(\si,\io,\pi)]\in\M(I,\tl)_\A$ with
a commutative diagram of $\K$-groups
\e
\begin{gathered}
\xymatrix@C=50pt@R=10pt{ *+[r]{\dot P} \ar[rrr] \ar[d]^{\,\subseteq}
&&& *+[l]{\Aut(\si,\io,\pi)} \ar[d]_{\prod_{i\in
I}\bs\si(\{i\})_*\,}
\\
*+[r]{\jmath_m((\K^\t)^I)} \ar[r]^(0.6){\jmath_m^{-1}}
& (\K^\t)^I \ar[rr]^(0.31){\de\mapsto\prod_{i\in I}\de(i)\id_{X_i}}
&&
*+[l]{\prod_{i\in I}\Aut(\si(\{i\})).} }
\end{gathered}
\label{ai5eq25}
\e
Now the bottom right corner extends naturally to algebra morphisms
between the algebras $\K^I,\prod_{i\in I}\End(\si(\{i\}))$ and
$\End(\si,\io,\pi)$. Let $A\subseteq\K^I$ be the subalgebra of
$\K^I$ generated by $\K$-subgroup $\jmath_m^{-1}(\dot P)$, and
$A^\t\subseteq(\K^\t)^I$ the $\K$-subgroup of invertible elements in
$A$. Then $\jmath_m(A^\t)$ is a $\K$-subgroup of
$\jmath_m((\K^\t)^I)$ with $\dot
P\subseteq\jmath_m(A^\t)\subseteq\jmath_m((\K^\t)^I)$, and it is
easy to see that we can extend \eq{ai5eq25} to replace $\dot P$ by
$\jmath_m(A^\t)$. Hence $y\in Y_{mp}(\K)$ is fixed
by~$\jmath_m(A^\t)$.

From \cite[Def.~5.3]{Joyc2} we now see that each $\dot
P\!\in\!\cP(Y_{mp},\jmath_m((\K^\t)^I))$ is of the form
$\jmath_m(A^\t)$ for some subalgebra $A\subseteq\K^I$. A related
proof shows the same holds for each $\dot Q\!\in\!\cQ(\jmath_m
((\K^\t)^I)\!\lt\!K_{mp},\jmath_m((\K^\t)^I))$. Let $K$ be a finite
set and $\phi:I\ra K$ a surjective map, and define
\begin{equation*}
A_{I,K,\phi}=\bigl\{\de\in\K^I:\text{$\phi(i)=\phi(j)$ implies
$\de(i)=\de(j)$, $i,j\in I$}\bigr\}.
\end{equation*}
Then $A_{I,K,\phi}$ is a $\K$-subalgebra of $\K^I$, and every
subalgebra is of this form. Thus
\begin{gather*}
\cP\bigl(Y_{mp},\jmath_m((\K^\t)^I)\bigr),\cQ\bigl(\jmath_m
((\K^\t)^I)\!\lt\!K_{mp},\jmath_m((\K^\t)^I)\bigr)\\
\subseteq\bigl\{A_{I,K,\phi}^\t:\text{$K$ finite, $\phi:I\ra K$
surjective}\bigr\}.
\end{gather*}

It is a consequence of the proof in \cite[\S 5]{Joyc2} that the
definition of $\Pi^\vi_k$ is independent of choices, that if in
\eq{ai5eq24} we replace $\cP\bigl(Y_{mp},\jmath_m((\K^\t)^I)\bigr),
\cQ\bigl(\jmath_m((\K^\t)^I)\!\lt\!K_{mp},\jmath_m((\K^\t)^I)\bigr)$,
$\cR\bigl(Y_{mp},\jmath_m((\K^\t)^I)\!\lt\!K_{mp},\jmath_m((\K^\t)^I)
\bigr)$ by larger finite sets closed under intersection, and compute
$M^{Y_{mp}}_{\jmath_m((\K^\t)^I)\!\lt\!K_{mp}}(\dot P,\dot Q,\dot
R)$ using these larger sets, then we get the same answer for
$\Pi^\vi_k[\cdots]$. Replacing all three sets by
$\{A_{I,K,\phi}^\t:K$ finite, $\phi:I\ra K$ surjective$\}$, equation
\eq{ai5eq24} becomes:
\e
\Pi^\vi_k\bigl[P_\sIt(f_i:i\!\in\! I)\bigr]\!=
\!\!\!\!\!\!\!\!\!\!\!\!\!\!\!\!\!\!\!\!
\sum_{\begin{subarray}{l}
\text{$m\!\in\!M$, $n\in N_{nm}$, $p\in P_m$,}\\
\text{$\dot P,\dot R\in\{A_{I,K,\phi}^\t:K$ finite,}\\
\text{$\phi:I\ra K$ surjective$\}$}\\ \text{$\dot R\!\subseteq\!\dot
P$, $M_I(\dot P,\dot P,\dot R)\ne 0$, $\dim\dot R\!=\!k$}
\end{subarray}}
\!\!\!\!\!\!\!\!\!\!\!\!\!\!\!\!
\begin{aligned}[t]
&M_I(\dot P,\dot P,\dot R)\cdot c_{mn}\cdot\\
&\bigl[\bigl([W_{mn}\!\t_{V_m}\!Y_{mp}^{\dot P}/
G_m\!\lt\!(K_{mp})^{\dot P}],\!\!\!\!\!\! \\
&\qquad\qquad\qquad \bs\si(I)\!\ci\!\xi_{mnp}\bigr)\bigr].
\end{aligned}
\label{ai5eq26}
\e

Here $M_I(\dot P,\dot P,\dot R)$ is $M^{Y_{mp}}_{\jmath_m((\K^\t)^I)
\!\lt\!K_{mp}}(\dot P,\dot P,\dot R)$ computed using $\{A_{I,K,
\phi}^\t:K$ finite, $\phi:I\ra K$ surjective$\}$ in place of
$\cP,\cQ,\cR(\cdots)$, and may be written explicitly as in
\eq{ai5eq37} below. We have also used the fact
\cite[Lem.~5.9]{Joyc2} that $M_I(\dot P,\dot Q,\dot R)=0$ unless
$\dot P,\dot Q$ are the smallest elements containing $\dot P\cap\dot
Q$, which forces $\dot P=\dot Q$ as $\dot P,\dot Q$ take values in
the same set. We may rewrite \eq{ai5eq26} as
\begin{gather}
\begin{gathered}
\Pi^\vi_k\bigl[P_\sIt(f_i:i\in I)\bigr]=
\sum_{\substack{\text{iso. classes of finite}\\
\text{sets $K$, $k\le\md{K}\le\md{I}$}}}\,\,
\sum_{\substack{\text{$\phi:I\ra K$}\\
\text{surjective}}}\!\!\! N_{I,K,\phi,k}\,\cdot
\\
\raisebox{-5pt}{\begin{Large}$\displaystyle\Bigl[$\end{Large}}
\sum_{m\in M,\; n\in N_{nm},\; p\in P_m,
\!\!\!\!\!\!\!\!\!\!\!\!\!\!\!\!\!\!\!\!\!\!\!\!\!\!\!\!\!\!\!\!
\!\!\!\!\!\!\!\!\!\!\!\! }
c_{mn}\bigl[\bigl([W_{mn}\!\t_{V_m}\!Y_{mp}^{A_{I,K,\phi}^\t}/
G_m\!\lt\!(K_{mp})^{A_{I,K,\phi}^\t}],\bs\si(I)\!\ci\!
\xi_{mnp}\bigr)\bigr]
\raisebox{-5pt}{\begin{Large}$\displaystyle\Bigr]$\end{Large}},
\end{gathered}
\label{ai5eq27}
\\
\begin{gathered}
\text{where}\quad N_{I,K,\phi,k}=\frac{1}{\md{K}!}\cdot
\sum_{\begin{subarray}{l}
\text{$\dot R\in\{A_{I,L,\psi}^\t:L$ finite,
$\psi:I\ra L$ surjective$\}:$}\\
\text{$\dot R\!\subseteq\!A_{I,K,\phi}^\t$,
$M_I(A_{I,K,\phi}^\t,A_{I,K,\phi}^\t,\dot R)\ne 0$, $\dim\dot
R\!=\!k$}
\end{subarray}
\!\!\!\!\!\!\!\!\!\!\!\!\!\!\!\!\!\!\!\!\!\!\!\!\!\!\!\!\!\!\!\!
\!\!\!\!\!\!\!\!\!\!\!\!\!\!\!\!\!\!\!\!\!\!\!\!\!\!\!\!\!\!\!\! }
\!\!\!\!\!\!\!\!\! M_I(A_{I,K,\phi}^\t,A_{I,K,\phi}^\t,\dot R),
\end{gathered}
\label{ai5eq28}
\end{gather}
setting $\dot P=A_{I,K,\phi}^\t$ and replacing the sum over $\dot P$
in \eq{ai5eq26} with sums over $K,\phi$. The factor $1/\md{K}!$ in
\eq{ai5eq28} compensates for the fact that if $\imath:K\ra K$ is a
bijection then $A_{I,K,\imath\ci\phi}^\t=A_{I,K,\phi}^\t$, and there
are $\md{K}!$ such bijections $\imath$, so each $A_{I,K,\phi}^\t$ is
represented by $\md{K}!$ choices of $K,\phi$ in \eq{ai5eq27}. Note
that $N_{I,K,\phi,k}$ in \eq{ai5eq28} lies in $\Q$ and depends only
on $I,K,\phi$ up to isomorphism and $k$, as we want.

Now let $I,K,\phi$ be as in \eq{ai5eq27} and define a partial order
$\pr$ on $I$ as in \eq{ai5eq13}. Then $\tl$ dominates $\pr$ and we
have a 1-morphism $Q(I,\pr,\tl):\fM(I,\pr)_\A\ra\fM(I,\tl)_\A$. In a
similar way to \eq{ai5eq15}, we claim there is a natural
1-isomorphism
\e
\fH_{mp}\t_{\io_{mp},\fM(I,\tl)_\A,Q(I,\pr,\tl)}\fM(I,\pr)_\A\cong
\bigl[Y_{mp}^{A_{I,K,\phi}^\t}/G_m\!\lt\!(K_{mp})^{A_{I,K,\phi}^\t}
\bigr],
\label{ai5eq29}
\e
where $\io_{mp}:\fH_{mp}\ra\fM(I,\tl)_\A$ is the inclusion. To see
this, observe that a point of the r.h.s.\ of \eq{ai5eq29} is
equivalent to a point $[(\si,\io,\pi)]$ in
$\fH_{mp}(\K)\subseteq\M(I,\tl)_\A$ together with a choice of
commutative diagram \eq{ai5eq25} in which~$\dot
P=\jmath_m^{-1}(A_{I,K,\phi}^\t)$.

That is, the $(I,\tl)$-configuration $(\si,\io,\pi)$ is equipped
with a choice of $\K$-group morphism $\rho:A_{I,K,\phi}^\t\ra
\Aut(\si,\io,\pi)$ such that $\prod_{i\in I}\bs\si(\{i\})_*\ci
\rho:A_{I,K,\phi}^\t\ra\prod_{i\in I}\Aut(\si(\{i\}))$ maps
$\de\mapsto(\de(i)\id_{\si(\{i\})})_{i\in I}$. It is not difficult
to show that choosing such a $\rho$ is equivalent, up to canonical
isomorphism, to choosing an $(I,\pr)$-{\it improvement\/} of
$(\si,\io,\pi)$, in the sense of \cite[\S 6]{Joyc3}, and
\eq{ai5eq29} follows. Using \eq{ai5eq16}, \eq{ai5eq17} and
\eq{ai5eq29} we see that the last line $[\cdots]$ of \eq{ai5eq27} is
\begin{equation*}
\ts\bs\si(I)_*\ci Q(I,\pr,\tl)^*\ci\bigl(\prod_{i\in
I}\bs\si(\{i\})\bigr)^*\bigl(\bigot_{i\in I}f_i\bigr)=
P_\sIp(f_i:i\in I),
\end{equation*}
and so \eq{ai5eq27} implies \eq{ai5eq13}, completing the proof of
Theorem~\ref{ai5thm3}.
\end{proof}

Write $[f,g]=f*g-g*f$ for $f,g\in\SFa(\fObj_\A)$ or
$\oSFa(\fObj_\A,*,*)$. Then $[\,,\,]$ satisfies the {\it Jacobi
identity} and is a {\it Lie bracket} by Theorem \ref{ai5thm1}. We
shall prove an analogue of Theorem~\ref{ai4thm4}:

\begin{thm} Let Assumptions \ref{ai2ass} and \ref{ai3ass} hold.
Then $\SFai(\fObj_\A)$ and\/ $\oSFai(\fObj_\A,*,*)$ are closed under
$[\,,\,],$ and are Lie algebras, and\/ \eq{ai5eq2} restricts to Lie
algebra morphisms $\SFai(\fObj_\A),\oSFai(\fObj_\A,*,*)\!\ra\!
\CFi(\fObj_\A)$ when $\K$ has characteristic zero.
\label{ai5thm4}
\end{thm}

\begin{proof} If $f,g\in\SFai(\fObj_\A)$ or $\oSFai
(\fObj_\A,*,*)$ then as $*=P_{\sst(\{1,2\},\le)}$ we have
\begin{align*}
\Pi^\vi_1&\bigl([f,g])=\Pi^\vi_1\bigl(P_{\sst(\{1,2\},\le)}(f,g)
\bigr)-\Pi^\vi_1\bigl(P_{\sst(\{1,2\},\le)}(g,f)\bigr)\\
&=\bigl(P_{\sst(\{1,2\},\le)}(f,g)\!-\!P_{\sst(\{1,2\},\bu)}(f,g)
\bigr)\!-\!\bigl(P_{\sst(\{1,2\},\le)}(g,f)
\!-\!P_{\sst(\{1,2\},\bu)}(g,f)\bigr)\\
&=P_{\sst(\{1,2\},\le)}(f,g)-P_{\sst(\{1,2\},\le)}(g,f)=[f,g]
\end{align*}
by Theorem \ref{ai5thm3}, where $\bu$ is the partial order on
$\{1,2\}$ with $i\bu j$ if $i=j$, so that by symmetry
$P_{\sst(\{1,2\},\bu)}(f,g)= P_{\sst(\{1,2\},\bu)}(g,f)$. Therefore
$[f,g]$ also lies in $\SFai(\fObj_\A)$ or $\oSFai(\fObj_\A,*,*)$,
proving the first part.

For the second part, combining Corollary \ref{ai5cor2} with the fact
that $\Pi^\vi_k$ is the identity on $[(U\t[\Spec\K/(\K^\t)^l],
\rho)]$ if $k=l$ and 0 otherwise, we see that $\oSFai(\fObj_\A,*,*)$
is generated over $\La,\La^\ci$ or $\Om$ by elements
$[(U\t[\Spec\K/\K^\t],\rho)]$ for $U$ a quasiprojective
$\K$-variety. If $u\in U(\K)$ with $\rho_*(u)=[X]\in\fObj_\A(\K)$
then $\rho_*:\K^\t\ra\Aut(X)$ is injective and $u$ contributes
$\chi\bigl(\Aut(X)/\rho_*(\K^\t)\bigr)$ to $\bar\pi_{\fObj_\A}^\stk
\bigl([(U\t[\Spec\K/\K^\t],\rho)]\bigr)$ at $[X]$. Since
$\K^\t\cong\rho(\K^\t)\subseteq\Aut(X)$ we see that $\rk\Aut(X)\ge
1$, and if $\rk\Aut(X)>1$ then the action of a maximal torus of
$\Aut(X)$ fibres $\Aut(X)/\rho_*(\K^\t)$ by tori, forcing
$\chi\bigl(\Aut(X)/\rho_*(\K^\t)\bigr)=0$.

Thus $u$ makes a nonzero contribution at $[X]$ only if
$\rk\Aut(X)=1$, that is, if $X$ is indecomposable. Hence
$\bar\pi_{\fObj_\A}^\stk\bigl([(U\t[\Spec\K/\K^\t],\rho)]\bigr)$ is
supported on indecomposables, and as such $[(U\t[\Spec\K/\K^\t],
\rho)]$ generate $\oSFai(\fObj_\A,*,*)$ we see that
$\bar\pi_{\fObj_\A}^\stk$ maps to $\CFi(\fObj_\A)$, and is a Lie
algebra morphism by Theorem \ref{ai5thm1}. The result for
$\pi_{\fObj_\A}^\stk$ follows as~$\pi_{\fObj_\A}^\stk\!
=\bar\pi_{\fObj_\A}^\stk\ci\bar\Pi^{*,*}_{\fObj_\A}$.
\end{proof}

\subsection{Relations between $*$ and $\Pi_{[I,\ka]}$ in
$\oSFa(\fObj_\A,*,*)$}
\label{ai53}

Theorem \ref{ai5thm4} gives a compatibility between multiplication
$*$ and the projections $\Pi^\vi_n,\Pi_{[I,\ka]}$, in that subspaces
$\SFai(\fObj_\A),\oSFai(\fObj_\A,*,*)$ defined using the $\Pi^\vi_n$
are closed under $(f,g)\mapsto f*g-g*f$. This is one consequence of
a deeper relationship, in which we can write $\Pi_{[K,\mu]}(f*g)$
explicitly in terms of the components $\Pi_{[I,\ka]}(f)$,
$\Pi_{[J,\la]}(g)$. This deeper relationship is very complicated to
write down, so for simplicity we do so only for
$f,g\in\oSFa(\fObj_\A,*,*)$, when we can modify Proposition
\ref{ai2prop2} to represent $f,g$ and $f\ot g$ in a special way.

\begin{thm} Suppose Assumptions \ref{ai2ass} and \ref{ai3ass} hold,
$[I,\ka],[J,\la],[K,\mu]$ are as in Definition \ref{ai5def4}, and\/
$f,g$ lie in $\oSFa(\fObj_\A,\Up,\La),\oSFa(\fObj_\A,\Up,\La^\ci)$
or $\oSFa(\fObj_\A,\Th,\Om)$ with\/ $\Pi_{[I,\ka]}(f)=f$ and\/
$\Pi_{[J,\la]}(g)=g$. Choose constructible
$T\subseteq\fObj_\A(\K)\t\fObj_\A(\K)$ with\/ $f\ot g$ supported on
$T$, and use the notation of Proposition \ref{ai5prop5}. Write
$(\K^\t)^I,\ldots$ for the $\K$-groups of
functions~$I,\ldots\ra\K^\t$.

Then we may represent $f\ot g\in\oSF(\fObj_\A\t\fObj_\A,*,*)$ as
\e
f\ot g=\ts\sum_{m\in M,\; n\in N_m}c_{mn}\bigl[\bigl(W_{mn}\t
[\Spec\K/(\K^\t)^I\t(\K^\t)^J],\tau_{mn}\bigr)\bigr],
\label{ai5eq30}
\e
where $N_m$ is finite, $c_{mn}\in\La,\La^\ci$ or $\Om$, $W_{mn}$ is
a quasiprojective $\K$-variety, and\/ $\tau_{mn}:W_{mn}\t[\Spec\K/
(\K^\t)^I\t(\K^\t)^J]\ra\fG_m\subseteq\fObj_\A\t\fObj_\A$ is
induced, using the $1$-isomorphism $\fG_m\cong[V_m/G_m]$, by an
injective $\K$-group morphism $\rho_{mn}:(\K^\t)^I\t(\K^\t)^J\ra
G_m$ and a morphism~$\si_{mn}:W_{mn}\ra V_m^{\rho_{mn}
((\K^\t)^I\t(\K^\t)^J)}\subseteq V_m$.

These have the property that if\/ $w\in W_{mn}(\K)$ projects to
$v\in V_m(\K)$ and\/ $([X],[Y])\in\fObj_\A(\K)\t\fObj_\A(\K)$ then
$\rho_{mn}$ maps to $\Stab_{G_m}(v)\subseteq G_m$, and there exist
splittings $X\cong\bigop_{i\in I}X_i$ and\/ $Y\cong\bigop_{j\in
J}Y_j$ in $\A$ with\/ $[X_i]=\ka(i)$ and\/ $[Y_j]=\la(j)$ in $C(\A)$
for all\/ $i,j,$ such that composing $\rho_{mn}$ with\/ \eq{ai5eq6}
yields the $\K$-group morphism $(\K^\t)^I\t(\K^\t)^J\ra\Aut(X)\t
\Aut(Y)$ given by~$(\ga,\de)\mapsto\bigl(\sum_{i\in I}\ga(i)
\id_{X_i},\sum_{j\in J}\de(j)\id_{Y_j}\bigr)$.

In this representation we have
\e
f*g=\sum_{\!\!\! m\in M,\; n\in N_m \!\!\!\!\!\!\!\!\!\!\!\!\!}
c_{mn}\bigl[\bigl(W_{mn}\t[E^1_m/((\K^\t)^I\t(\K^\t)^J)\lt E^0_m],
\bs\si(\{1,2\})\ci\xi_{mn}\bigr)\bigr],
\label{ai5eq31}
\e
where $(\K^\t)^I\t(\K^\t)^J$ acts on $E^0_m,E^1_m$ via $\rho_{mn}$
and the $G_m$-actions, and\/ $E^0_m$ acts trivially on $E^1_m$,
and\/ $\si_{mn},\rho_{mn}$ induce
\e
\xi_{mn}:W_{mn}\t[E^1_m/((\K^\t)^I\t(\K^\t)^J)\lt E^0_m]\ra
\bigl[V_m\t E^1_m/G_m\lt E^0_m\bigr],
\label{ai5eq32}
\e
using \eq{ai5eq7} to regard the r.h.s.\ of\/ \eq{ai5eq32} as a
substack of\/~$\fM(\{1,2\},\le)_\A$.

If\/ $L$ is a finite set and $\phi:I\ra L$, $\psi:J\ra L$ maps with
$\phi\amalg\psi:I\amalg J\ra L$ surjective, define
$T_{L,\phi,\psi}\subseteq(\K^\t)^I\t(\K^\t)^J$ to be the subtorus
of\/ $(\ga,\de)\in(\K^\t)^I\t(\K^\t)^J$ for which there exists
$\ep:L\ra\K^\t$ with\/ $\ga(i)=\ep\ci\phi(i)$ and\/
$\de(j)=(\ep\ci\psi(j))^{-1}$ for all\/ $i\in I$ and\/ $j\in J$.
Then $\ep$ determines $\ga,\de$ uniquely, so that\/
$T_{L,\phi,\psi}\cong(\K^\t)^L$. For each\/ $m\in M$, $n\in N_m$
and\/ $i=0,1$ write $(E^i_m)^{T_{L,\phi,\psi}}$ for the vector
subspace of\/ $E^i_m$ fixed by $\rho_{mn}(T_{L,\phi,\psi})$. Write
$\Aut(K,\mu)$ for the finite group of bijections $\io:K\ra K$ with\/
$\mu=\mu\ci\io$. Then
\e
\begin{split}
&\Pi_{[K,\mu]}(f*g)=
\\
&\frac{1}{\md{\Aut(K,\mu)}}\!
\sum_{\substack{\text{iso.}\\ \text{classes}\\ \text{of finite}\\
\text{sets $L$}}}\!\!\frac{(-1)^{\md{L}-\md{K}}}{\md{L}!}\!\!\!\!\!
\sum_{\substack{\text{$\phi:I\ra L$, $\psi:J\ra L$ and}\\
\text{$\th:L\!\ra\!K$: $\phi\!\amalg\!\psi$ surjective,}\\
\text{$\mu(k)=\ka((\th\ci\phi)^{-1}(k))+$}\\
\text{$\la((\th\ci\psi)^{-1}(k))$, $k\in K$}}}\,\prod_{k\in
K}\!(\md{\th^{-1}(k)}\!-\!1)!
\\
&\sum_{\!m\in M,\; n\in N_m \!\!\!\!\!\!\!}
\begin{aligned}[t]
c_{mn}\bigl[\bigl(W_{mn}\!\t\![(E^1_m)^{T_{L,\phi,\psi}}
/((\K^\t)^I\!\t\!(\K^\t)^J)\!\lt\!(E^0_m)^{T_{L,\phi,\psi}}],&\\
\bs\si(\{1,2\})\ci\xi_{mn}\bigr)\bigr].&
\end{aligned}
\end{split}
\label{ai5eq33}
\e
\label{ai5thm5}
\end{thm}

\begin{proof} By Corollary \ref{ai5cor3} we can write $f\ot g$ in
the form \eq{ai5eq9}. The proofs of the first part of Proposition
\ref{ai2prop2} in \cite[\S 5.3]{Joyc2} and Corollary \ref{ai5cor2}
then show we can write $f\ot g$ as a sum of terms
$c_{mn}\bigl[\bigl(W_{mn}\t[\Spec\K/T],\tau_{mn}\bigr)\bigr]$, where
$T\cong(\K^\t)^k$ and $\tau_{mn}$ maps to
$\fG_m\subseteq\fObj_\A\t\fObj_\A$ and is induced by a $\K$-group
morphism $\rho_{mn}:T\ra G_m$ and a morphism~$\si_{mn}:W_{mn}\ra
V_m^{\smash{\rho_{mn}(T)}}\subseteq V_m$.

As $\Pi_{[I,\ka]}(f)=f$ and $\Pi_{[J,\la]}(g)=g$ we can show using
Definition \ref{ai5def4} that the sum can be chosen such that when
$w\in W_{mn}(\K)$ projects to $v\in V_m(\K)$ and
$([X],[Y])\in\fObj_\A(\K)\t\fObj_\A(\K)$, there is an isomorphism
$T\cong(\K^\t)^I\t(\K^\t)^J$ for which the second and third
paragraphs hold. This isomorphism may depend on $w$, but it does so
constructibly, so refining the sum we can take the isomorphism to be
constant on $W_{mn}(\K)$, and identify $T$ with
$(\K^\t)^I\t(\K^\t)^J$. This gives the first three paragraphs of the
theorem.

Equation \eq{ai5eq31} now follows from \eq{ai5eq30} as for
\eq{ai5eq12}. To prove \eq{ai5eq33} we apply $\Pi_{[K,\mu]}$ to
\eq{ai5eq31} and use Definition \ref{ai5def4} and \cite[\S
5.2]{Joyc2}. Deleting terms with $W_{mn}(\K)=\emptyset$, let $m\in
M$ and $n\in N_m$ and pick $w\in W_{mn}(\K)$, projecting to $v\in
V_m(\K)$ and $([X],[Y])\in\fObj_\A(\K)\t\fObj_\A(\K)$. Then the
first part of the theorem gives splittings $X\cong\bigop_{i\in
I}X_i$ and $Y\cong\bigop_{j\in J}Y_j$. By Proposition
\ref{ai5prop5}(a) we have isomorphisms
\begin{align*}
&E^0_m\cong\Hom(Y,X)\cong\ts\bigop_{i\in I,\; j\in J}\Hom(Y_j,X_i),\\
&E^1_m\cong\Ext^1(Y,X)\cong\ts\bigop_{i\in I,\; j\in J}
\Ext^1(Y_j,X_i).
\end{align*}
Under these isomorphisms, $(\ga,\de)\in(\K^\t)^I\t(\K^\t)^J$ acts on
$E^0_m,E^1_m$ via $\rho_{mn}$ by multiplying by $\ga(i)\de(j)$ in
$\Hom(Y_j,X_i),\Ext^1(Y_j,X_i)$. It follows that for $L,\phi,\psi$
and $T_{L,\phi,\psi}$ as in the theorem we have
\e
(E^0_m)^{T_{L,\phi,\psi}}\!\cong\!\!\!\bigop_{\substack{i\in I,\;
j\in J:\\ \phi(i)=\psi(j)}}\!\!\!\Hom(Y_j,X_i),\;\>
(E^1_m)^{T_{L,\phi,\psi}}\!\cong\!\!\!\bigop_{\substack{i\in I,\;
j\in J:\\ \phi(i)=\psi(j)}}\!\!\!\Ext^1(Y_j,X_i).
\label{ai5eq34}
\e

In applying $\Pi_{[K,\mu]}$ to \eq{ai5eq31} we can take the $W_{mn}$
factor outside, as $((\K^\t)^I\t(\K^\t)^J)\lt E^0_m$ acts trivially
upon it. If $\Hom(Y_j,X_i)\not\cong 0\not\cong\Ext^1(Y_j,X_i)$ for
all $i,j$ then using \eq{ai5eq34} and the notation of \cite[\S
5.2]{Joyc2} we find that
\e
\begin{split}
\cP\bigl(&E^1_m,(\K^\t)^I\!\t\!(\K^\t)^J\bigr)\!=\!\cQ\bigl(
((\K^\t)^I\!\t\!(\K^\t)^J)\!\lt\!E^0_m,(\K^\t)^I\!\t\!(\K^\t)^J
\bigr)
\\
&=\cR\bigl(E^1_m,((\K^\t)^I\t(\K^\t)^J)\lt
E^0_m,(\K^\t)^I\t(\K^\t)^J\bigr)
\\
&=\bigl\{T_{L,\phi,\psi}:\text{$L$ a finite set, $\phi:I\!\ra\!L$,
$\psi:J\!\ra\!L$, $\phi\!\amalg\!\psi$ surjective}\bigr\}.
\end{split}
\label{ai5eq35}
\e

If some $\Hom(Y_j,X_i),\Ext^1(Y_j,X_i)$ are zero then
$\cP,\cQ,\cR(\cdots)$ may be subsets of the values in \eq{ai5eq35}.
However, since the formulae in \cite[\S 5.2]{Joyc2} give the same
answers if we replace $\cP(X,T^G),\cQ(G,T^G),\cR(X,G,T^G)$ by larger
finite sets of $\K$-subgroups of $T^G$ closed under intersections,
the values \eq{ai5eq35} give the correct answer for computing
$\Pi_{[K,\mu]}(\cdots)$, so we shall still use them.

Now let $P\in\cP(\cdots)$ and $R\in\cR(\cdots)$ with $R\subseteq P$
and $0\ne c\in\Q$, and $\nu$ be the $\fObj_\A$-weight function of
Definition \ref{ai5def4} defining $\Pi_{[K,\mu]}$. Then
\cite[Def.~5.15]{Joyc2} defines a constructible set $(W_{mn}\t
E^1_m)^{P,R}_{\nu,c}$ in $(W_{mn}\t E^1_m)(\K)$, which we will
evaluate. Let $(w,e)\in(W_{mn}\t E^1_m)(\K)$ project to
$([X],[Y])\in\fObj_\A(\K)\t\fObj_\A(\K)$ under
$(\bs\si(\{1\})\t\bs\si(\{2\}))_*$ and to $[Z]\in\fObj_\A(\K)$ under
$\bs\si(\{1,2\})_*$. Let $R=T_{L,\phi,\psi}$ for some $L,\phi,\psi$.
Then there is an exact sequence $0\ra X\ra Z\ra Y\ra 0$ invariant
under $R$. The splittings $X\cong\bigop_{i\in I}X_i$,
$Y\cong\bigop_{j\in J}Y_j$ correspond to a splitting
$Z\cong\bigop_{l\in L}Z_l$ with $Z_l$ in an exact sequence
\begin{equation*}
0\longra\ts\bigop_{i\in I:\;\phi(i)=l}X_i\longra Z_l\longra
\ts\bigop_{j\in J:\;\psi(j)=l}Y_j\longra 0.
\end{equation*}
It follows that $[Z_l]=\ka(\phi^{-1}(l))+\la(\psi^{-1}(l))$
in~$C(\A)$.

Under the natural isomorphism $R\cong(\K^\t)^L$, $\ep\in(\K^\t)^L$
acts on $Z$ as $\sum_{l\in L}\ep(l)\id_{Z_l}$. From Definition
\ref{ai5def4} we find $\nu\bigl(R,[Z],(\bs\si(\{1,2\})\ci
\xi_{mn})_*\bigr)$ is 1 if there is a bijection $\th:L\ra K$ with
$\mu(k)=\ka(\bar\phi^{-1}(k))+\la(\bar\psi^{-1}(k))$ for $k\in K$,
where $\bar\phi=\th\ci\phi$ and $\bar\psi=\th\ci\psi$, and 0
otherwise. Then $R=T_{K,\bar\phi,\bar\psi}$. So
\cite[Def.~5.15]{Joyc2} gives $(W_{mn}\t
E^1_m)^{P,R}_{\nu,c}=W_{mn}(\K)\t(E^1_m)^P(\K)$ if $c=1$ and
$R=T_{K,\bar\phi,\bar\psi}$ as above, and $(W_{mn}\t
E^1_m)^{P,R}_{\nu,c}=\emptyset$ otherwise for $c\ne 0$. The
definitions now yield:
\begin{gather}
\Pi_{[K,\mu]}\bigl(\bigl[\bigl(W_{mn}\t[E^1_m/((\K^\t)^I\t(\K^\t)^J)
\lt E^0_m],\bs\si(\{1,2\})\ci\xi_{mn}\bigr)\bigr]\bigr)=
\label{ai5eq36}
\\
\sum_{\substack{\text{$P,R\!\in\!
\cP(E^1_m,(\K^\t)^I\!\t\!(\K^\t)^J):R\!\subseteq\!P$,}\\
\text{$R=T_{K,\bar\phi,\bar\psi}$, $\bar\phi:I\ra
K$, $\bar\psi:J\ra K$,}\\
\text{$\mu(k)=\ka(\bar\phi^{-1}(k))+\la(\bar\psi^{-1}(k))$, $k\in
K$,}\\
\text{$M^{E^1_m}_{((\K^\t)^I\t(\K^\t)^J)\lt
E^0_m}(P,P,R)\!\ne\!0$,}}}
\begin{aligned}[t]
&M^{E^1_m}_{((\K^\t)^I\t(\K^\t)^J)\lt E^0_m}(P,P,R)\,\cdot\\
&\bigl[\bigl([W_{mn}\!\t\![(E^1_m)^P/((\K^\t)^I\!\t\!(\K^\t)^J)
\!\lt\!(E^0_m)^P],\\
&\qquad\qquad\qquad \bs\si(\{1,2\})\ci\xi_{mn}\bigr)\bigr].
\end{aligned}
\nonumber
\end{gather}

Here \cite[\S 5.2]{Joyc2} actually defines $\Pi^\nu_\fF(\cdots)$ as
a sum over triples $P\in\cP(X,T^G)$, $Q\in\cQ(G,T^G)$ and
$R\in\cR(X,G,T^G)$ with $R\subseteq P\cap Q$ and $M^X_G(P,Q,R)\ne
0$. But $P,Q$ are the smallest elements of $\cP(X,T^G),\cQ(G,T^G)$
containing $P\cap Q$ by \cite[Lem.~5.9]{Joyc2}, so \eq{ai5eq35}
implies that $P=Q$, which we have included in~\eq{ai5eq36}.

To get from \eq{ai5eq36} to \eq{ai5eq33} we set $P=T_{L,\phi,\psi}$
and $R=T_{K,\bar\phi,\bar\psi}$ in \eq{ai5eq33}. Then $R\subseteq P$
if and only if there exists a surjective $\th:L\ra K$ with
$\bar\phi=\th\ci\phi$ and $\bar\psi=\th\ci\psi$, which is then
unique. So we replace the sums over $P,R$ in \eq{ai5eq36} with sums
over isomorphism classes of $L$ and maps $\phi,\psi,\th$, as in
\eq{ai5eq33}. The combinatorial factors in the second line of
\eq{ai5eq33} are then the product of
\e
M^{E^1_m}_{((\K^\t)^I\t(\K^\t)^J)\lt E^0_m}(P,P,R)
=(-1)^{\md{L}-\md{K}}\ts\prod_{k\in K}(\md{\th^{-1}(k)}-1)!
\label{ai5eq37}
\e
in \eq{ai5eq36}, which can be computed explicitly using \eq{ai5eq35}
and \cite[Def.~5.8]{Joyc2}, and factors $1/\md{\Aut(K,\mu)}$,
$1/\md{L}!$ to compensate for the fact that each pair $(P,R)$ in
\eq{ai5eq36} is $(T_{L,\phi,\psi},T_{K,\th\ci\phi,\th\ci\psi})$ for
exactly $\md{\Aut(K,\mu)}\cdot\md{L}!$ quadruples
$(L,\phi,\psi,\th)$ in \eq{ai5eq33}. This completes the proof.
\end{proof}

This theorem will be very useful in \S\ref{ai6} when we impose extra
assumptions on $\A$ implying formulae for $\dim E^0_m-\dim E^1_m$
and $\dim(E^0_m)^{T_{L,\phi,\psi}}-\dim(E^1_m)^{T_{L,\phi,\psi}}$,
as then \eq{ai5eq33} will enable us to construct {\it algebra
morphisms} from $\oSFa(\fObj_\A,*,*)$ to certain explicit algebras
$B(\A,\La,\chi),B(\A,\La^\ci,\chi)$ and~$C(\A,\Om,\chi)$.

Theorems \ref{ai5thm3}--\ref{ai5thm5} are an important reason for
introducing {\it virtual rank\/} in \cite[\S 5]{Joyc2} and the
operators $\Pi^\vi_n,\Pi_{[I,\ka]}$. They show these operators have
a useful compatibility with $P_\sIp$ and $*$, although this is
difficult to state. For comparison, the simpler idea of {\it real
rank\/} and operators $\Pi^{\rm re}_n$ in \cite[\S 5.1]{Joyc2} have
no such compatibility with $*$, as far as the author knows.

\subsection{Generalization of other parts of \S\ref{ai4}}
\label{ai54}

So far we have generalized the material of \S\ref{ai41} and
\S\ref{ai48} to stack algebras. We can also generalize much of the
rest of \S\ref{ai4}, using the techniques of
\S\ref{ai51}--\S\ref{ai52}. This is straightforward, so we just
sketch the main ideas.

For \S\ref{ai42}, we can use the idea of {\it local stack functions}
$\uLSF,\LSF(\fF)$ in \cite[\S 4]{Joyc2} to define spaces
$\duLSF,\dLSF(\fObj_\A)$ and $\duLSF,\duoLSF,\doLSF(\fObj_\A,*,*)$
of stack functions $f$ in $\uLSF,\LSF(\fObj_\A)$ and
$\uLSF,\uoLSF,\oLSF(\fObj_\A,*,*)$ supported on $\coprod_{\al\in
S}\fObj_\A^\al$ for finite $S\subseteq\bar C(\A)$, and subspaces
$\dLSFa(\fObj_\A)$, $\doLSFa(\fObj_\A,*,*)$ analogous to
$\SFa(\fObj_\A)$, $\oSFa(\fObj_\A,*,*)$. Multiplication $*$ is
well-defined on all of these, making them into associative algebras.

For \S\ref{ai43}, define $*_L:\uSF(\fObj_\A)\!\t\!\uSF
(\fM(\{1,2\},\le)_\A)\!\ra\!\uSF(\fM(\{1,2\},\le)_\A)$~by
\begin{align*}
f*_Lr&=Q(\{1,2,3\},\le,\{1,2\},\le,\be)_*\\
&\qquad \bigl[\bigl(\bs\si(\{2\})\t Q(\{1,2,3\},\le,\{1,2\},\le,\al)
\bigr)^*(f\ot r)\bigr],
\end{align*}
by analogy with \eq{ai4eq6}, and define $*_R$ in a similar way. Then
the analogue of Theorem \ref{ai4thm3} holds and shows $*_L,*_R$ are
left and right representations of the algebra $\uSF(\fObj_\A)$ on
$\uSF\bigl(\fM(\{1,2\},\le)_\A)$. The same result holds if we work
with any of the spaces $\SF,\uSF,\uoSF,\oSF(*)$ or their local stack
function versions. If $\K$ has characteristic zero, the linear maps
$\pi^\stk_{\fObj_\A},\pi^\stk_{\fM(\{1,2\},\le)_\A}$ intertwine the
representations $*_L,*_R$ of $\SF(\fObj_\A)$ on $\SF(\fM
(\{1,2\},\le)_\A)$ with the representations $*_L,*_R$ of
$\CF(\fObj_\A)$ on $\CF(\fM(\{1,2\},\le)_\A)$ of~\S\ref{ai43}.

The generalization of \S\ref{ai45} is immediate, with Assumption
\ref{ai4ass} implying that the subspaces $\uSF_\fin,\SF_\fin
(\fObj_\A)$ and $\uSF_\fin,\uoSF_\fin,\oSF_\fin(\fObj_\A,*,*)$ of
stack functions with {\it finite support\/} are closed under $*$.
Each stack function in one of these subspaces is a sum of stack
functions supported over single points $[X]\in\fObj_\A(\K)$. So
using the relations in $\uSF_\fin,\uoSF_\fin,\oSF_\fin
(\fObj_\A,*,*)$ we can write down simple representations of these
subspaces. For instance:

\begin{lem} The subspace $\uSF_\fin(\fObj_\A,\Up,\La)$ of\/ $f$
in $\uSF(\fObj_\A,\Up,\La)$ with finite support has $\La$-basis
$[(\Spec\K,X)]$ for $[X]\in\fObj_\A(\K)$. Thus, $\io_{\fObj_\A}:
\CF_\fin(\fObj_\A)\ot_\Q\La\ra\uSF_\fin(\fObj_\A,\Up,\La)$ is an
isomorphism.
\label{ai5lem}
\end{lem}

If Assumption \ref{ai4ass} holds then $\CF_\fin(\fObj_\A)\ot_\Q\La$
and $\uSF_\fin(\fObj_\A,\Up,\La)$ are both $\La$-algebras, but in
general $\io_{\fObj_\A}$ is {\it not\/} an isomorphism between them.
Rather, $\uSF_\fin(\fObj_\A,\Up,\La)$ may be thought of as a
`quantized' version of the `classical' algebra
$\CF_\fin(\fObj_\A)\ot_\Q\La$.

To generalize \S\ref{ai46}, let Assumptions \ref{ai2ass},
\ref{ai3ass} and \ref{ai4ass} hold, and consider the $\La$-{\it
universal enveloping algebra} $U^{\sst\La}(\bar{\rm SF}{}^{\rm
ind}_{\rm al,fin}(\fObj_\A,\Up,\La))$ of the $\La$-Lie subalgebra
$\bar{\rm SF}{}^{\rm ind}_{\rm al,fin}(\fObj_\A,\Up,\La)$ of stack
functions with finite support in $\oSFai(\fObj_\A,\Up,\La)$. Then
$\Phi_\fin:U^{\sst\La}(\bar{\rm SF}{}^{\rm ind}_{\rm
al,fin}(\fObj_\A,\Up,\La))\ra\bar{\rm SF}{}_{\rm
al,fin}(\fObj_\A,\Up,\La)$ is an {\it isomorphism}. The same holds
for $U^{\sst\La^\ci}(\bar{\rm SF}{}^{\rm ind}_{\rm
al,fin}(\fObj_\A,\Up,\La^\ci))$ and $U^{\sst\Om}(\bar{\rm SF}{}^{\rm
ind}_{\rm al,fin}(\fObj_\A,\Th,\Om))$. The author is not sure
whether $\Phi$ is injective in the non-finite-support case.

We do not generalize \S\ref{ai47}, as there are problems with the
analogue of the proof that $\De$ is multiplicative in Theorem
\ref{ai4thm5}. We have already discussed the analogue of
\S\ref{ai48} in \S\ref{ai51}. In the quiver examples of
\S\ref{ai49}, we can use the algebra $\uSF_\fin(\fObj_\A, \Up,\La)$
to construct examples of {\it quantum groups}.

\begin{ex} Suppose Assumption \ref{ai2ass} holds for $\K,\Up,\La$,
and let $\g=\n_+\op\h\op\n_-$, $\Ga$ and $Q=(Q_0,Q_1,b,e)$ be as in
Example \ref{ai4ex1} or \ref{ai4ex2}. Then Assumption \ref{ai3ass}
holds for $\A=\modKQ$, and there are simple elements $V_i\in\A$ for
$i\in Q_0$. Write $\bde_{[V_i]}$ for the stack function in
$\uSF(\fObj_\A,\Up,\La)$ associated to the constructible set
$\{[V_i]\}\subseteq\fObj_\A(\K)$, and let $\bar\cC$ be the
subalgebra of $\uSF(\fObj_\A,\Up,\La)$ generated by the
$\bde_{[V_i]}$ for $i\in Q_0$. This is a stack function version of
the {\it composition algebra} $\cC$ of \S\ref{ai41}. Define $e_{ij}$
for $i,j\in Q_0$ by $e_{ii}=1$ and $-e_{ij}$ is the number of edges
$\smash{\mathop{\bu}\limits^{\sst i}\ra\mathop{\bu}\limits^{\sst
j}}$ in $Q$. Then $\dim\Ext^1(V_i,V_j)=-e_{ij}$ for $i\ne j$ in
$Q_0$. Define~$a_{ij}=e_{ij}+e_{ji}$.

Now the isomorphism $\cC\cong U(\n_+)$ in Example \ref{ai4ex2}
identifies $\de_{[V_i]}$ with $e_i\in U(\n_+)$, and holds because
the $\de_{[V_i]}$ satisfy the identity
\e
(\ad(\de_{[V_i]}))^{1-a_{ij}}\de_{[V_j]} \quad\text{in
$\CF(\fObj_\A)$ if $i\ne j\in Q_0$,}
\label{ai5eq38}
\e
known as the {\it Serre relations}. We will explain the
corresponding relations for the $\bde_{[V_i]}$ in
$\uSF(\fObj_\A,\Up,\La)$. For $0\le k\le n$ define the {\it Gauss
polynomial\/}
\begin{equation*}
\binom{n}{k}_\el=\frac{(\el^n-1)(\el^{n-1}-1)\cdots
(\el^{n-k+1}-1)}{(\el^k-1)(\el^{k-1}-1)\cdots(\el-1)}\,.
\end{equation*}
With this notation, for $i\ne j\in Q_0$ we claim that
\e
\sum_{k=0}^{1-a_{ij}}(-1)^k\binom{1-a_{ij}}{k}_\el
\el^{-k(1-a_{ij}-k)/2}\bde_{[V_i]}^{*^{\scriptstyle k}}*\bde_{[V_j]}*
\bde_{[V_i]}^{*^{\scriptstyle 1-a_{ij}-k}}=0.
\label{ai5eq39}
\e
Here $f^{*^k}$ means $f*f*\cdots*f$, with $f$ occurring $k$ times.
Equation \eq{ai5eq39} is known as the {\it quantum Serre relations},
and with $q$ in place of $\el$ and $X_i^+$ in place of
$\bde_{[V_i]}$ it was introduced by Drinfeld \cite[Ex.~6.2]{Drin} as
the defining relations of the {\it quantum group} $U_q(\n_+)$. We
recover \eq{ai5eq38} from \eq{ai5eq39} by replacing $\bde_{[V_i]}$
by $\de_{[V_i]}$ and taking the limit~$\el\ra 1$.

To prove \eq{ai5eq39} we can adapt the proof in Ringel \cite[\S
2]{Ring3}. Ringel works over finite fields with $q$ elements, so
that to define his Ringel--Hall multiplication he simply counts
numbers of filtrations satisfying some conditions, giving answers
which are polynomials in $q$. In our case, using the ideas of
\S\ref{ai51}--\S\ref{ai52} we translate Ringel's manipulation of
finite counts $q^k$ into addition and subtraction of constructible
sets of the form $\K^k$, which become factors $\el^k$ by the
relations in $\uSF(\fObj_\A,\Up,\La)$. We leave the details to the
reader.

Let the quantum group $U_\el(\n_+)$ be defined by the usual quantum
Serre relations over the algebra $\La$. Then from \eq{ai5eq39} we
obtain a unique, surjective algebra morphism
$U_\el(\n_+)\ra\bar\cC$. In the case of Example \ref{ai4ex1}, when
$\uSF_\fin(\fObj_\A,\Up,\La)= \uSF(\fObj_\A,\Up,\La)$, we can use
Lemma \ref{ai5lem} and $U(\n_+)\cong\cC$ to show this morphism is an
{\it isomorphism}. This is a Ringel--Hall-type realization of
$U_\el(\n_+)$ using stack functions, parallel to those of Ringel
\cite{Ring3} using finite fields, and Lusztig \cite[Th.~10.17]{Lusz}
using perverse sheaves.
\label{ai5ex}
\end{ex}

\section{Morphisms from stack (Lie) algebras}
\label{ai6}

If $\A=\modKQ$ for a quiver $Q$, or $\A=\coh(P)$ for a smooth
projective curve $P$, there is a biadditive $\chi:K(\A)\t
K(\A)\ra\Z$ called the {\it Euler form} with
\begin{equation*}
\dim_\K\Hom(X,Y)-\dim_\K\Ext^1(X,Y)=\chi\bigl([X],[Y]\bigr)
\quad\text{for all $X,Y\in\A$.}
\end{equation*}
Assuming this, \S\ref{ai61}--\S\ref{ai65} construct {\it algebra
morphisms} $\Phi^{\sst\La},\Psi^{\sst\La},\Psi^{\sst\La^\ci},
\Psi^{\sst\Om}$ from $\uSF(\fObj_\A),\oSFa(\fObj_\A,*,*)$ to {\it
explicit algebras} $A(\A,\La,\chi),B(\A,\La$ or $\La^\ci,\chi)$ and
$C(\A,\Om,\chi)$ depending only on $C(\A),\chi,\La,\La^\ci,\Om$,
which restrict to {\it Lie algebra morphisms} from
$\oSFai(\fObj_\A,*,*)$ to~$B^\ind(\A,\La$ or
$\La^\ci,\chi),C^\ind(\A,\Om,\chi)$.

In a similar way, if $\A=\coh(P)$ for $P$ a {\it Calabi--Yau
$3$-fold\/} then
\begin{align*}
&\bigl(\dim_\K\Hom(X,Y)-\dim_\K\Ext^1(X,Y)\bigr)-\\
&\bigl(\dim_\K\Hom(Y,X)-\dim_\K\Ext^1(Y,X)\bigr)=
\bar\chi\bigl([X],[Y]\bigr)\quad\text{for all $X,Y\in\A$,}
\end{align*}
for an antisymmetric bilinear $\bar\chi:K(\A)\t K(\A)\ra\Z$.
Assuming this, \S\ref{ai66} constructs a {\it Lie algebra
morphism}~$\Psi^{\sst\Om}:\oSFai(\fObj_\A,\Th,\Om)\ra
C^\ind(\A,\Om,\ha\bar\chi)$.

These (Lie) algebra morphisms will be essential tools in the sequels
\cite{Joyc4,Joyc5} on {\it stability conditions}. Given a
permissible weak stability condition $(\tau,T,\le)$ on $\A$, in
\cite{Joyc4} we will construct interesting subalgebras
$\bar{\mathcal H}^{\rm pa}(\tau),\bar{\mathcal H}^{\rm to}(\tau)$ of
$\SFa(\fObj_\A)$ and Lie subalgebras $\bar{\mathcal L}^{\rm
pa}(\tau),\bar{\mathcal L}^{\rm to}(\tau)$ of $\SFai(\fObj_\A)$.
Under mild conditions, \cite{Joyc5} shows these (Lie) subalgebras
are {\it independent of the choice of\/} $(\tau,T,\le)$, and gives
combinatorial {\it basis change formulae} relating bases of
$\bar{\mathcal H}^{\rm pa}(\tau),\ldots,\bar{\mathcal L}^{\rm
to}(\tau)$ associated to $(\tau,T,\le)$ and~$(\ti\tau,\ti T,\le)$.

Thus $\Phi^{\sst\La},\Psi^{\sst\La},\Psi^{\sst\La^\ci},
\Psi^{\sst\Om}$ induce morphisms $\bar{\mathcal H}^{\rm
pa}(\tau),\bar{\mathcal H}^{\rm to}(\tau)\!\ra\!A,B,C(\A,*,\chi)$
and $\bar{\mathcal L}^{\rm pa}(\tau), \bar{\mathcal L}^{\rm
to}(\tau)\ra B^\ind,C^\ind(\A,*,\chi)$. We shall regard these as
encoding interesting {\it families of invariants} of
$\A,(\tau,T,\le)$, which `count' $\tau$-(semi)stable objects and
configurations of objects in $\A$ with fixed classes in $C(\A)$. The
fact that $\Phi^{\sst\La},\ldots,\Psi^{\sst\Om}$ are (Lie) algebra
morphisms imply {\it multiplicative identities} on the invariants,
and the basis change formulae imply {\it transformation laws} for
the invariants from $(\tau,T,\le)$ to~$(\ti\tau,\ti T,\le)$.

\subsection{Identities relating $i_\La$ and $*,P_\sIp$}
\label{ai61}

Recall that a $\K$-linear abelian category $\A$ is called of {\it
finite type} if $\Ext^i(X,Y)$ is a finite-dimensional $\K$-vector
space for all $X,Y\in\A$ and $i\ge 0$, and $\Ext^i(X,Y)=0$ for $i\gg
0$. Then there is a unique biadditive map $\chi:K_0(\A)\t
K_0(\A)\ra\Z$ on the Grothendieck group $K_0(\A)$ known as the {\it
Euler form}, satisfying
\e
\chi\bigl([X],[Y]\bigr)=\ts\sum_{i\ge 0}(-1)^i\dim_\K\Ext^i(X,Y)
\quad\text{for all $X,Y\in\A$.}
\label{ai6eq1}
\e
We shall suppose $K(\A)$ in Assumption \ref{ai3ass} is chosen such
that $\chi$ factors through the projection $K_0(\A)\ra K(\A)$, and
so descends to $\chi:K(\A)\t K(\A)\ra\Z$. All this holds for the
examples $\A=\coh(P)$ in \cite[Ex.~9.1]{Joyc3} with $P$ a smooth
projective $\K$-scheme, for the quiver examples $\A=\modKQ$ in
\cite[Ex.~10.5]{Joyc3}, and for many of the other examples
of~\cite[\S 10]{Joyc3}.

Now assume $\Ext^i(X,Y)=0$ for all $X,Y\in\A$ and $i>1$. Then
\eq{ai6eq1} becomes
\e
\dim_\K\Hom(X,Y)-\dim_\K\Ext^1(X,Y)=\chi\bigl([X],[Y]\bigr)
\quad\text{for all $X,Y\in\A$.}
\label{ai6eq2}
\e
This happens for $\A=\coh(P)$ in \cite[Ex.~9.1]{Joyc3} with $P$ a
smooth projective curve, and for $\A=\modKQ$ in
\cite[Ex.~10.5]{Joyc3}. We shall prove that multiplication $*$ on
$\uSF(\fObj_\A,\Up,\La)$ and $i_\La$ in Proposition \ref{ai2prop1}
satisfy an important identity.

\begin{thm} Let Assumptions \ref{ai2ass} and \ref{ai3ass} hold,
and\/ $i_\La$ be as in Proposition \ref{ai2prop1}. Suppose
$\chi:K(\A)\t K(\A)\ra\Z$ is biadditive and satisfies \eq{ai6eq2}.
Let\/ $f,g\in\uSF(\fObj_\A,\Up,\La)$ be supported on
$\fObj_\A^\al,\fObj_\A^\be$ respectively, for $\al,\be\in\bar
C(\A)$. Write $\Pi:\fObj_\A\ra\Spec\K$ for the projection. Then
\e
i_\La^{-1}\ci\Pi_*(f*g)=\el^{-\chi(\be,\al)}\,\bigl(i_\La^{-1}
\ci\Pi_*(f)\bigr)\bigl(i_\La^{-1}\ci\Pi_*(g)\bigr)\quad\text{in
$\La$.}
\label{ai6eq3}
\e
\label{ai6thm1}
\end{thm}

\begin{proof} Choose a constructible set $T\subseteq\fObj^\al_\A
(\K)\t\fObj^\be_\A(\K)$ with $f\ot g$ supported on $T$, and use the
notation of Proposition \ref{ai5prop5}. Since
$\uSF(\fObj_\A,\Up,\La)$ is generated by $[(U,\rho)]$ with $U$ a
variety and $\rho$ representable, arguing as in Proposition
\ref{ai5prop3} and Corollary \ref{ai5cor3} we may write $f\ot g$ and
$f*g$ in the forms \eq{ai5eq9} and \eq{ai5eq12}, for $N_m$ finite
and~$c_{mn}\in\La$.

As $i_\La$ is an algebra isomorphism we have
\e
\bigl(i_\La^{-1}\ci\Pi_*(f)\bigr)\bigl(i_\La^{-1}\ci\Pi_*(g)\bigr)=
i_\La^{-1}(\Pi_*(f)\cdot\Pi_*(g))=i_\La^{-1}\ci(\Pi\t\Pi)_*(f\ot g).
\label{ai6eq4}
\e
Equations \eq{ai5eq9} and \eq{ai5eq12} and relations in
$\uSF(\Spec\K,\Up,\La)$ imply that
\ea
&i_\La^{-1}\ci(\Pi\t\Pi)_*(f\ot g)=\sum_{m\in M,\; n\in
N_m}c_{mn}\Up([W_{mn}])\Up([G_m])^{-1},
\label{ai6eq5}
\\
&i_\La^{-1}\ci\Pi_*(f*g)\!=\!\!\sum_{m\in M,\; n\in N_m}\!\!\!
c_{mn}\Up([W_{mn}])\Up([E^1_m])\Up([G_m])^{-1}\Up([E^0_m])^{-1},
\label{ai6eq6}
\ea
using Assumption \ref{ai2ass}(ii) and that $G_m\lt E^0_m\cong G_m\t
E^0_m$ as $\K$-varieties. But
\e
\Up([E^1_m])\Up([E^0_m])^{-1}=\el^{\dim E^1_m}\,\el^{-\dim E^0_m}
=\el^{-\chi(\be,\al)},
\label{ai6eq7}
\e
by Assumption \ref{ai2ass}, Proposition \ref{ai5prop5}(a),
\eq{ai6eq2} and the fact that $T\subseteq\fObj^\al_\A(\K)
\t\fObj^\be_\A(\K)$. Equation \eq{ai6eq3} now follows from
\eq{ai6eq4}--\eq{ai6eq7}.
\end{proof}

We generalize this to the operations $P_\sIp$ of Definition
\ref{ai5def2}. Theorem \ref{ai6thm1} is the case $(I,\pr)=
(\{1,2\},\le)$ of Theorem~\ref{ai6thm2}.

\begin{thm} Suppose Assumptions \ref{ai2ass} and \ref{ai3ass} hold,
and\/ $\chi:K(\A)\t K(\A)\ra\Z$ is biadditive and satisfies
\eq{ai6eq2}. Let\/ $(I,\pr)$ be a finite poset, $\ka:I\ra\bar C(\A)$
and\/ $f_i\in\uSF(\fObj_\A,\Up,\La)$ be supported on
$\fObj_\A^{\ka(i)}$ for all\/ $i\in I$. Then
\e
i_\La^{-1}\!\ci\!\Pi_*\bigl(P_\sIp(f_i:i\in I)\bigr)\!=\!
\raisebox{-5pt}{\begin{Large}$\displaystyle\Bigl[$\end{Large}}
\prod_{i\ne j\in I:\; i\pr j\!\!\!\!\!\!\!\!\!\!\!\!\!\!\!}
\el^{-\chi(\ka(j),\ka(i))}
\raisebox{-5pt}{\begin{Large}$\displaystyle\Bigr]$\end{Large}}\cdot
\raisebox{-5pt}{\begin{Large}$\displaystyle\Bigl[$\end{Large}}
\prod_{i\in I}i_\La^{-1}\!\ci\!\Pi_*(f_i)
\raisebox{-5pt}{\begin{Large}$\displaystyle\Bigr]$\end{Large}}.
\label{ai6eq8}
\e
\label{ai6thm2}
\end{thm}

\begin{proof} When $\md{I}=0$ or 1, equation \eq{ai6eq8} is
obvious. Suppose by induction that \eq{ai6eq8} holds for $\md{I}\le
n$, and let $I,\pr,\ka$ be as above with $\md{I}=n+1$. Choose $k\in
I$ to be $\pr$-maximal, and define $J=I\sm\{k\}$, $K=\{i\in I:i\pr
k\}$, $L=J\cap K$, and $\phi:K\ra\{1,2\}$ by $\phi(i)=1$ for $i\in
L$ and $\phi(k)=2$. Then a similar proof to \cite[Th.~7.10]{Joyc3}
shows the following is a {\it Cartesian square}:
\e
\begin{gathered}
\xymatrix@C=160pt@R=15pt{
*+[r]{\fM(I,\pr)_\A} \ar[r]_{Q(L,\pr,\{1,2\},\le,\phi)\ci S(I,\pr,K)}
\ar[d]^{\,S(I,\pr,J)\t\si(\{k\})}
& *+[l]{\fM(\{1,2\},\le)_\A} \ar[d]_{\bs\si(\{1\})\t\bs\si(\{2\})\,} \\
*+[r]{\fM(J,\pr,\ka)_\A\t\fObj_\A}
\ar[r]^{\bs\si(L)\t\id_{\fObj_\A}}
& *+[l]{\fObj_\A\t\fObj_\A.}
}
\end{gathered}
\label{ai6eq9}
\e
If $f_i\in\uSF(\fObj_\A,\Up,\La)$ for $i\in I$ are as in the theorem
we have
\ea
i_\La^{-1}&\ci\Pi_*\bigl(P_\sIp(f_i:i\in
I)\bigr)=i_\La^{-1}\ci\Pi_*\ci \bs\si(I)_*\bigl[\bigl(\ts\prod_{i\in
I}\bs\si(\{i\})\bigr)^* \bigl(\bigot_{i\in I}f_i\bigr)\bigr]
\nonumber\\
&=i_\La^{-1}\ci\Pi_*\ci\bs\si(\{1,2\})_*\ci\bigl(Q(L,\pr,\{1,2\},\le,\phi)
\ci S(I,\pr,K)\bigr)_*\ci
\nonumber\\
&\qquad\qquad
\bigl(S(I,\pr,J)\t\si(\{k\})\bigr)^*\bigl[\bigl(\ts\prod_{j\in J}
\bs\si(\{j\})\bigr)^*(\bigot_{j\in J}f_j)\bigr)\ot f_k\bigr]
\allowdisplaybreaks
\nonumber\\
&=i_\La^{-1}\ci\Pi_*\ci\bs\si(\{1,2\})_*\ci\bigl(\bs\si(\{1\})\t
\bs\si(\{2\})\bigr)^*\ci
\nonumber\\
&\qquad\qquad
\bigl(\bs\si(L)\t\id_{\fObj_\A}\bigr)_*\bigl[\bigl(\ts\prod_{j\in J}
\bs\si(\{j\})\bigr)^*(\bigot_{j\in J}f_j)\bigr)\ot f_k\bigr]
\allowdisplaybreaks
\nonumber\\
&=i_\La^{-1}\ci\Pi_*\bigl\{\bigl(\bs\si(L)_*\bigl[(\ts\prod_{j\in J}
\bs\si(\{j\}))^*(\bigot_{j\in J}f_j)\bigr]\bigr)*f_k\bigr\}
\allowdisplaybreaks
\nonumber\\
&=\el^{-\chi(\ka(k),\ka(L))}\bigl(i_\La^{-1}\!\ci\!\Pi_*\!\ci\!
\bs\si(L)_*\bigl[(\ts\prod_{j\in J}\bs\si(\{j\}))^* (\bigot_{j\in
J}f_j)\bigr]\bigr) \bigl(i_\La^{-1}\ci\Pi_*(f_k)\bigr)
\nonumber\\
&=\bigl[\ts\prod_{i\in
L}\el^{-\chi(\ka(k),\ka(i))}\bigr]\bigl(i_\La^{-1} \ci\Pi_*\ci
P_\sJp(f_j:j\in J)\bigr)\bigl(i_\La^{-1}\ci\Pi_*(f_k)\bigr).
\label{ai6eq10}
\ea

Here we have used \eq{ai5eq3} in the first step, 2-isomorphisms
\begin{align*}
\Pi\ci\bs\si(I)&\cong\Pi\ci\bs\si(\{1,2\})\ci Q(L,\pr,\{1,2\},\le,\phi)
\ci S(I,\pr,K)\quad\text{and}\\
\ts\prod_{i\in I}\bs\si(\{i\})&\cong\bigl((\ts\prod_{j\in J}\bs\si
(\{j\})))\t\id_{\fObj_\A}\bigr)\ci\bigl(S(I,\pr,J)\t\si(\{k\})\bigr)
\end{align*}
in the second, Theorem \ref{ai2thm2} and \eq{ai6eq9} Cartesian in
the third, \eq{ai5eq1} in the fourth, Theorem \ref{ai6thm1} and
$\bs\si(L)_*[(\prod_{j\in J}\bs\si(\{j\}))^*(\bigot_{j\in J}f_j)]$
supported on $\fObj_\A^{\ka(L)}$ in the fifth, and $\Pi\ci\bs\si(L)
\cong\Pi\ci\bs\si(J)$ and \eq{ai5eq3} in the sixth. Since $\md{J}=n$
we can expand $i_\La^{-1}\ci\Pi_*\ci P_\sJp(f_j:j\in J)$ in the last
line of \eq{ai6eq10} using \eq{ai6eq8} with $J$ in place of $I$, and
this proves \eq{ai6eq8} for $I$. The theorem follows by induction.
\end{proof}

\subsection{Algebras $A(\A,\La,\chi)$ and morphisms to them}
\label{ai62}

If the factor $\el^{-\chi(\be,\al)}$ were not there, equation
\eq{ai6eq3} would say $i_\La^{-1}\ci\Pi_*$ is a morphism of
$\Q$-algebras. We can make an algebra morphism
$\Phi^{\sst\La}:\uSF(\fObj_\A,\Up,\La)\ra A(\A,\La,\chi)$ by
introducing generators $a^\al$ in $A(\A,\La,\chi)$ for $\al\in \bar
C(\A)$, and twisting multiplication $a^\al\star a^\be$ in
$A(\A,\La,\chi)$ by~$\el^{-\chi(\be,\al)}$.

\begin{dfn} Let Assumptions \ref{ai2ass} and \ref{ai3ass} hold.
Then $K(\A)$ is an abelian group, $\bar C(\A)\subseteq K(\A)$ closed
under addition, $\La$ a commutative $\Q$-algebra, and $\el\in\La$ is
invertible. Suppose $\chi:K(\A)\t K(\A)\ra\Z$ is a biadditive map.
Using only this data $K(\A),\bar C(\A),\La,
\el,\chi$ we will define a $\Q$-{\it algebra}~$A(\A,\La,\chi)$.

Let $a^\al$ for $\al\in\bar C(\A)$ be formal symbols, and
$A(\A,\La,\chi)$ the $\La$-module with basis $\{a^\al:\al\in\bar
C(\A)\}$. That is, $A(\A,\La,\chi)$ is the set of sums
$\sum_{\al\in\bar C(\A)}\la^\al a^\al$ with $\la^\al\in\La$ nonzero
for only finitely many $\al$. Addition, and multiplication by $\Q$,
are defined in the obvious way. Define a {\it multiplication}
$\star$ on $A(\A,\La,\chi)$ by
\e
\bigl(\ts\sum_{i\in I}\la_i\,a^{\al_i}\bigr)\star
\bigl(\ts\sum_{j\in J}\mu_j\,a^{\be_j}\bigr)= \ts\sum_{i\in
I}\sum_{j\in J}\la_i\mu_j\el^{-\chi(\be_j,\al_i)}\, a^{\al_i+\be_j},
\label{ai6eq11}
\e
where $I,J$ are finite indexing sets, $\la_i,\mu_j\in\La$ and
$\al_i,\be_j\in\bar C(\A)$. Using the biadditivity of $\chi$ it
is easy to verify $\star$ is associative, and makes $A(\A,\La,\chi)$
into a $\Q$-algebra (in fact, a $\La$-algebra), with identity~$a^0$.

For $(I,\pr)$ a finite poset, define $P_\sIp:\prod_{i\in I}
A(\A,\La,\chi)\ra A(\A,\La,\chi)$ by
\begin{gather}
\text{
\begin{footnotesize}
$\displaystyle
P_\sIp\raisebox{-5pt}{\begin{Large}$\displaystyle\Bigl[$\end{Large}}
\sum_{c_i\in C_i\!\!\!\!}\mu^{c_i}_ia^{\al^{c_i}_i}\!:i\!\in\!I
\raisebox{-5pt}{\begin{Large}$\displaystyle\Bigr]$\end{Large}}\!=
\!\!\!\!
\sum_{\substack{\text{choices of}\\ \text{$c_i\in C_i$ for}\\
\text{all $i\in I$}}}
\raisebox{-5pt}{\begin{Large}$\displaystyle\Bigl[$\end{Large}}
\prod_{i\in I\!\!}\mu^{c_i}_i
\raisebox{-5pt}{\begin{Large}$\displaystyle\Bigr]$\end{Large}}
\!\cdot\!
\raisebox{-5pt}{\begin{Large}$\displaystyle\Bigl[$\end{Large}}
\prod_{i\ne j\in I:\; i\pr j\!\!\!\!\!\!\!\!\!\!\!\!\!\!\!}
\el^{-\chi(\al^{c_j}_j,\al^{c_i}_i)}
\raisebox{-5pt}{\begin{Large}$\displaystyle\Bigr]$\end{Large}}
a^{\sum_{i\in I}a^{\al_i^{c_i}}},
$
\end{footnotesize}
}
\nonumber \\[-20pt]
\label{ai6eq12}
\end{gather}
for $C_i$ finite indexing sets, $\mu^{c_i}_i\in\La$ and
$\al^{c_i}_i\in\bar C(\A)$. These $P_\sIp$ satisfy~\eq{ai4eq23}.

Now for each $\al\in\bar C(\A)$ write
$i_\al:\fObj_\A^\al\ra\fObj_\A$ for the inclusion and
$\Pi^\al:\fObj_\A^\al\ra\Spec\K$ for the projection 1-morphisms, and
let $i_\La$ be as in Proposition \ref{ai2prop1}. Define
$\Phi^{\sst\La}:\uSF(\fObj_\A,\Up,\La)\ra A(\A,\La,\chi)$ by
\e
\Phi^{\sst\La}(f)=\ts\sum_{\al\in\bar
C(\A)}\bigl[i_\La^{-1}\ci\Pi^\al_*\ci
i_\al^*(f)\bigr]a^\al\quad\text{for $f\in\uSF(\fObj_\A,\Up,\La)$.}
\label{ai6eq13}
\e
This is well-defined as $f$ is supported on the disjoint union of
$\fObj_\A^\al$ over finitely many $\al\in\bar C(\A)$, so
$[i_\La^{-1}\ci\Pi^\al_*\ci i_\al^*(f)]\ne 0$ in $\La$ for only
finitely many~$\al$.
\label{ai6def1}
\end{dfn}

We can think of $\Phi^{\sst\La}(f)$ as encoding the `integral' of
$f$ over $\fObj_\A^\al$ for all~$\al$.

\begin{thm} Let Assumptions \ref{ai2ass} and \ref{ai3ass} hold
and\/ $\chi:K(\A)\t K(\A)\ra\Z$ be biadditive and satisfy
\eq{ai6eq2}. Then $\Phi^{\sst\La}:\uSF(\fObj_\A,\Up,\La)\ra
A(\A,\La,\chi)$ is a $\La$-algebra morphism. If\/ $(I,\pr)$ is a
finite poset and $f_i\in\uSF(\fObj_\A,\Up,\La)$ for $i\in I$
then~$\Phi^{\sst\La}\bigl(P_\sIp[f_i:i\in
I]\bigr)=P_\sIp\bigl[\Phi^{\sst\La}(f_i):i\in I\bigr]$.
\label{ai6thm3}
\end{thm}

\begin{proof} Suppose $f^\al,g^\be\in\uSF(\fObj_\A,\Up,\La)$ are
supported on $\fObj_\A^\al,\fObj_\A^\be$ respectively, for
$\al,\be\in\bar C(\A)$. Then for $\ga\in\bar C(\A)$ we have
\begin{equation*}
\Pi^\ga_*\ci i_\ga^*(f^\al)=\begin{cases} \Pi_*(f^\al), & \al=\ga, \\
0, & \al\ne\ga,\end{cases} \;\>\text{so}\;\>
\Phi^{\sst\La}(f^\al)=\bigl[i_\La^{-1}\ci\Pi_*(f^\al)\bigr]a^\al
\;\>\text{by \eq{ai6eq13}.}
\end{equation*}
Similarly
$\Phi^{\sst\La}(g^\be)=\bigl[i_\La^{-1}\ci\Pi_*(g^\be)\bigr] a^\be$
and
$\Phi^{\sst\La}(f^\al*g^\be)=\bigl[i_\La^{-1}\ci\Pi_*(f^\al*g^\be)\bigr]
a^{\al+\be}$, as $f^\al*g^\be$ is supported on $\fObj_\A^{\al+\be}$.
Thus
$\Phi^{\sst\La}(f^\al*g^\be)=\Phi^{\sst\La}(f^\al)\star\Phi^{\sst\La}(g^\be)$
follows from equations \eq{ai6eq3} and \eq{ai6eq11}. For the general
case, any $f,g\in\uSF(\fObj_\A,\Up,\La)$ may be written as
$f=\sum_{\al\in S}f^\al$, $g=\sum_{\be\in T}g^\be$ with
$S,T\subset\bar C(\A)$ finite and $f^\al,g^\be$ supported on
$\fObj_\A^\al,\fObj_\A^\be$, so
$\Phi^{\sst\La}(f*g)=\Phi^{\sst\La}(f)\star\Phi^{\sst\La}(g)$
follows by linearity. Clearly $\Phi^{\sst\La}$ is $\La$-linear and
$\Phi^{\sst\La}(\bde_{[0]})=a^0$, so $\Phi^{\sst\La}$ is an algebra
morphism. The $P_\sIp$ equation is proved in the same way, using
\eq{ai6eq8} rather than~\eq{ai6eq3}.
\end{proof}

\begin{rem}{\bf(a)} Since $\Pi^{\Up,\La}_{\fObj_\A}:\uSF(\fObj_\A)\ra
\uSF(\fObj_\A,\Up,\La)$ is an algebra morphism, $\Phi^{\sst\La}\ci
\Pi^{\Up,\La}_{\fObj_\A}:\uSF(\fObj_\A)\ra A(\A,\La,\chi)$ is also
an algebra morphism, which commutes with the $P_\sIp$. The same
applies to morphisms $\Phi^{\sst\La}\ci\Pi^{\Up,\La}_{\fObj_\A}$
from $\uoSF(\fObj_\A,\Up,\La)$ and $\uoSF(\fObj_\A,\Up,\La^\ci)$,
and all the subalgebras $\SF(\fObj_\A)$, $\SFa(\fObj_\A),\ldots$. We
defined $\Phi^{\sst\La}$ on $\uSF(\fObj_\A,\Up,\La)$ as it is the
coarsest choice.

\noindent{\bf(b)} Suppose $\el$ admits a {\it square root\/} $\wp$
in $\La$, that is, $\wp^2=\el$. In Example \ref{ai2ex} we can take
$\wp=z$. Define elements $\ti a^\al=\wp^{-\chi(\al,\al)}a^\al$ in
$A(\A,\La,\chi)$ for $\al\in\bar C(\A)$. Then the $\ti a^\al$ are an
alternative basis for $A(\A,\La,\chi)$ over $\La$, and
\begin{equation*}
\ti a^\al\star\ti a^\be=\wp^{\chi(\al,\be)-\chi(\be,\al)}\ti a^{\al+\be},
\end{equation*}
by \eq{ai6eq11}. This depends only on the {\it antisymmetrization}
of $\chi$. (This is not true for the $P_\sIp$ in the $\ti a^\al$
basis, though). If $\chi$ is symmetric, with
$\chi(\al,\be)\equiv\chi(\be,\al)$, then $A(\A,\La,\chi)$ is
commutative with~$\ti a^\al\star\ti a^\be=\ti a^{\al+\be}$.
\label{ai6rem1}
\end{rem}

\subsection{Algebras $B(\A,\La,\chi)$ and morphisms to them}
\label{ai63}

Theorem \ref{ai5thm5} gave a compatibility between multiplication
$*$ in $\oSFa(\fObj_\A,\Up,\La)$ and the projections
$\Pi_{[I,\ka]}$. We now exploit this to construct a larger algebra
$B(\A,\La,\chi)$, with an algebra
morphism~$\Psi^{\sst\La}:\oSFa(\fObj_\A,\Up,\La)\ra B(\A,\La,\chi)$.

\begin{dfn} Let Assumptions \ref{ai2ass} and \ref{ai3ass} hold, and
$\chi:K(\A)\t K(\A)\ra\Z$ be a biadditive map. Using only the data
$K(\A),C(\A),\La,\el,\chi$ we will define a $\Q$-{\it algebra}
$B(\A,\La,\chi)$. Consider pairs $(I,\ka)$ for $\ka:I\ra C(\A)$ and
$\approx$-equivalence classes $[I,\ka]$ as in Definition
\ref{ai5def4}. Introduce formal symbols $b_{[I,\ka]}$ for all such
equivalence classes $[I,\ka]$. Let $B(\A,\La,\chi)$ be the
$\La$-module with basis the $b_{[I,\ka]}$. That is, $B(\A,\La,\chi)$
is the set of sums $\sum_{\text{classes
$[I,\ka]$}}\be_{[I,\ka]}b_{[I,\ka]}$ with $\be_{[I,\ka]}\in\La$
nonzero for only finitely many $[I,\ka]$. Addition, and
multiplication by $\Q$, are defined in the obvious way. Define a
{\it multiplication} $\star$ on $B(\A,\La,\chi)$ by
\begin{gather}
b_{[I,\ka]}\star b_{[J,\la]}=
\sum_{\text{eq. classes $[K,\mu]$}}b_{[K,\mu]}\,\cdot\,
\frac{(\el-1)^{\md{K}-\md{I}-\md{J}}}{\md{\Aut(K,\mu)}}\,\cdot
\label{ai6eq14}
\\
\raisebox{-9pt}{\begin{Large}$\displaystyle\biggl[$\end{Large}}
\sum_{\substack{\text{iso.}\\ \text{classes}\\ \text{of finite}\\
\text{sets $L$}}}\!\!
\frac{(-1)^{\md{L}-\md{K}}}{\md{L}!}\!\!\!\!\!\!\!\!
\sum_{\substack{\text{$\phi:I\ra L$, $\psi:J\ra L$ and}\\
\text{$\th:L\!\ra\!K$: $\phi\!\amalg\!\psi$ surjective,}\\
\text{$\mu(k)=\ka((\th\ci\phi)^{-1}(k))+$}\\
\text{$\la((\th\ci\psi)^{-1}(k))$, $k\in K$}}}\,\prod_{k\in
K}\!(\md{\th^{-1}(k)}\!-\!1)!
\prod_{\begin{subarray}{l} i\in I,\; j\in J:\\
\phi(i)=\psi(j)\end{subarray} \!\!\!\!\!\!\!\!\!\!\!\!\!\!\!\! }
\el^{-\chi(\la(j),\ka(i))}
\raisebox{-9pt}{\begin{Large}$\displaystyle\biggr]$\end{Large}}.
\nonumber
\end{gather}
extended $\La$-bilinearly. An elementary but lengthy combinatorial
calculation shows $\star$ is {\it associative}, and makes
$B(\A,\La,\chi)$ into a $\Q$-{\it algebra} (in fact, a $\La$-{\it
algebra}), with identity $b_{[\emptyset,\emptyset]}$, writing
$\emptyset$ for the trivial map $\emptyset\ra C(\A)$. Define
$B^\ind(\A,\La,\chi)$ to be the subspace of
$\sum_{[I,\ka]}\be_{[I,\ka]}b_{[I,\ka]}$ in $B(\A,\La,\chi)$ with
$\be_{[I,\ka]}=0$ unless $\md{I}=1$. Equation \eq{ai6eq14} implies
$B^\ind(\A,\La,\chi)$ is closed under the Lie bracket $[b,c]=b\star
c-c\star b$, and so is a $\Q$- or $\La$-{\it Lie algebra}.

Let $f\in\oSFa(\fObj_\A,\Up,\La)$, so that $\Pi_{[I,\ka]}(f)\in
\oSFa(\fObj_\A,\Up,\La)$ and $\Pi_*\ci\Pi_{[I,\ka]}(f)\in
\uoSF(\Spec\K,\Up,\La)$. Using the explicit form
\cite[Prop.~5.23]{Joyc2} for $\uoSF(\Spec\K,\Up,\La)$ and properties
of the $\Pi_{[I,\ka]}$ we find that
\e
\Pi_*\ci\Pi_{[I,\ka]}(f)=\be_{[I,\ka]}[\Spec\K/(\K^\t)^{\md{I}}],
\label{ai6eq15}
\e
for some unique $\be_{[I,\ka]}\in\La$. Now define
$\Psi^{\sst\La}:\oSFa(\fObj_\A,\Up,\La)\ra B(\A,\La,\chi)$ by
\e
\Psi^{\sst\La}(f)=\ts\sum_{\text{eq. classes
$[I,\ka]$}}\be_{[I,\ka]}b_{[I,\ka]}.
\label{ai6eq16}
\e
Definition \ref{ai5def5} implies that $\Psi^{\sst\La}$ maps
$\oSFai(\fObj_\A,\Up,\La)$ to~$B^\ind(\A,\La,\chi)$.
\label{ai6def2}
\end{dfn}

\begin{thm} Let Assumptions \ref{ai2ass} and \ref{ai3ass} hold
and\/ $\chi:K(\A)\t K(\A)\ra\Z$ be biadditive and satisfy
\eq{ai6eq2}. Then $\Psi^{\sst\La}:\oSFa(\fObj_\A,\Up,\La)\ra
B(\A,\La,\chi)$ and\/ $\Psi^{\sst\La}:\oSFai(\fObj_\A,\Up,\La)
\!\ra\!B^\ind(\A,\La,\chi)$ are (Lie) algebra morphisms.
\label{ai6thm4}
\end{thm}

\begin{proof} Clearly $\Psi^{\sst\La}$ is bilinear and
$\Psi^{\sst\La}(\bde_{[0]})=b_{[\emptyset,\emptyset]}$, so that
$\Psi^{\sst\La}$ takes the identity to the identity. Thus the
theorem follows from
$\Psi^{\sst\La}(f*g)=\Psi^{\sst\La}(f)\star\Psi^{\sst\La}(g)$ for
all $f,g\in\oSFa(\fObj_\A,\Up,\La)$. By Proposition \ref{ai5prop4}
and bilinearity it is enough to show that $\Psi^{\sst\La}(f*g)=
\Psi^{\sst\La}(f)\star\Psi^{\sst\La}(g)$ when $\Pi_{[I,\ka]}(f)=f$
and $\Pi_{[J,\la]}(g)=g$ for some $[I,\ka]$ and $[J,\la]$. Thus we
can apply Theorem \ref{ai5thm5} to get representations \eq{ai5eq30}
and \eq{ai5eq33} for $f\ot g$ and~$\Pi_{[K,\mu]}(f*g)$.

We have $\Psi^{\sst\La}(f)=\be_{[I,\ka]}b_{[I,\ka]}$ and
$\Psi^{\sst\La}(g)=\ga_{[J,\la]}b_{[J,\la]}$ for some
$\be_{[I,\ka]},\ga_{[J,\la]}\in\La$. Then $\Pi_*(f\ot
g)=\be_{[I,\ka]}\ga_{[J,\la]}[\Spec\K/(\K^\t)^I\t(\K^\t)^J]$ in
$\uoSF(\Spec\K,\Up,\La)$. Projecting \eq{ai5eq30} to
$\uoSF(\Spec\K,\Up,\La)$, we deduce that
\e
\be_{[I,\ka]}\ga_{[J,\la]}=\ts\sum_{m\in M,\; n\in
N_m}c_{mn}\Up([W_{mn}]).
\label{ai6eq17}
\e
Now write
$\Psi^{\sst\La}(f*g)=\sum_{[K,\mu]}\de_{[K,\mu]}b_{[K,\mu]}$ for
$\de_{[K,\mu]}\in\La$. Then
$\Pi_*\ci\Pi_{[K,\mu]}(f*g)=\de_{[K,\mu]}[\Spec\K/(\K^\t)^K]$ in
$\uoSF(\Spec\K,\Up,\La)$. Applying $\Pi^{\Up,\La}_{\Spec\K}$ to map
to $\uSF(\Spec\K,\Up,\La)$ and $i_\La^{-1}$ to map to $\La$ gives
\begin{equation*}
i_\La^{-1}\ci\Pi^{\Up,\La}_{\Spec\K}\ci\Pi_*\ci\Pi_{[K,\mu]}(f*g)
=\de_{[K,\mu]}(\el-1)^{-\md{K}}.
\end{equation*}
So applying $i_\La^{-1}\ci\Pi^{\Up,\La}_{\Spec\K}\ci\Pi_*$ to
\eq{ai5eq33} and relations in $\uSF(\Spec\K,\Up,\La)$ gives
\begin{gather}
\de_{[K,\mu]}\!=\!\frac{(\el-1)^{\md{K}-\md{I}-\md{J}}}{
\md{\Aut(K,\mu)}}\!\!\!
\sum_{\substack{\text{iso.}\\ \text{classes}\\ \text{of finite}\\
\text{sets $L$}}}\!\!
\frac{(-1)^{\md{L}-\md{K}}}{\md{L}!}\!\!\!\!\!\!
\sum_{\substack{\text{$\phi:I\ra L$, $\psi:J\ra L$ and}\\
\text{$\th:L\!\ra\!K$: $\phi\!\amalg\!\psi$ surjective,}\\
\text{$\mu(k)=\ka((\th\ci\phi)^{-1}(k))+$}\\
\text{$\la((\th\ci\psi)^{-1}(k))$, $k\in K$}}}
\label{ai6eq18}
\\
\prod_{k\in K}\!(\md{\th^{-1}(k)}\!-\!1)!
\sum_{\!m\in M,\; n\in N_m \!\!\!\!\!\!\!\!\!\!\!}c_{mn}
\Up\bigl([W_{mn}]\bigr)\Up\bigl([(E^1_m)^{T_{L,\phi,\psi}}]\bigr)
\Up\bigl([(E^0_m)^{T_{L,\phi,\psi}}]\bigr)^{-1}. \nonumber
\end{gather}

Let $m,n,L,\phi,\psi$ be as in \eq{ai6eq18}, and pick $w\in
W_{mn}(\K)$ projecting to $v\in V_m(\K)$ and $([X],[Y])\in
\fObj_\A(\K)\t\fObj_\A(\K)$. Then Theorem \ref{ai5thm5} gives
splittings $X\cong\bigop_{i\in I}X_i$ and $Y\cong\bigop_{j\in J}Y_j$
in $\A$ with $[X_i]=\ka(i)$ and $[Y_j]=\la(j)$ in $C(\A)$ for all
$i,j$, and isomorphisms \eq{ai5eq34}. Combining these with
\eq{ai6eq2} yields
\e
\Up\bigl([(E^1_m)^{T_{L,\phi,\psi}}]\bigr)
\Up\bigl([(E^0_m)^{T_{L,\phi,\psi}}]\bigr)^{-1}=\ts\prod_{i\in I,\;
j\in J:\phi(i)=\psi(j)}
\el^{-\chi(\la(j),\ka(i))}.
\label{ai6eq19}
\e
We now see that
$\Psi^{\sst\La}(f*g)=\Psi^{\sst\La}(f)\star\Psi^{\sst\La}(g)$ by
comparing \eq{ai6eq14} and~\eq{ai6eq17}--\eq{ai6eq19}.
\end{proof}

\begin{rem}{\bf(a)} We can also generalize \eq{ai6eq12} to
$P_\sIp:\prod_{i\in I}B(\A,\La,\chi)\ra B(\A,\La,\chi)$ satisfying
\eq{ai4eq23} and $\Psi^{\sst\La}\bigl(P_\sIp[f_i:i\in
I]\bigr)=P_\sIp \bigl[\Psi^{\sst\La}(f_i):i\in I\bigr]$ for
$f_i\in\oSFa(\fObj_\A,\Up,\La)$, as in Theorem \ref{ai6thm3}. But
since the definition and proof are rather complicated, we omit them.

\noindent{\bf(b)} As $\bar\Pi^{\Up,\La}_{\fObj_\A}:
\SFa(\fObj_\A)\!\ra\!\oSFa(\fObj_\A,\Up,\La)$ is an algebra morphism
taking $\SFai(\fObj_\A)\!\ra\!\oSFai(\fObj_\A,\Up,\La)$,
$\Psi^{\sst\La}\ci
\bar\Pi^{\Up,\La}_{\fObj_\A}:\SFa(\fObj_\A)\!\ra\!B(\A,\La,\chi)$ is
too, and restricts to a Lie algebra morphism
$\SFai(\fObj_\A)\!\ra\!B^\ind(\A,\La,\chi)$. The same holds for
$\bar\Pi^{\Up,\La}_{\fObj_\A}:\oSFa(\fObj_\A,\Up,\La^\ci)\!\ra\!
\oSFa(\fObj_\A,\Up,\La)$. We defined $\Psi^{\sst\La}$ on
$\oSFa(\fObj_\A,\Up,\La)$ as it is the coarsest choice.

\noindent{\bf(c)} Here is an alternative description of
$B(\A,\La,\chi)$, $B^\ind(\A,\La,\chi)$. Define $b^\al\in
B^\ind(\A,\La,\chi)$ for $\al\in C(\A)$ by $b^\al=b_{[\{1\},\al']}$,
where $\al'(1)=\al$. Then $B^\ind(\A,\La,\chi)$ is the $\La$-module
with basis $b^\al$ for $\al\in C(\A)$, and \eq{ai6eq14} yields
\e
[b^\al,b^\be]=\frac{\el^{-\chi(\be,\al)}-\el^{-\chi(\al,\be)}}{
\el-1}\,b^{\al+\be}.
\label{ai6eq20}
\e
Given $\ka:\{1,\ldots,n\}\ra C(\A)$, we can use \eq{ai6eq14} to show
that
\begin{equation*}
b^{\ka(1)}\star b^{\ka(2)}\star\cdots\star b^{\ka(n)}=
\frac{1}{\md{\Aut(\{1,\ldots,n\},\ka)}}\cdot
b_{[\{1,\ldots,n\},\ka]}+
\begin{aligned}\bigl(&\text{terms in $b_{[J,\la]}$}\\
&\text{for $\md{J}<n$.}\bigr)
\end{aligned}
\end{equation*}
By induction on $n$ we find the $b^\al$ generate $B(\A,\La,\chi)$
over $\La$, and \eq{ai6eq20} are the only relations on the $b^\al$
over $\La$. Also \eq{ai6eq20} satisfies the Jacobi identity.

Thus, $B(\A,\La,\chi)$ is the $\La$-algebra generated by the $b^\al$
for $\al\in C(\A)$, with relations \eq{ai6eq20}. Equivalently,
$B^\ind(\A,\La,\chi)$ is the $\La$-Lie algebra with basis $b^\al$
for $\al\in C(\A)$ and relations \eq{ai6eq20}, and $B(\A,\La,\chi)$
is the {\it universal enveloping $\La$-algebra} of
$B^\ind(\A,\La,\chi)$. Note this does {\it not\/} mean
$B(\A,\La,\chi)$ is the universal enveloping $\Q$-algebra of
$B^\ind(\A,\La,\chi)$ as a $\Q$-Lie algebra.

\noindent{\bf(d)} Suppose as in Remark \ref{ai6rem1}(b) that $\el$
admits a square root $\wp$ in $\La$. Define elements $\ti
b^\al=\wp^{\,1-\chi(\al,\al)}b^\al$ in $B^\ind(\A,\La,\chi)$ for
$\al\in C(\A)$. These are another $\La$-basis for
$B^\ind(\A,\La,\chi)$, and
\begin{equation*}
[\ti b^\al,\ti b^\be]=\frac{\wp^{\,\chi(\al,\be)-\chi(\be,\al)}-
\wp^{\,\chi(\be,\al)-\chi(\al,\be)}}{\wp-\wp^{\,-1}}\,\ti
b^{\al+\be}
\end{equation*}
by \eq{ai6eq20}, which depends only on the antisymmetrization of
$\chi$, and is also unchanged by replacing $\wp$ by~$\wp^{\,-1}$.

\noindent{\bf(e)} Define a $\La$-algebra morphism
$\De:B(\A,\La,\chi)\!\ra\!A(\A,\La,\chi)$ by
$\De(b_{[I,\ka]})\ab\!=\!(\el-1)^{-\md{I}}a^{\ka(I)}$. Then
$\Phi^{\sst\La}\ci\Pi^{\Up,\La}_{\fObj_\A}\!=\!\De\ci\Psi^{\sst\La}:
\oSFa(\fObj_\A,\Up,\La)\!\ra\!A(\A,\La,\chi)$.
\label{ai6rem2}
\end{rem}

\subsection{Multiplication in $B(\A,\La,\chi)$ as a sum over graphs}
\label{ai64}

We now rewrite the multiplication law \eq{ai6eq14} in
$B(\A,\La,\chi)$ as a sum over {\it directed graphs} $\Ga$, with
vertex set $I\amalg J$. In \eq{ai6eq21}, `no multiples' means there
are {\it no multiple edges}, that is, at most one edge joins any two
vertices in $\Ga$. By the {\it connected components} of $\Ga$ we
mean the sets of vertices of connected components, which are subsets
of $I\amalg J$. And $b_1(\Ga)$ is the {\it first Betti number}
of~$\Ga$.

We find graphs helpful as we can use topological ideas like
connected, simply-connected and $b_1(\Ga)$. The transformation laws
for Calabi--Yau 3-fold invariants in \cite{Joyc5} will also be
written in terms of sums over graphs, and the author believes these
may have something to do with {\it Feynman diagrams} in physics.
Since $b_1(\Ga)\ge 0$, the rational functions of $\el$ appearing in
\eq{ai6eq21} are {\it continuous} at $\el=1$, and lie in $\La^\ci$.
This will be important in~\S\ref{ai65}.

\begin{thm} Equation \eq{ai6eq14} is equivalent to
\begin{gather}
b_{[I,\ka]}\star b_{[J,\la]}=
\sum_{\text{eq. classes $[K,\mu]$}}b_{[K,\mu]}\,\cdot\,
\frac{1}{\md{\Aut(K,\mu)}}\!\!\!\!
\sum_{\substack{\text{$\eta:I\!\ra\!K$, $\ze:J\!\ra\!K$:}\\
\text{$\mu(k)=\ka(\eta^{-1}(k))+\la(\ze^{-1}(k))$}}} \nonumber
\\
\raisebox{-9pt}{\begin{Large}$\displaystyle\biggl[$\end{Large}}
\sum_{\substack{\text{directed graphs $\Ga$: vertices $I\amalg J$,}\\
\text{edges $\mathop{\bu} \limits^{\sst
i}\ra\mathop{\bu}\limits^{\sst j}$, $i\in I$, $j\in J$, no multiples,}\\
\text{conn. components $\eta^{-1}(k)\amalg\ze^{-1}(k)$, $k\in
K$}}}\!\!\!\!\!\!\!\!\!\!\!\!\!\!\!\!
(\el-1)^{b_1(\Ga)} \prod_{\substack{\text{edges}\\
\text{$\mathop{\bu}\limits^{\sst i}\ra\mathop{\bu}\limits^{\sst
j}$}\\ \text{in $\Ga$}}} \frac{\el^{-\chi(\la(j),\ka(i))}-1}{\el-1}
\raisebox{-9pt}{\begin{Large}$\displaystyle\biggr]$\end{Large}}.
\label{ai6eq21}
\end{gather}
\label{ai6thm5}
\end{thm}

\begin{proof} First we rewrite \eq{ai6eq14} as far as possible as a
product over $k\in K$. To do this, for $L,\phi,\psi,\th$ as in
\eq{ai6eq14} write $\eta=\th\ci\phi$ and $\ze=\th\ci\psi$, and for
all $k\in K$ set $L_k=\th^{-1}(\{k\})$ and $\phi_k:\eta^{-1}(\{k\})
\ra L_k$, $\psi_k:\ze^{-1}(\{k\})\ra L_k$ to be the restrictions of
$\phi,\psi$ to $\eta^{-1}(\{k\}),\ze^{-1}(\{k\})$. Then
$L=\coprod_{k\in K}L_k$. The number of surjective maps $\th:L\ra K$
such that $\md{\th^{-1}(\{k\})}=\md{L_k}$ for $k\in K$ is
$\md{L}!/\prod_{k\in K}\md{L_k}!$. Thus, replacing the choice of
$L,\th$ in \eq{ai6eq14} by the choice of sets $L_k$ for $k\in K$, we
replace the factor $1/\md{L}!$ in \eq{ai6eq14} by $1/\prod_{k\in
K}\md{L_k}!$. Writing other terms in \eq{ai6eq14} as products over
$k\in K$, we deduce:
\begin{gather}
b_{[I,\ka]}\star b_{[J,\la]}=
\sum_{\text{eq. classes $[K,\mu]$}}b_{[K,\mu]}\,\cdot\,
\frac{1}{\md{\Aut(K,\mu)}}\!\!\!\!
\sum_{\substack{\text{$\eta:I\!\ra\!K$, $\ze:J\!\ra\!K$:}\\
\text{$\mu(k)=\ka(\eta^{-1}(k))+\la(\ze^{-1}(k))$}}} \nonumber
\\
\prod_{k\in K}\left[
\begin{aligned}
&(\el-1)^{1-\md{\eta^{-1}(k)}-\md{\ze^{-1}(k)}}\,\cdot
\\[-4pt]
&\sum_{\substack{\text{iso. classes}\\ \text{of finite}\\
\text{sets $L_k$}}}\!\!\!\!\!
\raisebox{-4pt}{$\displaystyle\frac{(-1)^{\md{L_k}-1}}{\md{L_k}}$}
\!\!\!
\sum_{\substack{\text{$\phi_k\!:\!\eta^{-1}(k)\!\ra\!L_k$,}\\
\text{$\psi_k\!:\!\ze^{-1}(k)\!\ra\!L_k$:}\\
\text{$\phi_k\amalg\psi_k$ surjective}}}
\,\prod_{\begin{subarray}{l} i\in \eta^{-1}(k),\; j\in \ze^{-1}(k):\\
\phi_k(i)=\psi_k(j)\end{subarray} \!\!\!\!\!\!\!\!\!\!\!\! }\!
\el^{-\chi(\la(j),\ka(i))}
\end{aligned}
\right].
\label{ai6eq22}
\end{gather}

Next we prove the bottom lines of \eq{ai6eq22} and \eq{ai6eq21} are
equal. We can rewrite the bottom line of \eq{ai6eq21} as a product
over $k\in K$ of sums of weighted, connected digraphs $\Ga_k$ with
vertices $\eta^{-1}(\{k\})\amalg\ze^{-1}(\{k\})$, so it is enough to
show the terms in each product over $k\in K$ are equal. This is
equivalent to proving the case $\md{K}=1$, so that $K=\{k\}$.
Dropping subscripts $k$, we have to prove that
\begin{gather}
(\el-1)^{1-\md{I}-\md{J}}\cdot \!\!\!\!\sum_{\substack{\text{iso.
classes of}\\ \text{finite sets $L$}}}\!\!\!\!\!
\frac{(-1)^{\md{L}-1}}{\md{L}} \!\!\!\!\!\!\!
\sum_{\substack{\text{$\phi\!:\!I\!\ra\!L$, $\psi\!:\!J\!\ra\!L$:}\\
\text{$\phi\amalg\psi$ surjective}}}
\,\,\,\prod_{\begin{subarray}{l} i\in I,\; j\in J:\\
\phi(i)=\psi(j)\end{subarray} \!\!\!\!\!\!\!\!\!\!\!\!\!\! }
\,\,\,\el^{-\chi(\la(j),\ka(i))} \nonumber
\\
=\sum_{\substack{\text{connected directed graphs $\Ga$:}\\
\text{vertices $I\amalg J$, edges $\mathop{\bu}\limits^{\sst
i}\ra\mathop{\bu}\limits^{\sst j}$,}\\ \text{$i\in I$, $j\in J$, no
multiples}}} \!\!\!\!\!\!\!\!\!\!\!\!\!\!\!\!\!\!\!
(\el-1)^{b_1(\Ga)} \,\,\prod_{\substack{\text{edges}\\
\text{$\mathop{\bu}\limits^{\sst i}\ra\mathop{\bu}\limits^{\sst
j}$}\\ \text{in $\Ga$}}}
\frac{\el^{-\chi(\la(j),\ka(i))}-1}{\el-1}\,.
\label{ai6eq23}
\end{gather}

Rewrite the top line of \eq{ai6eq23} as a sum over $\Ga$ as follows.
Replace $\el^{-\chi(\la(j),\ka(i))}$ by
$(\el^{-\chi(\la(j),\ka(i))}-1)+1$ and multiply out the product in
$i,j$ to get a sum of products of $\el^{-\chi(\la(j),\ka(i))}-1$ or
1. Associate a digraph $\Ga$ to each of these by putting in an edge
$\smash{\mathop{\bu}\limits^{\sst i}\ra\mathop{\bu}\limits^{\sst
j}}$ for a factor $\el^{-\chi(\la(j),\ka(i))}-1$, and no edge for a
factor 1. Then edges only join $i,j$ with $\phi(i)=\psi(j)$, so each
connected component of $\Ga$ lies in
$\phi^{-1}(\{l\})\amalg\psi^{-1}(\{l\})$ for some $l\in L$. Thus,
reversing the sums over $\Ga$ and $L,\phi,\psi$, the top line of
\eq{ai6eq23} becomes
\e
\begin{gathered}
(\el-1)^{1-\md{I}-\md{J}}
\sum_{\begin{subarray}{l} \text{directed graphs $\Ga$:}\\
\text{vertices $I\amalg J$, edges $\smash{\mathop{\bu}\limits^{\sst
i}\ra\mathop{\bu}\limits^{\sst j}}$,}\\ \text{$i\in I$, $j\in J$, no
multiples}\end{subarray}} \,\,\,
\prod_{\begin{subarray}{l} \text{edges}\\
\text{$\mathop{\bu}\limits^{\sst i}\ra\mathop{\bu}\limits^{\sst j}$
in $\Ga$}\end{subarray}\!\!\!\!\!\!\!\!\!\!\!\!\!\!\!\!\!
}\,\,\,(\el^{-\chi(\la(j),\ka(i))}-1)\cdot
\\
\raisebox{-9pt}{\begin{Large}$\displaystyle\biggl[$\end{Large}}
\sum_{\substack{\text{iso. classes of}\\ \text{finite sets
$L$}}}\!\!\!\!\! \frac{(-1)^{\md{L}-1}}{\md{L}} \!\!\!\!\!\!\!
\sum_{\substack{\text{$\phi\!:\!I\!\ra\!L$, $\psi\!:\!J\!\ra\!L$:
$\phi\amalg\psi$ surjective, conn.}\\
\text{components of $\Ga$ lie in $\phi^{-1}(l)\amalg\psi^{-1}(l)$,
$l\in L$}}} \!\! 1\,
\raisebox{-9pt}{\begin{Large}$\displaystyle\biggr]$\end{Large}}.
\end{gathered}
\label{ai6eq24}
\e

We shall show the bottom line $[\cdots]$ of \eq{ai6eq24} is 1 if
$\Ga$ is connected, and 0 otherwise. If $\Ga$ is connected then
$I\amalg J$ must lie in some $\phi^{-1}(\{l\})\amalg
\psi^{-1}(\{l\})$, so as $\phi\amalg\psi$ is surjective the only
possibility is $L=\{l\}$ and $\md{L}=1$, and the sum reduces to 1.
If $\Ga$ is not connected, fix one connected component $I_0\amalg
J_0$ of $\Ga$, $I_0\subseteq I$, $J_0\subseteq J$. If a triple
$(L,\phi,\psi)$ appears in the bottom line of \eq{ai6eq24} then for
some $l\in L$ we have $\phi\vert_{I_0}\equiv l$,
$\psi\vert_{J_0}\equiv l$. Divide such $(L,\phi,\psi)$ into two
cases: (a) $\phi^{-1}(\{l\})=I_0$ and $\psi^{-1}(\{l\})=J_0$, and
(b) otherwise.

If $(L,\phi,\psi)$ satisfy (b) then define $(L',\phi',\psi')$
satisfying (a) by $L'=L\amalg\{l'\}$ for some $l'\notin L$, and
$\phi'(i)=l'$ for $i\in I_0$ and $\phi'(i)=\phi(i)$ for $i\notin
I_0$, and $\psi'(j)=l'$ for $j\in J_0$ and $\psi'(j)=\psi(j)$ for
$j\notin J_0$. Conversely, if $(L',\phi',\psi')$ satisfy (a) with
$\phi'\vert_{I_0}\equiv\psi'\vert_{J_0}\equiv l'$, then define
$(L,\phi,\psi)$ satisfying (b) by $L=L'\sm\{l'\}$, and choosing some
$l\in L$ define $\phi(i)=l$ for $i\in I_0$ and $\phi(i)=\phi'(i)$
for $i\notin I_0$, and $\psi(j)=l$ for $j\in J_0$ and
$\psi(j)=\psi'(j)$ for~$j\notin J_0$.

This establishes a 1-1 correspondence between $(L,\phi,\psi)$
satisfying (b), and quadruples $L',\phi',\psi',l$ with special
element $l'\in L'$, such that $(L',\phi',\psi')$ satisfies (a) with
$\phi'\vert_{I_0}\equiv\psi' \vert_{J_0}\equiv l'\in L'$, and $l\in
L'\sm\{l'\}$. We also have $\md{L'}=\md{L}+1$. Thus, the terms in
the bottom line of \eq{ai6eq24} corresponding to $(L,\phi,\psi)$ and
$(L',\phi',\psi')$ differ by a factor~$-\md{L}/\md{L'}$.

Now the equation $\phi'\vert_{I_0}\equiv\psi' \vert_{J_0}\equiv
l'\in L$ constrains the choice of $\phi',\psi'$, fixing $l'$ out of
$\md{L'}$ choices of points in $L'$, and this accounts for the
factor $1/\md{L'}$. And given $L',\phi',\psi'$ there are $\md{L}$
possible choices for $l\in L'\sm\{l'\}$, which accounts for the
factor $\md{L}$. Because of these, for disconnected $\Ga$ the
contributions of $L,\phi,\psi$ of types (a) and (b) in the bottom
line of \eq{ai6eq24} cancel, giving zero.

Thus \eq{ai6eq24} reduces to the bottom line of \eq{ai6eq23}, except
for the powers of $\el-1$, which are $b_1(\Ga)-n$ in the bottom line
of \eq{ai6eq23}, where $n$ is the number of edges in $\Ga$, and
$1-\md{I}-\md{J}$ in \eq{ai6eq24}. But $\Ga$ is connected with
$\md{I}+\md{J}$ vertices and $n$ edges, so $b_1(\Ga)-n=
1-\md{I}-\md{J}$. This proves \eq{ai6eq23}, and hence~\eq{ai6eq21}.
\end{proof}

\subsection{Algebras $B(\A,\La^\ci,\chi),C(\A,\Om,\chi)$ and
morphisms to them}
\label{ai65}

We define algebras $B(\A,\La^\ci,\chi),C(\A,\Om,\chi)$ and a
morphism $\Pi$ between them.

\begin{dfn} Let Assumptions \ref{ai2ass} and \ref{ai3ass} hold,
$\chi:K(\A)\t K(\A)\ra\Z$ or $\Q$ be biadditive, and use the
notation of Definition \ref{ai6def2}. Define $B(\A,\La^\ci,\chi)$ to
be the subspace of $\sum_{[I,\ka]}\be_{[I,\ka]}b_{[I,\ka]}$ in
$B(\A,\La,\chi)$ with all $\be_{[I,\ka]}\in\La^\ci$. Since the
bottom line of \eq{ai6eq21} lies in $\La^\ci$ we see that
$B(\A,\La^\ci,\chi)$ is closed under $\star$, and so is a
$\La^\ci$-{\it subalgebra} of $B(\A,\La,\chi)$, with $\La^\ci$-{\it
basis} the $b_{[I,\ka]}$. Define
$B^\ind(\A,\La^\ci,\chi)=B^\ind(\A,\La,\chi)\cap
B(\A,\La^\ci,\chi)$. Then $B^\ind(\A,\La^\ci,\chi)$ is a
$\La^\ci$-{\it Lie subalgebra} of $B(\A,\La^\ci,\chi)$, since
$B^\ind(\A,\La,\chi)$ is a $\La$-Lie subalgebra in~$B(\A,\La,\chi)$.

As in Definition \ref{ai6def2}, introduce symbols $c_{[I,\ka]}$ for
all equivalence classes $[I,\ka]$, and let $C(\A,\Om,\chi)$ be the
$\Om$-module with basis the $c_{[I,\ka]}$. That is, $C(\A,\Om,\chi)$
is the set of sums $\sum_{\text{classes
$[I,\ka]$}}\ga_{[I,\ka]}c_{[I,\ka]}$ with $\ga_{[I,\ka]}\in\Om$
nonzero for only finitely many $[I,\ka]$. Define a {\it
multiplication} $\star$ on $C(\A,\Om,\chi)$ by
\begin{gather}
c_{[I,\ka]}\star c_{[J,\la]}=
\sum_{\text{eq. classes $[K,\mu]$}}c_{[K,\mu]}\,\cdot\,
\frac{1}{\md{\Aut(K,\mu)}}\!\!\!\!
\sum_{\substack{\text{$\eta:I\!\ra\!K$, $\ze:J\!\ra\!K$:}\\
\text{$\mu(k)=\ka(\eta^{-1}(k))+\la(\ze^{-1}(k))$}}} \nonumber
\\
\raisebox{-9pt}{\begin{Large}$\displaystyle\biggl[$\end{Large}}
\sum_{\substack{\text{simply-connected directed graphs $\Ga$:}\\
\text{vertices $I\amalg J$, edges $\mathop{\bu} \limits^{\sst
i}\ra\mathop{\bu}\limits^{\sst j}$, $i\in I$, $j\in J$,}\\
\text{conn. components $\eta^{-1}(k)\amalg\ze^{-1}(k)$, $k\in K$}}}
\,\,\,\prod_{\substack{\text{edges}\\
\text{$\mathop{\bu}\limits^{\sst i}\ra\mathop{\bu}\limits^{\sst
j}$}\\ \text{in $\Ga$}}}\bigl(-\chi(\la(j),\ka(i))\bigr)
\raisebox{-9pt}{\begin{Large}$\displaystyle\biggr]$\end{Large}},
\label{ai6eq25}
\end{gather}
extended $\Om$-bilinearly in the usual way. This is still
well-defined if $\chi$ takes values in $\Q$ rather than $\Z$, and in
\S\ref{ai66} we allow $\chi=\ha\bar\chi$ to take values in~$\ha\Z$.

Define a $\Q$-linear map $\Pi:B(\A,\La^\ci,\chi)\longra
C(\A,\Om,\chi)$ by
\e
\Pi:\ts\sum_{\text{eq. classes $[I,\ka]$}}\be_{[I,\ka]}b_{[I,\ka]}
\longmapsto\ts\sum_{\text{eq. classes
$[I,\ka]$}}\pi(\be_{[I,\ka]})c_{[I,\ka]},
\label{ai6eq26}
\e
for $\pi:\La^\ci\ra\Om$ as in Assumption \ref{ai2ass}. Then
$\Pi(b_{[I,\ka]})=c_{[I,\ka]}$, and $\Pi$ is {\it surjective} as
$\pi$ is. Comparing \eq{ai6eq21} and \eq{ai6eq25} shows that
\e
\Pi\bigl(b_{[I,\ka]}\bigr)\star\Pi\bigl(b_{[J,\la]}\bigr)
=c_{[I,\ka]}\star c_{[J,\la]}=\Pi\bigl(b_{[I,\ka]}\star
b_{[J,\la]}\bigr).
\label{ai6eq27}
\e
Here we effectively take the limit $\el\ra 1$ in the bottom line
$[\cdots]$ of \eq{ai6eq21} to get the bottom line $[\cdots]$ of
\eq{ai6eq25}, since $\pi(\el)=1$ in $\Om$. The factor
$(\el-1)^{b_1(\Ga)}$ in \eq{ai6eq21} shows that only $\Ga$ with
$b_1(\Ga)=0$ contribute in the limit, that is, only {\it
simply-connected\/} $\Ga$. We drop the condition `no multiples' from
\eq{ai6eq21}, as this is implied by $\Ga$ simply-connected.

Using \eq{ai6eq27}, the associativity of $\star$ in
$B(\A,\La^\ci,\chi)$, and $\pi$ an algebra morphism, we see that
$\star$ in $C(\A,\Om,\chi)$ is {\it associative}, so
$C(\A,\Om,\chi)$ is an $\Om$-{\it algebra}, with identity
$c_{[\emptyset,\emptyset]}$, and $\Pi$ in \eq{ai6eq26} is a
$\Q$-{\it algebra morphism}. Define $C^\ind(\A,\Om,\chi)$ to be the
subspace of $\sum_{[I,\ka]}\ga_{[I,\ka]}c_{[I,\ka]}$ in
$C(\A,\Om,\chi)$ with $\ga_{[I,\ka]}=0$ unless $\md{I}=1$. From
\eq{ai6eq25} we see $C^\ind(\A,\Om,\chi)$ is an $\Om$-{\it Lie
subalgebra} of $C(\A,\Om,\chi)$. Also, $\Pi$ restricts to a
$\Q$-{\it Lie algebra morphism}~$\Pi:B^\ind(\A,\La^\ci,\chi)\ra
C^\ind(\A,\Om,\chi)$.

Now let $f\in\oSFa(\fObj_\A,\Up,\La^\ci)$. As in Definition
\ref{ai6def2} equation \eq{ai6eq15} holds, but this time for
$\be_{[I,\ka]}\in\La^\ci$. Define $\Psi^{\sst\La^\ci}:
\oSFa(\fObj_\A,\Up,\La^\ci)\ra B(\A,\La^\ci,\chi)$ by
$\Psi^{\sst\La^\ci}(f)=\sum_{[I,\ka]}\be_{[I,\ka]}b_{[I,\ka]}$, as
in \eq{ai6eq16}. In the same way, if $f\in\oSFa(\fObj_\A, \Th,\Om)$
then \eq{ai6eq15} holds for $\be_{[I,\ka]}\in\Om$. Define
$\Psi^{\sst\Om}:\oSFa(\fObj_\A,\Th,\Om)\ra C(\A,\Om,\chi)$ by
$\Psi^{\sst\Om}(f)=\sum_{[I,\ka]}\be_{[I,\ka]}c_{[I,\ka]}$. These
restrict to $\Psi^{\sst\La^\ci}:\oSFai(\fObj_\A,\Up,\La^\ci)\ra
B^\ind(\A,\La^\ci,\chi)$ and~$\Psi^{\sst\Om}:\oSFai
(\fObj_\A,\Th,\Om)\ra C^\ind(\A,\Om,\chi)$.
\label{ai6def3}
\end{dfn}

Here is the analogue of Theorem~\ref{ai6thm4}.

\begin{thm} Let Assumptions \ref{ai2ass} and \ref{ai3ass} hold
and\/ $\chi:K(\A)\t K(\A)\ra\Z$ be biadditive and satisfy
\eq{ai6eq2}. Then $\Psi^{\sst\La^\ci}:\oSFa(\fObj_\A,\Up,\La^\ci)\ra
B(\A,\La^\ci,\chi)$ and\/ $\Psi^{\sst\Om}:\oSFa(\fObj_\A,\Th,\Om)\ra
C(\A,\Om,\chi)$ are $\La^\ci,\Om$-algebra morphisms, and\/
$\Psi^{\sst\La^\ci}: \oSFai(\fObj_\A,\Up,\La^\ci)\ra
B^\ind(\A,\La^\ci,\chi)$ and\/ $\Psi^{\sst\Om}:\oSFai
(\fObj_\A,\Th,\Om)\ra C^\ind(\A,\Om,\chi)$ are $\La^\ci,\Om$-Lie
algebra morphisms.
\label{ai6thm6}
\end{thm}

\begin{proof} $\Psi^{\sst\La^\ci}$ is the restriction of
$\Psi^{\sst\La}$ to $\oSFa(\fObj_\A,\Up,\La^\ci)\subset
\oSFa(\fObj_\A,\Up,\La)$, so $\Psi^{\sst\La^\ci}$ is a (Lie) algebra
morphism by Theorem \ref{ai6thm4}. For $\Psi^{\sst\Om}$, note that
\e
\Pi\ci\Psi^{\sst\La^\ci}=\Psi^{\sst\Om}\ci\bar
\Pi^{\Th,\Om}_{\fObj_\A}:\oSFa(\fObj_\A,\Up,\La^\ci)\ra
C(\A,\Om,\chi).
\label{ai6eq28}
\e
Let $f,g\in\oSFa(\fObj_\A,\Th,\Om)$. Since $\bar
\Pi^{\Th,\Om}_{\fObj_\A}$ is surjective we can lift them to
$f',g'\in\oSFa(\fObj_\A,\Up,\La^\ci)$ with $f,g=\bar\Pi^{\Th,
\Om}_{\fObj_\A}(f',g')$. Then
\begin{align*}
\Psi^{\sst\Om}(f*g)&=\Psi^{\sst\Om}\bigl(\bar\Pi^{\Th,
\Om}_{\fObj_\A}(f')*\bar\Pi^{\Th,
\Om}_{\fObj_\A}(g')\bigr)=\Psi^{\sst\Om}\ci\bar\Pi^{\Th,
\Om}_{\fObj_\A}(f'*g')\\
&=\Pi\ci\Psi^{\sst\La^\ci}(f'*g')=
\bigl(\Pi\ci\Psi^{\sst\La^\ci}(f')\bigr)\star
\bigl(\Pi\ci\Psi^{\sst\La^\ci}(g')\bigr)\\
&=\bigl(\Psi^{\sst\Om}\ci\bar\Pi^{\Th,\Om}_{\fObj_\A}(f')\bigr)
\star\bigl(\Psi^{\sst\Om}\ci\bar\Pi^{\Th,\Om}_{\fObj_\A}(g')\bigr)
=\Psi^{\sst\Om}(f)\star\Psi^{\sst\Om}(g),
\end{align*}
using \eq{ai6eq28} and that $\Pi,\Psi^{\sst\La^\ci}$ and $\bar
\Pi^{\Th,\Om}_{\fObj_\A}$ are algebra morphisms. Also
$\Psi^{\sst\Om}(\bde_{[0]})=c_{[\emptyset,\emptyset]}$, so
$\Psi^{\sst\Om}$ is a (Lie) algebra morphism.
\end{proof}

The analogue of Remark \ref{ai6rem2} holds. In particular,
$\Psi^{\sst\La^\ci}\ci\bar\Pi^{\Up,\La^\ci}_{\fObj_\A},
\Psi^{\sst\Om}\ci\bar\Pi^{\Th,\Om}_{\fObj_\A}$ are algebra morphisms
from $\SFa(\fObj_\A)\ra B(\A,\La^\ci,\chi),C(\A,\Om,\chi)$. Define
$c^\al\in C^\ind(\A,\Om,\chi)$ for $\al\in C(\A)$ by
$c^\al=c_{[\{1\},\al']}$, where $\al'(1)=\al$. Then
$C^\ind(\A,\Om,\chi)$ is the $\Om$-module with basis $c^\al$ for
$\al\in C(\A)$, and \eq{ai6eq25} yields
\e
[c^\al,c^\be]=(\chi(\al,\be)-\chi(\be,\al))\,c^{\al+\be}.
\label{ai6eq29}
\e
It will be important in \S\ref{ai66} that this depends only on the
{\it antisymmetrization} of $\chi$. The argument of Remark
\ref{ai6rem2}(c) shows that $B(\A,\La^\ci,\chi)$ and
$C(\A,\Om,\chi)$ are the $\La^\ci$- and $\Om$-enveloping algebras of
$B^\ind(\A,\La^\ci,\chi)$ and~$C^\ind(\A,\Om,\chi)$.

\subsection{Calabi-Yau 3-folds and Lie algebra morphisms}
\label{ai66}

Let $P$ be a smooth projective $\K$-scheme of dimension $m$, and
$\A=\coh(P)$ and $\fF_\A,K(\A)$ be as in \cite[Ex.~9.1]{Joyc3}. Then
as in \S\ref{ai61} there is a bilinear map $\bar\chi:K(\A)\t
K(\A)\ra\Z$ called the {\it Euler form} satisfying
\e
\bar\chi\bigl([X],[Y]\bigr)=\ts\sum_{i=0}^m(-1)^i\dim_\K\Ext^i(X,Y)
\quad\text{for all $X,Y\in\A$.}
\label{ai6eq30}
\e
We denote it $\bar\chi$ to distinguish it from $\chi$ in
\S\ref{ai61}--\S\ref{ai65}. By {\it Serre duality} we have natural
isomorphisms
\e
\Ext^i(X,Y)^*\cong\Ext^{m-i}(Y,X\ot K_P) \quad\text{for all
$X,Y\in\A$ and $i=0,\ldots,m$,}
\label{ai6eq31}
\e
where $K_P$ is the {\it canonical line bundle} of $P$. We call $P$ a
{\it Calabi--Yau $m$-fold\/} if $K_P$ is trivial, so that
\eq{ai6eq31} reduces to~$\Ext^i(X,Y)^*\cong\Ext^{m-i}(Y,X)$.

In particular, if $P$ is a {\it Calabi--Yau $3$-fold\/} then
\eq{ai6eq30} and \eq{ai6eq31} imply that
\e
\begin{split}
&\bigl(\dim_\K\Hom(X,Y)-\dim_\K\Ext^1(X,Y)\bigr)-\\
&\bigl(\dim_\K\Hom(Y,X)-\dim_\K\Ext^1(Y,X)\bigr)=
\bar\chi\bigl([X],[Y]\bigr)\quad\text{for all $X,Y\in\A$,}
\end{split}
\label{ai6eq32}
\e
where $\bar\chi$ is antisymmetric. This is similar to equation
\eq{ai6eq2}, which we used to construct the algebra morphisms
$\Phi^{\sst\La},\Psi^{\sst\La}, \Psi^{\sst\La^\ci},\Psi^{\sst\Om}$
in \S\ref{ai62}--\S\ref{ai65}. In fact, \eq{ai6eq2} implies
\eq{ai6eq32} with $\bar\chi(\al,\be)=\chi(\al,\be)-\chi(\be,\al)$,
so \eq{ai6eq32} is a {\it weakening} of \eq{ai6eq2}, which holds
when $\A=\coh(P)$ for $P$ a smooth projective curve and when
$\A=\modKQ$, as in \S\ref{ai61}, and also when $\A=\coh(P)$ for $P$
a Calabi--Yau 3-fold.

We shall show that \eq{ai6eq32} implies $\Psi^{\sst\Om}:
\oSFai(\fObj_\A,\Th,\Om)\ra C^\ind(\A,\Om,\ha\bar\chi)$ is a Lie
algebra morphism, generalizing Theorem \ref{ai6thm6}. First we
explain why this {\it only} works for the restriction of
$\Psi^{\sst\Om}$ to $\oSFai(\fObj_\A,\Th,\Om)$. That is,
\eq{ai6eq32} is too weak an assumption to make $\Psi^{\sst\Om}$ an
algebra morphism on $\oSFa(\fObj_\A,\Th,\Om)$, nor to make
$\Phi^{\sst\La},\Psi^{\sst\La}$ or $\Psi^{\sst\La^\ci}$ into (Lie)
algebra morphisms.

Consider whether \eq{ai6eq32} could imply $\Phi^{\sst\La}:
\uSF(\fObj_\A,\Up,\La)\ra A(\A,\La,\ha\bar\chi)$ is an algebra
morphism. Let $\chi:K(\A)\t K(\A)\ra\Z$ be bilinear with
$\bar\chi(\al,\be)=\chi(\al,\be)-\chi(\be,\al)$; note that for fixed
$\bar\chi$ there will be {\it many} such $\chi$, differing by
symmetric bilinear forms. As a $\La$-module $A(\A,\La,\chi)$ depends
only on $\A,\La$, so that $A(\A,\La,\chi)=A(\A,\La,\ha\bar\chi)$,
and $\Phi^{\sst\La}$ also depends only on $\A,\La$. It is only the
multiplication $\star$ in $A(\A,\La,\chi)$ which depends on the
choice of~$\chi$.

Now \eq{ai6eq2} for $\chi$ implies \eq{ai6eq32} for $\bar\chi$, as
above. If \eq{ai6eq2} holds then $\Phi^{\sst\La}:
\uSF(\fObj_\A,\Up,\La)\ra A(\A,\La,\chi)$ is an algebra morphism by
Theorem \ref{ai6thm3}. Therefore $\Phi^{\sst\La}:
\uSF(\fObj_\A,\Up,\La)\ra A(\A,\La,\ha\bar\chi)$ {\it cannot\/} be
an algebra morphism in general, since the multiplications
\eq{ai6eq11} on $A(\A,\La,\chi)=A(\A,\La,\ha\bar\chi)$ associated to
$\chi$ and $\ha\bar\chi$ are different. The many choices of $\chi$
giving the same $\bar\chi$ mean there is no one choice of $\star$
for which $\Phi^{\sst\La}$ is an algebra morphism whenever
\eq{ai6eq32} holds.

Similarly $\Psi^{\sst\La},\Psi^{\sst\La^\ci},\Psi^{\sst\Om}$ cannot
be algebra morphisms on $\oSFa(\fObj_\A,*,*)$, nor $\Psi^{\sst\La},
\Psi^{\sst\La^\ci}$ Lie algebra morphisms on $\oSFai(\fObj_\A,*,*)$,
as in each case $\star$ or $[\,,\,]$ in the image varies
nontrivially with $\chi$ giving the same $\bar\chi$. But
$\Psi^{\sst\Om}:\oSFai(\fObj_\A,\Th,\Om)\ra
C^\ind(\A,\Om,\ha\bar\chi)$ is different, since by \eq{ai6eq29} the
Lie bracket $[\,,\,]$ on $C^\ind(\A,\Om,\chi)$ depends only on
$\bar\chi(\al,\be)=\chi(\al,\be)-\chi(\be,\al)$, as we want.

\begin{thm} Let Assumptions \ref{ai2ass} and \ref{ai3ass} hold
and\/ $\bar\chi:K(\A)\t K(\A)\ra\Z$ be biadditive and satisfy
\eq{ai6eq32}. Then $\Psi^{\sst\Om}:\oSFai(\fObj_\A,\Th,\Om)\ra
C^\ind(\A,\Om,\ha\bar\chi)$ in Definition \ref{ai6def3} is an
$\Om$-Lie algebra morphism.
\label{ai6thm7}
\end{thm}

\begin{proof} We must show $\Psi^{\sst\Om}([f,g])=
[\Psi^{\sst\Om}(f),\Psi^{\sst\Om}(g)]$ for
$f,g\in\oSFai(\fObj_\A,\Th,\Om)$. It is enough to prove this for
$f,g$ supported on $\fObj_\A^\al(\K),\fObj_\A^\be(\K)$ respectively
for $\al,\be\in C(\A)$. Lift $f,g$ to $f',g'\in\oSFai(\fObj_\A,\Up,
\La^\ci)$ with $\bar\Pi_{\fObj_\A}^{\Th,\Om}(f',g')=f,g$, which is
possible as $\bar\Pi_{\fObj_\A}^{\Th,\Om}$ is surjective.

Choose constructible $T\subseteq\fObj_\A^\al(\K)\t\fObj_\A^\be(\K)$
with $f\ot g$ and $f'\ot g'$ supported on $T$. Generalizing
Proposition \ref{ai5prop5}, we can find a finite decomposition
$T=\coprod_{m\in M}\fG_m(\K)$, 1-isomorphisms $\fG_m\cong[V_m/G_m]$,
and finite-dimensional $G_m$-representations $E^0_m,E^1_m,\ti
E^0_m,\ti E^1_m$ for $m\in M$, satisfying analogues of Proposition
\ref{ai5prop5}(a)--(d), where in (a) we have isomorphisms:
\e
\begin{aligned}
\Hom(Y,X)&\cong E^0_m,\;& \Ext^1(Y,X)&\cong E^1_m,\\
\Hom(X,Y)&\cong\ti E^0_m,\;& \Ext^1(X,Y)&\cong\ti E^1_m.
\end{aligned}
\label{ai6eq33}
\e

The proof of equations \eq{ai5eq30}--\eq{ai5eq31} in Theorem
\ref{ai5thm5} now shows we may write
\begin{gather}
f'\ot g'=\ts\sum_{m\in M,\; n\in N_m}c_{mn}\bigl[\bigl(W_{mn}\t
[\Spec\K/(\K^\t)^2],\tau_{mn}\bigr)\bigr],
\label{ai6eq34}
\\
[f',g']=\sum_{\!\!\! m\in M,\; n\in N_m \!\!\!\!\!\!\!\!\!}
\begin{aligned}[t]c_{mn}
\bigl\{ &\bigl[\bigl(W_{mn}\!\t\![E^1_m/(\K^\t)^2\!\lt\!E^0_m],
\bs\si(\{1,2\})\!\ci\!\xi_{mn}\bigr)\bigr]\\
-&\bigl[\bigl(W_{mn}\!\t\![\ti E^1_m/(\K^\t)^2\!\lt\!\ti E^0_m],
\bs\si(\{1,2\})\!\ci\!\ti\xi_{mn}\bigr)\bigr]\bigr\},
\end{aligned}
\label{ai6eq35}
\end{gather}
for $N_m$ finite, $c_{mn}\in\La^\ci$, and $W_{mn}$ quasiprojective
$\K$-varieties.

Since $f',g',[f',g']$ lie in $\oSFai(\fObj_\A,\Up,\La^\ci)$ and are
supported on $\fObj_\A^\al(\K)$, $\fObj_\A^\be(\K)$ and
$\fObj_\A^{\al+\be}(\K)$ respectively, we have
\e
\Psi^{\sst\La^\ci}(f')=\de'\,b^\al,\quad \Psi^{\sst\La^\ci}(g')=
\ep'\,b^\be,\quad \Psi^{\sst\La^\ci}([f',g'])=\ze'\,b^{\al+\be},
\label{ai6eq36}
\e
for some $\de',\ep',\ze'\in\La^\ci$, and $b^\al,b^\be,b^{\al+\be}$
as in Remark \ref{ai6rem2}(c). Projecting \eq{ai6eq34} and
\eq{ai6eq35} to $\uoSF(\Spec\K,\Up,\La)$ as in the proofs of
\eq{ai6eq17} and \eq{ai6eq18} then shows that
\ea
\de'\,\ep'&=\ts\sum_{m\in M,\; n\in N_m}c_{mn}\Up([W_{mn}]),
\label{ai6eq37}
\\
\ze'&=\!\!\sum_{m\in M,\; n\in N_m\!\!\!\!\!\!\!\!\!\!\!}
c_{mn}\Up([W_{mn}])\frac{\Up([E^1_m])\Up([E^0_m])^{-1}-\Up([\ti
E^1_m])\Up([\ti E^0_m])^{-1}}{\el-1}\,.
\label{ai6eq38}
\ea

For $m\in M$ we can find $[X]\in\fObj^\al_\A(\K)$ and
$[Y]\in\fObj_\A^\be(\K)$ such that \eq{ai6eq33} holds, and then in
\eq{ai6eq38} we have
\begin{gather*}
\frac{\Up([E^1_m])\Up([E^0_m])^{-1}-\Up([\ti E^1_m])\Up([\ti E^0_m]
)^{-1}}{\el-1}=\el^{\dim\Ext^1(X,Y)-\dim\Hom(X,Y)}\cdot\\
\frac{\el^{\dim\Ext^1(Y,X)- \dim\Hom(Y,X)-
\dim\Ext^1(X,Y)+\dim\Hom(X,Y)}-1}{\el-1}\,.
\end{gather*}
Composing with $\pi:\La^\ci\ra\Om$ and using \eq{ai6eq32} and
$\pi(\el)=1$ then gives
\e
\pi\Bigl(\frac{\Up([E^1_m])\Up([E^0_m])^{-1}-\Up([\ti E^1_m])\Up
([\ti E^0_m])^{-1}}{\el-1}\Bigr)=\bar\chi(\al,\be)\in\Z\subset\Om.
\label{ai6eq39}
\e

Applying $\Pi$ in \eq{ai6eq26} to \eq{ai6eq36} and using
\eq{ai6eq28}, that $\bar\Pi^{\Th,\Om}_{\fObj_\A}$ is an algebra
morphism, and $\bar\Pi_{\fObj_\A}^{\Th,\Om}(f')=f$,
$\bar\Pi_{\fObj_\A}^{\Th,\Om}(g')=g$, we find that
\e
\Psi^{\sst\Om}(f)=\de\,c^\al,\quad \Psi^{\sst\Om}(g)=
\ep\,c^\be,\quad \Psi^{\sst\Om}([f,g])=\ze\,c^{\al+\be},
\label{ai6eq40}
\e
where $\pi(\de')=\de$, $\pi(\ep')=\ep$, $\pi(\ze')=\ze$ and
$c^\al,c^\be,c^{\al+\be}$ are as in \eq{ai6eq29}. Applying $\pi$ to
\eq{ai6eq37}--\eq{ai6eq38} and using \eq{ai6eq39} then shows that
$\ze=\de\,\ep\,\bar\chi(\al,\be)$. Combining this with \eq{ai6eq29}
and \eq{ai6eq40} now shows that $\Psi^{\sst\Om}([f,g])=
[\Psi^{\sst\Om}(f),\Psi^{\sst\Om}(g)]$, where the bracket $[\,,\,]$
in $C^\ind(\A,\Om,\ha\bar\chi)$ is defined using the form
$\chi=\ha\bar\chi$, so that $\chi(\al,\be)-\chi(\be,\al)$ in
\eq{ai6eq29} is $\bar\chi(\al,\be)$ as $\bar\chi$ is antisymmetric.
\end{proof}

As in the previous cases, $\Psi^{\sst\Om}\ci\bar
\Pi^{\Th,\Om}_{\fObj_\A}:\SFai(\fObj_\A)$ or
$\oSFai(\fObj_\A,\Up,\La^\ci)\ra C^\ind(\A,\Om,\ha\bar\chi)$ are
also Lie algebra morphisms. We explain Theorem \ref{ai6thm7} using
the following example. For simplicity, let $X,Y\in\A$ be
indecomposable, with $[X]=\al$ and $[Y]=\be$ in $C(\A)$, and form
$\bde_{[X]},\bde_{[Y]}$ in $\oSFai(\fObj_\A,\Th,\Om)$. Then
$\Psi^{\sst\Om}(\bde_{[X]})=c^\al$ and $\Psi^{\sst\Om}(\bde_{[Y]})
=c^\be$. Consider the commutator~$[\bde_{[X]},\bde_{[Y]}]$.

We find $\Psi^{\sst\Om}(\bde_{[X]}*\bde_{[Y]})\!=\!
\bigl(\dim\Ext^1(Y,X)\!-\!\dim\Hom(Y,X)\bigr)c^{\al+\be}\!+
\!c_{[\{1,2\},\ka]}$, where $\ka(1)=\al$ and $\ka(2)=\be$. This is
because $\bde_{[X]}*\bde_{[Y]}$ is essentially the characteristic
function of all $Z$ in short exact sequences $0\ra X\ra Z\ra Y\ra
0$. The effect of applying $\Psi^{\sst\Om}$ is to `count' such
sequences in a special way. The nontrivial extensions are
parametrized by $P(\Ext^1(Y,X))$ and contribute $\dim\Ext^1(Y,X)$ to
the `number' of such $Z$, and the trivial extension $0\ra X\ra X\op
Y\ra Y\ra 0$ has stabilizer group
$\bigl(\Aut(X)\t\Aut(Y)\bigr)\lt\Hom(Y,X)$ and contributes
$\dim\Hom(Y,X)$ to the number of `virtual indecomposables'
multiplying $c^{\al+\be}$, and 1 to the number of `virtual
decomposables' multiplying~$c_{[\{1,2\},\ka]}$.

Exchanging $X,Y$ and using \eq{ai6eq32} gives $\Psi^{\sst\Om}
\bigl([\bde_{[X]},\bde_{[Y]}]\bigr)=\bar\chi([X],[Y])c^{\al+\be}$.
By \eq{ai6eq29}, this is $[c^\al,c^\be]$ in the Lie algebra
$C^\ind(\A,\Om,\ha\bar\chi)$, so $\Psi^{\sst\Om}\bigl([\bde_{[X]},
\bde_{[Y]}]\bigr)=\bigl[\Psi^{\sst\Om}(\bde_{[X]}),\Psi^{\sst\Om}
(\bde_{[Y]})\bigr]$, as we want. Thus we see that Theorem
\ref{ai6thm7} relies on \eq{ai6eq32} and a very particular way of
`counting' stack functions on $\fObj_\A$, such that the `number' of
extensions of $Y$ by $X$ is~$\dim\Ext^1(Y,X)-\dim\Hom(Y,X)$.

\medskip

\noindent{\small\sc The Mathematical Institute, 24-29 St. Giles,
Oxford, OX1 3LB, U.K.}

\noindent{\small\sc E-mail: \tt joyce@maths.ox.ac.uk}

\end{document}